\documentclass[12pt,letterpaper]{article}

\pdfoutput=1

\usepackage{graphics,epsfig,graphicx,color}
\graphicspath{{/EPSF/}{../Figures/}{Figures/}}
\usepackage[caption = false]{subfig}
\usepackage{float}
\usepackage{capt-of}

\usepackage{amssymb,latexsym}

\usepackage{setspace}

\usepackage{amsmath}

\usepackage{mathrsfs,amsfonts}

\usepackage{amsthm,amsxtra}

\usepackage[final]{showkeys}

\usepackage{stmaryrd}

\usepackage{algorithm}
\usepackage{algorithmicx}
\usepackage{algpseudocode}

\usepackage[openbib]{currvita}

\newif\ifPDF
\ifx\pdfoutput\undefined
	\PDFfalse
\else
	\ifnum\pdfoutput > 0
		\PDFtrue
	\else
		\PDFfalse
	\fi
\fi

\ifPDF
	\usepackage{pdftricks}
	\begin{psinputs}
		\usepackage{pstricks}
		\usepackage{pstcol}
		\usepackage{pst-plot}
		\usepackage{pst-tree}
		\usepackage{pst-eps}
		\usepackage{multido}
		\usepackage{pst-node}
		\usepackage{pst-eps}
	\end{psinputs}
\else
	\usepackage{pstricks}
\fi

\ifPDF
	\usepackage[debug,pdftex,colorlinks=true, 
	linkcolor=blue, bookmarksopen=false,
	plainpages=false,pdfpagelabels]{hyperref}
\else
	\usepackage[dvips]{hyperref}
\fi

\usepackage{fancyhdr}

\newcommand{\eps}{\varepsilon}

\newcommand{\fT}{\mathfrak{T}}

\newcommand{\bzero}{{\mathbf 0}}
\newcommand{\btheta}{{\boldsymbol\theta}}

 \newcommand{\bff}{\mathbf f}

 \newcommand{\bn}{\mathbf n}
 
\newcommand{\bq}{\mathbf q}

\newcommand{\bw}{\mathbf w} \newcommand{\bx}{\mathbf x}

\newcommand{\bE}{\mathbf E} 
 \newcommand{\bH}{\mathbf H}

 \newcommand{\bR}{\mathbf R}
 \newcommand{\bT}{\mathbf T}
 
\newcommand{\bW}{\mathbf W}

\newcommand{\cE}{\mathcal E}

 \newcommand{\cL}{\mathcal L}
\newcommand{\cM}{\mathcal M} \newcommand{\cN}{\mathcal N}
 \newcommand{\cP}{\mathcal P} 
 
 \newcommand{\cT}{\mathcal T}

\newenvironment{keywords}
{\noindent{\bf Key words.}\small}{\par\vspace{1ex}}
\newenvironment{AMS}
{\noindent{\bf AMS subject classifications 2010.}\small}{\par}

\newcommand{\red}[1]{\textcolor{red}{#1}}
\newcommand{\blue}[1]{\textcolor{blue}{#1}}

\title{Numerical algorithms based on Galerkin methods for the modeling of reactive interfaces in photoelectrochemical (PEC) solar cells}

\author{
	Michael Harmon\thanks{
		Institute for Computational Engineering and Sciences (ICES), 
		University of Texas, Austin, TX 78712;
		\href{mailto:mharmon@ices.utexas.edu}{mharmon@ices.utexas.edu}
	}
	\and
	Irene M. Gamba\thanks{
		Department of Mathematics and ICES,
		University of Texas, Austin, TX 78712;
		\href{mailto:gamba@math.utexas.edu}{gamba@math.utexas.edu}
	}
	\and
	Kui Ren\thanks{
		Department of Mathematics and ICES,
		University of Texas, Austin, TX 78712;
		\href{mailto:ren@math.utexas.edu}{ren@math.utexas.edu}
	}
}

\begin{document}

\maketitle



\begin{abstract}

This work concerns  the numerical solution of a coupled system of self-consistent reaction-drift-diffusion-Poisson equations that describes the macroscopic dynamics of charge transport in photoelectrochemical (PEC) solar cells with reactive semiconductor and electrolyte interfaces. We present three numerical algorithms, mainly based on a mixed finite element and a local discontinuous Galerkin method for spatial discretization, with carefully chosen numerical fluxes, and implicit-explicit time stepping techniques, for solving the time-dependent nonlinear systems of partial differential equations. We perform computational simulations under various model parameters to demonstrate the performance of the proposed numerical algorithms as well as the impact of these parameters on the solution to the model.
\end{abstract}






\section{Introduction}
\label{SEC:Intro}

The phenomena of reaction and transport of charged particles in spatial regions with semiconductor and electrolyte interfaces has been extensively investigated in the past few decades~\cite{BeFaPeWi-EMAC03,GaGeMa-JCP00,GaMa-JCP00,KaTvBaRa-CR10,Lewis-ACR90,Memming-Book15,NoMe-JPC96,PeFaPl-SEMSC08,PeFaWiBe-JPP04,SiGuMuAlFi-JAP09}; see~\cite{Lewis-JPCB98} for recent reviews on the subject. There are three main components in any such semiconductor-electrolyte interface problem: the physical processes inside the semiconductor, the physical processes inside the electrolyte and the mechanisms of charge transfer across the interface. While the physical mechanisms and mathematical modeling of charge generation, recombination, and transport in both semiconductors and electrolytes are well-understood~\cite{AnAlRi-Book03,BaThAj-PRE04,BiFuBa-RJE12,Eisenberg-SM00,Fawcett-Book04,Galler-Book05,Grasser-Book03,HoLiLiEi-JPCB12,Jerome-Book96,Jungel-Book09,KiBaAj-PRE07B,Liu-JDE09,LuZh-BJ11,MaPeAg-JCP88,MaRiSc-Book90,MaMu-IJHMT09,PaJe-SIAM97,Schroeder-Book94,ScNaEi-PRE01,Selberherr-Book01,SiNo-SIAM09} by now, the physics of the charge reaction and transfer processes on the semiconductor-electrolyte interface is still under intensive studies~\cite{GaGeMa-JCP00,GaMa-JCP00,Lewis-ACR90,Memming-Book15,NoMe-JPC96}, and the mathematical modeling and numerical simulation of these processes is far less studied~\cite{HeGaLeRe-SIAM15,SiGuMuAlFi-JAP09}.

In recent years, the fast development in the field of photoelectrochemical (PEC) solar cell design has greatly boosted research interests in charge transfer across semiconductor-electrolyte interfaces. There are many different versions of PEC cells with very different specifics, but their working mechanism are all very similar~~\cite{Gratzel-Nature01,Gratzel-JPPC03,KaTvBaRa-CR10}. When sunlight shines on the semiconductor, photons are absorbed. Electron-hole pairs are then generated. An applied electric field splits the electron-hole pair which causes a build up of charges at the interface. Electrons are then transferred across the interface from the semiconductor to the electrolyte or vice versa. The transfer of electrons across the interface induces the generation and elimination of reductant-oxidant (redox) pairs in the electrolyte.  The transport of charged particles within each domain and their transfer across the interface leads to a continuous current throughout the device. There are many possible semiconductor-liquid combinations; see for instance, ~\cite{FaLe-JPCB97} for ${\rm Si/viologen}^{2+/+}$ junctions, ~\cite{PoLe-JPCB97} for n-type ${\rm InP/Me_2Fc}^{+/0}$ junctions, and ~\cite[Tab. 1]{FoPrFeMa-EES12} for a summary of many other possibilities. 

Computational simulations of PEC cells have been used extensively by researchers to study macroscopic behaviors, such as current-voltage characteristics, of PEC cells. Accurate simulations of PEC cells requires  not only good understanding of the complex electrochemical processes involved in the cells, but also accurate mathematical models that can describe the macroscopic effects of these electrochemical processes. There have been great efforts in constructing such mathematical models for different types of PEC cells. In most of the previous studies~\cite{FoPrFeMa-EES12, KaAtLe-JAP05, SiGuMuAlFi-JAP09}, the complicated interface processes are replaced with empirical approximations, such as the ``Shottky approximation''. In~\cite{HeGaLeRe-SIAM15}, a new mathematical model is proposed to more accurately model the interface charge transport process. Numerical simulations showed that the new model can produce more accurate simulations in regimes where classical approximations are not valid.
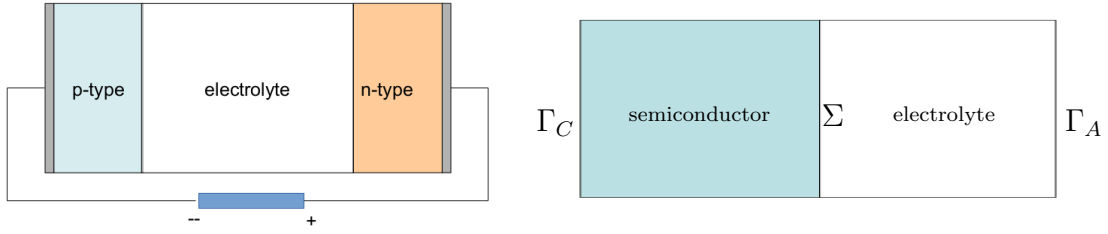
\begin{figure}[!ht]
\centering
\includegraphics[angle=0,width=0.4\textwidth]{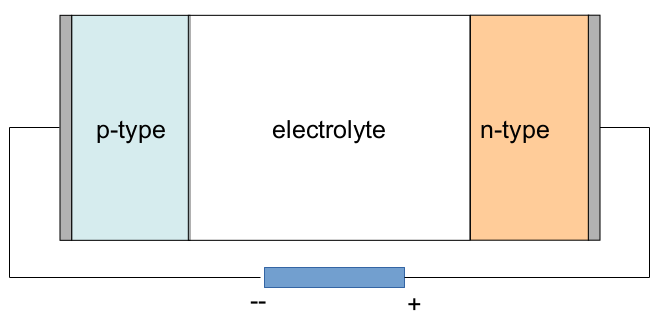}\hskip 1cm
\raisebox{0.15\height}{\includegraphics[angle=0,width=0.4\textwidth,height=0.15\textwidth]{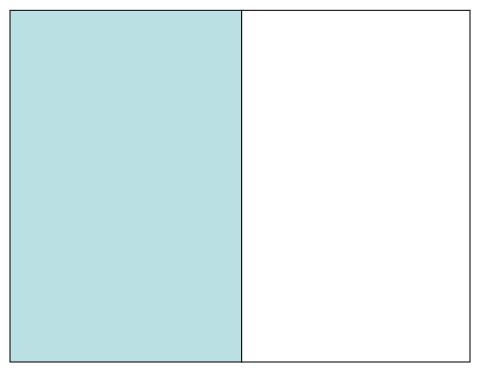}}
\put(-92,39){$\Sigma$}
\put(-165,41){{\scriptsize\bf \rm semiconductor}}
\put(-65,41){{\scriptsize \bf \rm electrolyte}}
\put(-200,37){$\Gamma_C$}
\put(0,37){$\Gamma_A$}
\caption{Left: Two-dimensional sketch of a typical semiconductor-electrolyte PEC cell. Right: The half-cell geometry, the semiconductor and electrolyte components are separated by the interface $\Sigma$.}
\label{FIG:PEC Cell}
\end{figure}

In this paper we propose numerical algorithms for the simulation of semiconductor-electrolyte PEC cells based on the mathematical model developed in~\cite{HeGaLeRe-SIAM15}. To simplify the presentation, we will only consider the half-cell setup which has three components: the semiconductor, the electrolyte, and the interface between the two; see the sketch in Fig.~\ref{FIG:PEC Cell}. To simulate the full cell, we only need to apply the same procedure to the second semiconductor-electrolyte interface at the other electrode.

The rest of this work is structured as follows. We first review briefly in Section~\ref{SEC:Model} the mathematical model for the dynamics of charge transport in PEC cells. We then propose in Section~\ref{SEC:MFEM-DG} a coupled discretization scheme based on a mixed finite element method and a local discontinuous Galerkin method for the spatial component of our mathematical model. In Section~\ref{SEC:IMEX} we combine our spatial discretization methods with some implicit-explicit time-stepping strategies to complete our numerical schemes for the system. We then present in Section~\ref{SEC:Num} numerical simulations under various conditions to demonstrate the performance of the algorithms we proposed. Concluding remarks are offered in Section~\ref{SEC:Concl}.

\section{The mathematical model}
\label{SEC:Model}

We review in this section the mathematical model we constructed for the simulation of semiconductor-electrolyte interfaces in photoelectrochemical cells in the half cell setup depicted in Fig.\ref{FIG:PEC Cell}. The model consists of three main components: a system of equations for transport of electrons and holes in the semiconductor $\Omega_{S}$, a system of equations for charge transport in the electrolyte $\Omega_{E}$ and a set of interface conditions on the semiconductor-electrolyte interface $\Sigma$. We will present the model directly in nondimensionalized form so that we do not need to keep repeating the equations. The characteristic scales used in the nondimensionalization process are (see, for instance, ~\cite{KiBe-Book91,LaBa-JES76A,LaBa-JES76B,Nelson-Book03,PeFaPl-SEMSC08}),
\begin{equation*}\label{EQ:Characteristic}
l^*=10^{-4}\ {\rm (cm)},\ \ \  t^*=10^{-12}\ {\rm (s)},\ \ \ \Phi^*=25.85\ {\rm (mV)},\ \ \ C^*= 10^{16}\ {\rm (cm^{-3})},
\end{equation*}
where $l^*$ is the characteristic length scale of the device, $t^*$ is the characteristic time scale, $\Phi^*$ is the characteristic voltage, and $C^*$ is the characteristic charge density. The unit for each quantity is given in the corresponding parenthesis. We refer interested reader to~\cite{HeGaLeRe-SIAM15} for detailed discussions on the model and its nondimensionalization.

\paragraph{The semiconductor equations.} In the semiconductor component the transport of electrons and holes is described by a standard system of drift-diffusion-Poisson equations. Let $\rho_n$ and $\rho_p$ denote the (nondimensionalized) densities of the electrons and holes respectively, and $\Phi$ denote the electric potential field. Then we have,
\begin{align}
\frac{\partial \rho_{n}}{\partial t} \, + \, 
\nabla  \cdot \,\mu_n\left( \, - \alpha_n \rho_{n} \nabla \Phi  
\, - \,  \nabla \, \rho_{n}  \, \right) 
\; &=- \; 
R(\rho_{n}, \rho_{p})
\, + \, 
\gamma \,  G(\bx) , && \text{in} \ (0,T] \times \Omega_{S}, \nonumber \\
\frac{\partial \rho_{p}}{\partial t} 
\, + \, 
\nabla \cdot \, \mu_p \left( -\alpha_p \rho_p \nabla \Phi 
\, - \, 
\nabla \, \rho_{p} \, \right) 
\; &= -\;
R(\rho_{n}, \rho_{p})
\, + \, 
\gamma G(\bx) ,  && \text{in} \ (0,T] \times \Omega_{S} \label{EQ:DDP S},\\
-\nabla \cdot \left( \, \lambda_S^2 \, \nabla \Phi \right)
\; &= \cN(\bx) +\alpha_p\rho_p+\alpha_n\rho_n,  && \text{in} \ (0,T]\times \Omega_{S}.  \nonumber
\end{align}
Here the coefficients $\mu_n$ and $\mu_p$ are the mobilities of the electrons and holes, respectively. The Einstein relation between the mobilities and the diffusivity has been used. The coefficients $\alpha_n=-1$ and $\alpha_p=1$ are respectively the charge numbers of the electrons and the holes. $\lambda_S=\dfrac{1}{l^*}\sqrt{\frac{\Phi^*\eps_S}{qC^*}}$ is the (rescaled) Debye length inside the semiconductor with $\eps_S$ being the permittivity of the semiconductor and $q$ being the charge of an electron. The doping contribution consists of a component $\cN_{D}$ from the donor and a component $\cN_{A}$ from the acceptor and $\cN=\cN_D-\cN_A$. The function $R(\rho_n,\rho_p)$ describes the rate at which electron and hole pairs are eliminated through recombination (when it is positive) and generated (when it is negative) due to thermal excitation. In this paper, we use Shockley-Reed-Hall recombination~\cite{ShRe-PhysRev1952} function which is believed to be the most dominant recombination mechanism in real materials~\cite{Nelson-Book03}. The function reads,
\begin{equation*}
R(\rho_{n}, \rho_{p})
\; = \;
\frac{\rho_{n}\rho_{p}-\rho_{i}^{2}}{\tau_{n} \,  (\rho_{n} \, + \, \rho_{i}) \, + \, \tau_{p} \, ( \rho_{p} \, + \, \rho_{i})},
\end{equation*}
where $\rho_{i}$ denotes the intrinsic electron density (again rescaled), while $\tau_{n}$ and $\tau_{p}$ are the rescaled electron and hole lifetimes. 

The generation of electrons and holes due to excitation of sunlight is modeled using the macroscopic source function $G(\bx)$. The parameter $\gamma \in \{0,1\}$ is used here to control the on ($\gamma=1$) and off ($\gamma=0$) condition of the sunlight illumination. The exact form of the function $G$ depends on the illumination geometry. In most of the applications, it is assumed that photons from sun light travel inside the device in straight lines (i.e., without being scattered). Therefore, $G$ takes the following form:
\begin{equation*}
G(\bx) = \sigma_a (\bx) G(\bx_{0}) \exp\big(-\int_{0}^{s} \sigma_a(\bx_{0}+s'\btheta_{0}) ds'\big), \qquad \bx = \bx_{0} + s \btheta_{0}
\end{equation*} 
where the point $\bx_{0}$ is the photon's incident location and $\btheta_{0}$ is the incident direction. The function $G(\bx_0)$ represents the surface photon flux at the point $\bx_0$ and the function $\sigma_a(\bx)$ is the absorption coefficient at $\bx$.

On the portion of the semiconductor boundary where an Ohmic metal contact is located, $\Gamma_C$, the charge densities take on their equilibrium values, while the potential is the sum of the applied voltage $\Phi_{\text{app}} $ and the so-called ``built-in'' potential $\Phi_{\text{bi}}$. That is,
\begin{align}\label{EQ:BC Ohmic}
\rho_{n} \vert_{\Gamma_C}  \; = \;  \rho_{n}^{e} \, , 
&&
\rho_{p} \vert_{\Gamma_C} \; = \; \rho_{p}^{e}  \, , 
&& \Phi \vert_{\Gamma_C} \;  = \; \Phi_{\text{bi}} \, + \, \Phi_{\text{app}}\, , && \text{on} \ (0,T] \times \Gamma_C.
\end{align}
We assume that the rest of the semiconductor boundary, $\Gamma_N^s\equiv\partial\Omega_S\backslash(\Sigma\cup\Gamma_C)$ is insulated. Therefore, we have the following conditions,
\begin{equation}\label{EQ:BC Neumman S}
\bn\cdot \mu_n(-\alpha_n\rho_n\nabla\Phi-\nabla\rho_n)=0, \quad 
\bn\cdot \mu_p(-\alpha_p\rho_p\nabla\Phi-\nabla\rho_p)=0, \quad \bn\cdot\nabla\Phi=0,\quad \mbox{on}\ (0, T]\times\Gamma_N^s,
\end{equation}
where $\bn$ is the unit outward normal vector to the semiconductor domain.

\paragraph{The electrolyte equations.} In the electrolyte component the transport of reductants and oxidants is described by a similar system of reaction-diffusion-Poisson equations. Let $\rho_{r}$  and $\rho_o$ be the (nondimensionalized) densities of the reductants and oxidants respectively, then we have,
\begin{align}
\frac{\partial \rho_{r}}{\partial t}
\, + \, 
\nabla  \cdot \, \mu_r \left(-\alpha_r \rho_r \nabla \Phi \, 
\, - \, 
\nabla \, \rho_{r}  \, \right) 
\; &= \; 0 , && \text{in} \ (0,T] \times \Omega_{E},  \nonumber \\
\frac{\partial \rho_{o}}{\partial t} 
\, + \, 
\nabla \cdot \, \mu_o \left( -  \, \alpha_o \rho_o \nabla \Phi 
\, - \, 
\nabla \, \rho_{o} \, \right) 
\; &= \; 0, && \text{in} \ (0,T] \times \Omega_{E}  \label{EQ:DDP E}, \\
- \nabla \cdot \left( \, \lambda_E^2 \, \nabla \Phi \right)
\; &= \alpha_o \rho_o +\alpha_r \rho_r, && \text{in} \ (0,T]\times\Omega_{E}.   \nonumber
\end{align}
Here $\mu_{r}$ and $\mu_{o}$ are respectively the reductant and oxidant mobilities, and $\lambda_E=\dfrac{1}{l^*}\sqrt{\frac{\Phi^*\eps_E}{qC^*}}$ is the rescaled Debye length inside the electrolyte with $\eps_E$ being the permittivity of the electrolyte. The charge carriers numbers $\alpha_r$ and $\alpha_o$ for the reductants and oxidants are their oxidation numbers respectively.  For our model here to conserve the total charges in the system, we have to have $ \alpha_o - \alpha_r =1$. We have assumed that besides the local imbalance of reductants and oxidants, the electrolyte is charge neutral. Therefore, we do not have a fixed background charge density (i.e.  doping profile $\cN$) in the Poisson equation for the electrolyte. Moreover, our model does not include any generation or recombination mechanisms in the electrolyte domain since we are only considering so-called ``heterogeneous reactions". That is to say that chemical reactions can only occur at the interface and not within the bulk of the electrolyte. However, we can indeed generalize a bit in that direction; see, for instance, the more complete model in~\cite{HeGaLeRe-SIAM15}.

On the electrolyte boundary $\Gamma_A$, the location of the anode, the reductant and oxidant densities as well as the electric potential take their bulk values~\cite{BaFa-Book00}. That is,
\begin{align}
\rho_{r} \vert_{\Gamma_A}  \; = \;  \rho_{r}^{\infty},
&&
\rho_{o} \vert_{\Gamma_A} \; = \; \rho_{o}^{\infty},
&&
\Phi \vert_{\Gamma_A} \; = \; \Phi^{\infty}, 
&& \text{on} \ (0,T]\times\Gamma_A
\end{align}
We assume that the rest of the electrolyte boundary, $\Gamma_N^e\equiv\partial\Omega_E\backslash(\Sigma\cup\Gamma_A)$ is insulated. Therefore, we have the following conditions,
\begin{equation}\label{EQ:BC Neumman E}
\bn\cdot \mu_r(-\alpha_r\rho_r\nabla\Phi-\nabla\rho_r)=0, \quad 
\bn\cdot \mu_o(-\alpha_o\rho_o\nabla\Phi-\nabla\rho_o)=0, \quad \bn\cdot\nabla\Phi=0,\quad \mbox{on}\ (0, T]\times\Gamma_N^e,
\end{equation}
where again $\bn(\bx)$ is the unit outward normal vector to the electrolye domain.

\paragraph{The interface conditions.} The third component of the mathematical model is the conditions for the unknown variables on the semiconductor-electrolyte interface. For the electric potential, the usual continuities of the potential and the electric flux are imposed as follows,
\begin{align}\label{EQ:Int Phi}
\Phi|_{\Sigma^S} - \Phi|_{\Sigma^E} \, = \, 0, &&  \bn_{\Sigma^S} \cdot \lambda_S^2 \nabla \, \Phi + \bn_{\Sigma^E} \cdot\lambda_E^2\nabla\Phi \, = \, 0, 
&& \text{on} \ (0,T]\times\Sigma,
\end{align}
where $\Sigma^S$ and $\Sigma^E$ are respectively the semiconductor and the electrolyte sides of $\Sigma$, while $\bn_{\Sigma^S}$ and $\bn_{\Sigma^E}$ are the unit normal vector of the semiconductor and the electrolyte domain on the interface. The conditions for the densities of electron-hole pairs and reductant-oxidant pairs on the interface are determined by the electron transfer dynamics across the interface. More precisely, we have the following Robin type of interface conditions,
\begin{align}
\label{EQ:Int Rho n}
\bn_{\Sigma^S} \cdot \mu_n\left(-\alpha_n \rho_n \nabla \Phi -\nabla \rho_n \right) &=I_{et}(\rho_{n}-\rho_{n}^{e},\rho_{o}), & \mbox{on}\ (0,T]\times\Sigma,\\ 
\label{EQ:Int Rho p}
\bn_{\Sigma^S} \cdot \mu_p\left(-\alpha_p \rho_p \nabla \Phi -\nabla \rho_p \right) &=I_{ht}(\rho_p-\rho_p^e, \rho_{r}), & \mbox{on}\ (0,T]\times\Sigma,\\ 
\label{EQ:Int Rho r}
\bn_{\Sigma^E} \cdot \mu_r\left(-\alpha_r \rho_r \nabla \Phi - \nabla \rho_r \right)&=I_{ht}(\rho_{p}-\rho_{p}^{e}, \rho_{r}) - I_{et} (\rho_{n} - \rho_{n}^{e}, \rho_{o}), & \mbox{on}\ (0,T]\times\Sigma,\\
\label{EQ:Int Rho o}
\bn_{\Sigma^E} \cdot \mu_o \left(-\alpha_o \rho_o \nabla \Phi -\nabla \rho_o \right) &=-I_{ht}(\rho_{p}- \rho_{p}^{e}, \rho_{r}) + I_{et}(\rho_{n}-\rho_{n}^{e}, \rho_{o}), & \mbox{on}\ (0,T]\times\Sigma,
\end{align}
where the functions $I_{et}(x,y)=k_{et} x y$, $I_{ht}(x,y)=k_{ht} x y$, and the constants $k_{et}$ and $k_{ht}$ are the rescaled electron and hole transfer rates respectively.  

The current density for each carrier is defined as the total flux of that carrier,
\begin{align}
\textbf{J}_{i} &= \mu_{i}\left( -\alpha_{i} \rho_{i} \nabla \Phi -\nabla \rho_i \right), & \text{for} \  i = \{n,p,r,o\}.
\end{align}
The total current density through the device will then be defined as,
\begin{equation}
\textbf{J}(\textbf{x}) = \left\{
\begin{array}{cc}
\alpha_{n} \textbf{J}_{n} + \alpha_{p} \textbf{J}_{p}, & \textbf{x}
\in \Omega_{S} \\
\alpha_{r} \textbf{J}_{r} + \alpha_{o} \textbf{J}_{o}, & \textbf{x}
\in \Omega_{E}
\end{array}  \right. \label{EQ:Current}
\end{equation}
Now using the fact $\alpha_{r} - \alpha_{o} = 1$, then substituting~\eqref{EQ:Int Rho n}, ~\eqref{EQ:Int Rho p}, ~\eqref{EQ:Int Rho r} and~\eqref{EQ:Int Rho o} into~\eqref{EQ:Current}  we can see that we have a continuous current through the interface $\Sigma$.

The system of equations in~\eqref{EQ:DDP S} and~\eqref{EQ:DDP E}, together with the interface conditions~\eqref{EQ:Int Phi}, ~\eqref{EQ:Int Rho n}-~\eqref{EQ:Int Rho o} and the boundary conditions can now be complemented with an initial condition to get a complete mathematical model for the charge transport process in semiconductor-electrolyte system. More complicated effects can be included in this model. However, all key ingredients are included here. Moreover, the current model is more general than the commonly-used Butler-Volmer model~\cite{BaFa-Book00,NeTh-Book04} in the electrochemistry literature. Indeed, we can derive the Butler-Volmer model from our model in the setting where the Butler-Volmer model is valid. We refer interested reader to~\cite{HeGaLeRe-SIAM15} for detailed discussions on this model and some preliminary simulations of one-dimensional devices based this model. 

The initial-boundary value problem associated to system ~\eqref{EQ:BC Ohmic}-\eqref{EQ:Int Rho o}, modeling an interface reaction across an interior surface whose transfer condition is determined by the non-linear by the system~\eqref{EQ:Int Rho n}-\eqref{EQ:Int Rho o} on such surface interface, is a strictly parabolic problem. These type of parabolic problems with  such  reaction terms are expected to have a comparison principle in the sense that any two solutions $\rho_i^1(t,x) $ and $ \rho_i^2(t,x), i=n,p,o,r$, satisfying $\rho_i^1(t,x) \leq \rho_i^2(t,x)$ in the set $\Omega_T:=[0,T)\times \overline{(\Omega_S\cup\Omega_E)}$, can not have a contact point $(t_0,x_0)$ where $\rho_i^1(t, x) = \rho_i^2(t, x)$ in the interior of the set $\Omega_T$ since the relations \eqref{EQ:Int Rho n}-\eqref{EQ:Int Rho o}  on the normal derivatives at the transmission surface $\Sigma$ yield monotone relations.  In fact, for a given external electric field,  these conditions are expected to be sufficient for an existence and uniqueness  theory in a variational formulation, or  by means of the classical continuity method. Details of these results will be given in a forthcoming  paper. In addition, such  results on existence and uniqueness would yield  error estimates for the proposed DG Method generated by the implicit-explicit-explicit alternating Schwarz (AS-IMEXEX ) scheme algorithm to be presented in Section~\ref{SEC:IMEX}.

However, the main challenges in constructing a numerical algorithm that produces a reliable  approximating numerical solution  to system~\eqref{EQ:DDP S}-\eqref{EQ:DDP E} and ~\eqref{EQ:Int Rho n}-\eqref{EQ:Int Rho o} are due to the nonlinearities arising, not only, from the coupling of the corresponding drift diffusion (transport) equations for the pairs $\rho_n,\rho_p$  and  $\rho_o,\rho_r$ to their corresponding  Poisson equations for charges  in both the semiconductor region $\Omega_S$ and the electrolyte region $\Omega_E$, respectively, but also, with a nonlinear interface condition on the common interface set $\Sigma=\overline{\Omega}_S \cap \overline{\Omega}_E$ as well as recombination-generation mechanisms added to  the corresponding transport equations.

While there have been many contributions in the last forty years on the approximations to  solution of a single  drift-diffusion-Poisson system in the semiconductor region, there is no previous work, to the best of our knowledge, that deals with the whole coupled system that incorporates the nonlinear interface conditions such as~\eqref{EQ:Int Rho n}-\eqref{EQ:Int Rho o}.  In addition the evolution problem under consideration is effectively multi-scale  in the sense that the evolution of the system in the semiconductor and the corresponding one in the electrolyte evolve at different time scales due to the quantitative scaling differences in the relevant  physical parameters such as mobilities or characteristic charge densities.

In particular, in order to ensure computational stability, the size of the time step (chosen for any time discretized scheme) is dominated by the value in the fast-varying component, that is, the semiconductor component. Using small time steps in the numerical simulations for the problem in the whole domain results in regions of stiffness caused by boundary layer formation where sharp transitions in  densities and  electric potential occur (i.e. near interfaces). This problem can be tackled by using fine enough meshes around the interface that would resolve those boundary layers. Such considerations limit even further  the time step sizes on the  employed computational method.

\section{Mixed schemes for spatial discretization}
\label{SEC:MFEM-DG}

We consider in this section the discretization of the system of nonlinear drift-diffusion-Poisson equations in the spatial variable. There are many different discretization strategies for such systems in the literature; see, for instance, ~\cite{BaRoFi-SIAM83,BaHe-JCP13,BoMiSa-JCP00,BrMaPi-SIAM89,BrFeMaWo-M3AS13,CadeGiSa-CMAME06,Chen-SIAM95,ChCo-NM95,deJeSa-JCP09,DeEl-JCP02,DeGaMe-JCP07,DoRu-SIAM82,DoGaSq-MAC86,GaSq-SIAM89,JeKe-SIAM91,JiPa-JCP00,Kulikovsky-JCP95,LeGrPe-JCP10,LuHoMcZh-JCP10,MaRiSc-Book90,MaSiShAtIs-JCP14,MoWe-JCP15,Selberherr-Book01,SpGuAlMaBy-JCP15,StRoNi-JCP13,ThNa-JCP98} and references therein. In this work, we take a hybrid discretization strategy for our problem. To be precise, we discretize the drift-diffusion components of the system using a local discontinuous Galerkin (LDG) method~\cite{BrFeMaWo-M3AS13,CoDa-MAFELAP00,CoSh-SIAM98,CoSh-JSC02,WaShZh-SIAM15,XuSh-CiCP10} and discretize the Poisson equation in the system with a mixed finite element method (MFEM)~\cite{Babuska-NM73,BoBrFo-Book13,Brezzi-RAIRO74,BrFo-Book91,Gatica-Book14,GlWh-DDM87,RaTh-MAFEM77}. Such a hybrid discretization strategy has also been employed by many other authors in solving coupled reactive transport problems~\cite{ChRi-AMAME09,DaSuWh-CMAME04,GiSuWhYo-SIAM08,VaYo-SIAM09}.

In the rest of the paper, we denote by $\left(\cdot , \cdot \right)_{\Omega}$ and $\langle \cdot , \cdot \rangle_{\partial\Omega}$ the $L^{2}$ inner products over the domain $\Omega$  and its boundary $\partial\Omega$ respectively. We denote by $\cT_h$ a finite element triangulation of the domain $\Omega$, with $\cE_h$ the set of its faces (or edges in two dimension). We assume that the partition is non-degenerate and quasi-uniform. Moreover, we assume that no element of the partition intersects the interface, and at most one face of an element can lie along the interface. We denote by $\cE_h^i$ the union of all interior faces in the partition that do not intersect the interface or boundaries. Let $K^+\in\cT_h$ and $K^-\in\cT_h$ be two adjacent elements that have a common face $\mathfrak e\in\cE_h$, we denote by $\bn^+$ and $\bn^-$ respectively their unit outward normal vectors on $\mathfrak e$, and we denote by $f^+$ and $f^-$ respectively traces (when they exist) of the function $f$ on $\mathfrak e$ from the interiors of the elements. We define the usual average and jump of a function across a face of an element as,
\[
	\{f\}=\dfrac{1}{2}(f^-+f^+) \qquad \mbox{and}\qquad \llbracket f \rrbracket = f^+\bn^++f^-\bn^-,
\]
where $f$ is a scalar function, and
\[
	\{\bff\}=\dfrac{1}{2}(\bff^-+\bff^+) \qquad \mbox{and}\qquad  \llbracket \bff \rrbracket = \bff^+\cdot\bn^+ + \bff^-\cdot\bn^-,
\]
when $\bff$ is vector-valued. 

\subsection{Mixed finite element discretization of Poisson}

To use the mixed finite element method on the Poisson equation, we first introduce an auxiliary variable $\bE$ and rewrite equation in the following mixed form: for any $t\in(0, T]$,
\begin{equation}\label{EQ:Poisson FO}
\begin{array}{rcll}
\lambda^{-2} \bE + \nabla \Phi & =& 0, & \mbox{in}\ (0,T]\times\Omega,\\
\nabla \cdot \bE &=& f, & \mbox{in}\ (0,T]\times \Omega,\\
 \bn_{\Sigma^S}\cdot \bE|_{\Sigma^S}+\bn_{\Sigma^E}\cdot \bE|_{\Sigma^E} =\bzero, & &\Phi|_{\Sigma^S}-\Phi|_{\Sigma^E} = 0, & \mbox{on}\  (0,T]\times\Sigma,\\
\Phi &=& \Phi_O, &  \mbox{on}\  (0,T]\times(\Gamma_C\cup \Gamma_A),\\
\bn\cdot\bE &=&0, & \mbox{on}\  (0,T]\times(\Gamma_N^s\cup \Gamma_N^e),
\end{array}
\end{equation}
where we have used the notations,
\[
	\lambda(\bx)=\left\{
		\begin{array}{cc}
		\lambda_S, &\bx\in\Omega_S\\
		\lambda_E, &\bx\in\Omega_E
		\end{array}\right.,
	\qquad 
	\Phi_O(\bx)=\left\{
		\begin{array}{cc}
		\Phi_{\rm app}+\Phi_{\rm bi}, &\bx\in\Gamma_C\\
		\Phi^\infty, &\bx\in\Gamma_A
		\end{array}\right. ,
\]
\[
	\mbox{and}\ \ f(t,\bx) = \left\{
		\begin{array}{cc}
		\cN(\bx)+\alpha_p\rho_p+\alpha_n\rho_n,&  \bx\in\Omega_S, \\
		\alpha_o\rho_o + \alpha_p \rho_p,  & \bx\in\Omega_E.
		\end{array}\right. .
\]
Note that this formulation combines the Poisson equations in~\eqref{EQ:DDP S} and~\eqref{EQ:DDP E} together so that the new equation is posed in domain $\Omega=\Omega_S\cup\Omega_E$.

We employ a standard Raviart-Thomas mixed finite element discretization of $(\Phi,\bE)$. Let $W_{h,k}(\Omega)$ and $\bW_{h,k}(\Omega)$ be the continuous, piecewise polynomial approximating spaces
\[
	W_{h,k}(\Omega) = \left\{ w \in L^2(\Omega) : w|_{K} \in  \cP^{k}(K),\ \forall K \in \cT_{h} \right\},
\]
and
\[
	\bW_{h,k}(\Omega) = \left\{\bw\in \left(L^2(\Omega)\right)^d : \bw|_K \in \bR\bT_k(K),\  \forall K \in \cT_{h} \right\},
\]
where $\cP^k$ $(k\ge 1)$ denotes the space of complete polynomials of degree $k$ and $\bR\bT_k(K)$ $(k\ge 1)$ denotes the standard Raviart-Thomas element of degree $k$~\cite{RaTh-MAFEM77}. It is known that $\bW_{h,k}(\Omega)$ is a subspace of $H({\rm div},\Omega)$. Functions in $\bW_{h,k}$ have a continuous normal component across faces of elements. Let $\bW_{h,k}^N(\Omega)=\bW_{h,k}(\Omega) \cap\{\bw: \bn \cdot \bw=0,\ \mbox{on}\ \Gamma_N^s\cup\Gamma_N^e\}$. The Raviart-Thomas mixed finite element method seeks an approximation $(\Phi, \bE)\in C\left((0, T]; W_{h,k}\right)\times C\left((0, T]; \bW_{h,k}^{N}\right)$ that satisfies
\begin{equation}\label{EQ:MFEM I}
\left(\bw, \lambda^{-2} \bE \right)_{\Omega}-\left(\nabla\cdot\bw, \Phi \right)_{\Omega}=-\langle   \bn\cdot\bw, \Phi_{\rm app} + \Phi_{\rm bi} \rangle_{\Gamma_C} - \langle \bn\cdot\bw, \Phi^{\infty} \rangle_{\Gamma_A},\ \ \forall \bw\in\bW_{h,k}^{N}
\end{equation}
\begin{equation}\label{EQ:MFEM II}
\left(w, \nabla\cdot \bE \right)_{\Omega} = \left(w, f \right)_{\Omega}, \ \ \forall w \in W_{h,k}   
\end{equation} 
for any $t\in(0, T]$. Here we have assumed that the functions in the boundary conditions $\Phi_{\rm app}+\Phi_{\rm bi} \in L^2(\Gamma_C)$ and $\Phi^\infty \in L^2(\Gamma_A)$.

This Raviart-Thomas mixed finite element discretization leads to the so-called ``saddle point problem'' in the MFEM literature. It results in the following linear system for the degree of freedom vectors (which we denote by the same symbols),
\begin{equation}\label{EQ:Poisson MAT}
\left[ \begin{matrix}
A & B^{\fT} \\
B & 0 \\
\end{matrix} \right]
\, \left[ \begin{matrix}\
\bE(t,\cdot) \\
\Phi(t,\cdot) \\
\end{matrix}
\right] \; = \;
\left[ 
\begin{matrix}
\bH(t,\cdot) \\
F(t,\cdot)
\end{matrix} 
\right],\ \ \forall t\in(0, T],
\end{equation}
where $\bH$ and $F$ are vectors that incorporate the boundary conditions and right-hand side of Poisson's equation respectively. It has been verified that this saddle point problem with our choice of approximation spaces satisfies the Babuska-Brezzi \emph{inf-sup} conditions~\cite{Babuska-NM73,Brezzi-RAIRO74} and is therefore well-posed.

The system matrix in~\eqref{EQ:Poisson MAT} is symmetric but indefinite. More importantly, it is independent of time. Therefore, this matrix is assembled and factorized only once at the beginning of the time evolution. The right hand side vector, $[\bH, F]^\fT$, depends on time and has to be assembled at every time step using the current values of the charge densities.

\subsection{Local discontinuous Galerkin discretizations}
\label{SUBSEC:DG}

We use a local discontinuous Galerkin (LDG) method~\cite{BrFeMaWo-M3AS13,CoDa-MAFELAP00,CoSh-SIAM98,CoSh-JSC02,WaShZh-SIAM15,XuSh-CiCP10} to discretize the drift-diffusion equations for charge transport. To simplify the presentation of the discretization, we consider in this section only the drift-diffusion equation for the electrons in the system:
\begin{align}
\partial_{t} \rho_n + \nabla \cdot \mu_n\left(-\alpha_n \bE \rho_n -\nabla \rho_n \right)&=-R(\rho_n,\rho_p) +\gamma G(\bx), && \mbox{in} \ (0, T]\times\Omega_S, \\
\bn_{\Sigma^S} \cdot \mu_n \left(- \alpha_n\bE\rho_n-\nabla\rho_n\right)&=I_{et}(\rho_n-\rho_n^e,\rho_o), && \mbox{on} \ (0, T]\times\Sigma, \\
\bn \cdot \mu_n\left(-\alpha_n\bE\rho_n-\nabla\rho_n \right)&= 0, && \mbox{on}\ (0, T]\times\Gamma_{N}^s, \\
\rho_n &= \rho_n^e, && \mbox{on}\ (0, T]\times\Gamma_C,
\end{align}
where $\alpha_n$ is now the new notation for $\alpha_n/\lambda_S^2$. Equations for other charge carriers are discretized in the same way.

As in the mixed finite element formulation, the LDG discretization requires the drift-diffusion equation be written as a first-order system. We do this by introducing an auxiliary variable, which we call the current-flux variable,  $\bq_n$:
\begin{align}\label{EQ:DDP S FO}
\partial_{t} \rho_n + \nabla \cdot \bq_n &=-R(\rho_n,\rho_p) +\gamma G(\bx), && \mbox{in} \ (0, T]\times\Omega_S, \\
\mu_n^{-1}\bq_n &= -\alpha_n\bE\rho_n-\nabla\rho_n, && \mbox{in} \ (0, T]\times\Omega_S, \\
\bn_{\Sigma^S} \cdot \bq_n &=I_{et}(\rho_n-\rho_n^e,\rho_o), && \mbox{on} \ (0, T]\times\Sigma, \\
\bn \cdot \bq_n &= 0, && \mbox{on}\ (0, T]\times\Gamma_{N}^s, \\
\rho_n &= \rho_n^e, && \mbox{on}\ (0, T]\times\Gamma_C.
\end{align}

We define the function space $\widetilde\bW_{h,k}(\Omega_S)=\big(W_{h,k}(\Omega_S)\big)^d$. The LDG method consists of finding approximations $\rho_n\in C((0,T]; W_{h,k})$ and $\bq_n\in C((0, T]; \widetilde\bW_{h,k})$ such that, for all $t\in(0, T]$, for all $(w, \bw)\in W_{h,k}\times\widetilde\bW_{h,k}$, we have,
\begin{eqnarray}
\label{EQ:DDP LDG1}
\left(w, \partial_t \rho_n \right)_{\Omega_S}
-\left(\nabla w, \bq_n \right)_{\Omega_S} + \langle \llbracket w \rrbracket, \widehat{\bq}_n  \rangle_{\cE_{h}^{i,s}} +\langle w, \bn\cdot \widehat{\bq}_n  \rangle_{\partial\Omega_S} = -\left(w, R(\rho_n, \rho_p)-\gamma G\right)_{\Omega_S}, \\
\label{EQ:DDP LDG2}
\left(\bw, \mu_n^{-1} \bq_n \right)_{\Omega_S}-\left(\nabla \cdot \bw, \rho_n \right)_{\Omega_S}
+ \left(\bw, \alpha_n\bE\rho_n\right)_{\Omega_S}
+ \langle \llbracket \bw \rrbracket, \widehat{\rho}_n \rangle_{\cE_{h}^{i,s}}+\langle \bn\cdot \bw, \widehat{\rho}_n \rangle_{\partial\Omega_S} =0,
\end{eqnarray}
where $\cE_h^{i,s}$ refers to the interior faces in the semiconductor domain, excluding the interface $\Sigma$ which is part of $\partial\Omega_S$. The terms $\widehat{\rho}_n$ and $\widehat{\bq}_n$ are the numerical fluxes to be selected to ensure uniqueness, stability and accuracy of the solution to the resulting system. Boundary conditions are also imposed through the definition of the numerical fluxes. In our implementation, we take the flux $\widehat{\rho}_n$ on $\partial K$ as
\begin{equation} \label{EQ:U Flux}
\widehat{\rho}_n = \left\{
\begin{array}{cl}
\left\{\rho_n \right\} \ + \ \boldsymbol\beta \cdot \llbracket \rho_n \rrbracket, & \partial K\in \mathcal{E}_{h}^{i,s} \\
\rho_n, & \partial K\in \Sigma \\
\rho_n, & \partial K\in \Gamma_{N}^s \\
\rho_n^e, & \partial K \in \Gamma_C\\
\end{array}
\right.
\end{equation}
and the flux $\widehat{\bq}_n$ on $\partial K$ as,
\begin{equation} \label{EQ:Q Flux}
\widehat{\bq}_n  = \left\{
\begin{array}{cl}
\left\{ \bq_n \right\} \ - \ \llbracket \bq_n \rrbracket \boldsymbol \beta + \tau \llbracket \rho_n \rrbracket, & \partial K \in \cE_{h}^{i,s} \\
I_{et}(\rho_n-\rho_n^e,\rho_o)\bn_{\Sigma^S}, & \partial K \in \Sigma \\
\bzero, & \partial K \in \Gamma_{N}^s \\
\bq_n + \tau \left( \rho_n -\rho_n^{e} \right)\bn,  & \partial K \in \Gamma_C \\
\end{array}
\right.
\end{equation}
Here the term $\boldsymbol \beta$ is a constant unit vector which does not lie parallel to any element face $\partial K \in \cE_{h}^{i,s}$.   The penalty parameter that is defined as,,
\begin{equation}
\tau \; = \; \left\{ 
\begin{array}{cc}
\tilde{\tau} \, \min \left( h^{-1}_{1}, h^{-1}_{2} \right) & \textbf{x} \in \langle K_{1}, K_{2} \rangle \\
\tilde{\tau}  \, h^{-1} & \textbf{x} \in \partial K \cap \in \Gamma_{C} 
\end{array}
\right. 
\label{eq:Penalty}
\end{equation}
\noindent
where $\tilde{\tau}$ is a postive number and $h$ is the diameter of the element $K$.

We can now substitute~\eqref{EQ:U Flux} and~\eqref{EQ:Q Flux} into~\eqref{EQ:DDP LDG1} and ~\eqref{EQ:DDP LDG2} to obtain the solution pair $(\rho_n, \bq_n)$ to the semi-discrete LDG approximation for the drift-diffusion equation given by the following problem: find $(\rho_n, \bq_n) \in \left( C\left((0,T]; W_{h,k}\right) \times C((0, T]; \widetilde \bW_{h,k}) \right)$ such that,
\begin{eqnarray}
\label{EQ:LDG I}\left(w, \partial_{t} \rho_n \right)_{\Omega_S}
+ \cL_n\left(w, \bw; \rho_n, \bq_n \right)+\langle w , I_{et}(\rho_n-\rho_n^e, \rho_o) \rangle_{\Sigma} 
&=& -\big( w ,  \widetilde R(\rho_n,\rho_p) \big)_{\Omega_S}\\
\label{EQ:LDG II} \cM_n\left(w, \bw; \rho_n, \bq_n\right) + \left(\bw, \alpha_n\bE\rho_n\right)_{\Omega_S}
&=&-\langle \bn\cdot \bw, \rho_n^e \rangle_{\Gamma_C},
\end{eqnarray}
for all $(w,\bw) \in W_{h,k} \times \widetilde \bW_{h,k}$, where the quad-linear forms $\cL_n(w, \bw; \rho_n, \bq_n)$ and $\cM_n(w, \bw; \rho_n, \bq_n)$ are defined as
\begin{eqnarray}
\cL_n= \langle \llbracket w \rrbracket,  \tau \llbracket \rho_n \rrbracket \rangle_{\cE_{h}^{i,s}}
 - \left(\nabla w, \bq_n \right)_{\Omega_S}+
\langle \llbracket w \rrbracket, 
\left\{\bq_n \right\}-\llbracket \bq_n \rrbracket \boldsymbol \beta \rangle_{\cE_{h}^{i,s}} 
+\langle w, \bn\cdot \bq_n \rangle_{\Gamma_C},\\
\cM_n= - \left(\nabla \cdot \bw, \rho_n\right)_{\Omega_S} + 
\langle  \llbracket \bw \rrbracket, \left\{\rho_n \right\} + \boldsymbol \beta \cdot \llbracket \rho_n \rrbracket \rangle_{\cE_h^{i,s}} +
\langle \bn\cdot \bw, \rho_n \rangle_{\Gamma_{N}^s \cup \Sigma} + \left(\bw, \mu_n^{-1}\bq_n \right)_{\Omega_S},
\end{eqnarray}
with the function $\widetilde R (\rho_n, \rho_p) =R(\rho_n, \rho_p)- \gamma G(\bx)$.

We note that in our formulation the auxillary variable $\bq_n$ is defined as the total flux since our boundary and interface conditions are all directly given in terms of the total flux.  This choice of auxillary variable eliminates the need to introduce the up-winding numerical fluxes for the drift component.  However, choosing the value of $\boldsymbol \beta$ properly one can recover the up-winding flux for the drift component.

\section{Implicit-explicit time stepping algorithms}
\label{SEC:IMEX}

We can combine the mixed finite element discretization of the Poisson equation~\eqref{EQ:MFEM I}-\eqref{EQ:MFEM II} with the semi-discrete LDG discretization of the drift-diffusion equation~\eqref{EQ:LDG I}-\eqref{EQ:LDG II} to get a discretization scheme for the whole system of coupled equations in~\eqref{EQ:DDP S} and~\eqref{EQ:DDP E}. The resulting system will be nonlinear due to the various nonlinearities in the PDE system, including for instance, the nonlinear recombination-generation effects, the coupling between the drift-diffusion components of the densities through drift terms such as $\alpha_n\bE \rho_n$, as well as the nonlinear interface transfer conditions.

To deal with the nonlinearities in the system, we combine different time-stepping schemes with domain decomposition techniques. More precisely, we use explicit time-stepping to treat the recombination-generation nonlinearity, use a ``time lagging'' technique~\cite{ArJuPoWh-CG13,Coats-SPE03,Selberherr-Book01,ShCa-SPE59,StGa-SPE61} to treat the coupling between the densities and the electric potential, and use a Schwarz domain decomposition technique~\cite{CaBaKlKo-JCP12,ChELiSh-JCP07,Lions-DDM87,Lions-DDM89, Lions-DDM90,MiQuSa-JCP95,QuVa-Book99,ToWi-Book05}, coupled with explicit stepping, to treat the coupling between the semiconductor and electrolyte domains. To overcome the limitation on the size of the time steps imposed by the Courant-Friedrichs-Lewy (CFL) condition when an explicit time stepping scheme is used, we use implicit scheme on the diffusion terms whenever it is possible. Therefore, our marching in time is done by an overall implicit-explicit (IMEX) time stepping scheme~\cite{Selberherr-Book01}.

We present two classes of time stepping algorithms that are based on IMEX time stepping and domain decomposition methods. The first class, in Sections~\ref{SUBSEC:AS-IMIMEX} and ~\ref{SUBSEC:AS-IMEXEX}, uses the alternating Schwarz method. The second class, in Sections~\ref{SUBSEC:PS-IMEXEX} and ~\ref{SUBSEC:TsPS-IMEXEX}, uses the parallel Schwarz strategy. The main difference between the two classes of methods lies in how they treat the carrier densities on the semiconductor-electrolyte interface.  We remark that in this work we only consider first-order time discretization schemes.  The reason being is that we are interested in the steady state characteristics of PECs and the accuracy of the temporal discretization did not appreciably affect the accuracy of the steady state numerical solutions.

\subsection{The implicit-implicit-explicit alternating Schwarz scheme}
\label{SUBSEC:AS-IMIMEX}

In the implicit-implicit-explicit alternating Schwarz (AS-IMIMEX) scheme, we first use a ``time lagging'' technique to decouple the Poisson equations from the drift-diffusion equations. We then apply an implicit scheme on the diffusion terms, an implicit scheme on the drift terms, and an explicit scheme on the recombination-generation terms. To decouple the equations in the semiconductor domain from these in the electrolyte domain, we use a strategy based on the alternating Schwarz domain decomposition method.

The algorithm works as follows. At time step $k$, we first solve the Poisson equation for the electric potential (and therefore the electric field) using charge densities at the previous time step $k-1$. That is, we solve for $(\Phi^{\blue{k-1}},\bE^{\blue{k-1}})$ from the Poisson equation:
\begin{equation}\label{EQ:MFEM III}
\begin{array}{rcl}
\left(\bw, \lambda^{-2} \bE^{\blue{k-1}} \right)_{\Omega}-\left(\nabla\cdot\bw, \Phi^{\blue{k-1}} \right)_{\Omega} &=& -\langle   \bn\cdot\bw, \Phi_{\rm app} + \Phi_{\rm bi} \rangle_{\Gamma_C} - \langle \bn\cdot\bw, \Phi^{\infty} \rangle_{\Gamma_A},\\
\left(w, \nabla\cdot \bE^{\blue{k-1}} \right)_{\Omega} &=& \left(w, f^{\red{k-1}} \right)_{\Omega},
\end{array}
\end{equation}
where $f^{\red{k-1}}$ is evaluated with the densities $(\rho_n^{\red{k-1}}, \rho_p^{\red{k-1}}, \rho_r^{\red{k-1}},\rho_o^{\red{k-1}})$. We then use $\bE^{k-1}$ in the drift-diffusion equations for the charge densities to find updates of these densities. This time lagging technique is called the IMPES (implicit pressure, explicit saturation) method in the subsurface flow simulation literature~\cite{ArJuPoWh-CG13,Coats-SPE03,Selberherr-Book01,ShCa-SPE59,StGa-SPE61}, and is conceptually the same as the Gummel iteration method~\cite{BeCaCaVe-JCP09,BuPi-M3AS09,Kulikovsky-JCP95} which is often used semiconductor device simulations to decouple the Poisson component and the drift-diffusion component of the system in steady state.

To update the charge densities using the drift-diffusion equations, we use a alternating Schwarz domain decomposition idea~\cite{CaBaKlKo-JCP12,ChELiSh-JCP07,Lions-DDM87,Lions-DDM89, Lions-DDM90,MiQuSa-JCP95,QuVa-Book99,ToWi-Book05} to decouple the equations in the semiconductor domain with those in the electrolyte domain. To be precise, we first update the electron and hole densities. This would require the densities of the reductant and the oxidant at the semiconductor-electrolyte interface $\Sigma$, in the interface conditions~\eqref{EQ:Int Rho n} and ~\eqref{EQ:Int Rho p}. We use the $k-1$ step values of these densities, $\rho_{r}^{\red{k-1}}$ and $\rho_{o}^{\red{k-1}}$. Therefore, we solve for $(\rho_n^{\blue{k}}, \rho_p^{\blue{k}})$ from:
\begin{multline}\label{EQ:DDP Disc n1}
(w, \rho_{n}^{\blue{k}})_{\Omega_{S}} 
+ \Delta t_k \cL_n \left(w, \bw; \rho_{n}^{\blue{k}}, \bq_{n}^{\blue{k}} \right)+ \Delta t_k \langle w,  I_{et}(\rho_{n}^{\blue{k}}, \rho_{o}^{\red{k-1}})  \rangle_{\Sigma} \\ 
= (w, \rho_{n}^{\red{k-1}})_{\Omega_{S}} +\Delta t_k \langle w, I_{et}(\rho_{n}^{e}, \rho_{o}^{\red{k-1}}) \rangle_{\Sigma}-\Delta t_k (w, \widetilde R(\rho_n^{\red{k-1}},\rho_p^{\red{k-1}}))_{\Omega_S},
\end{multline}
\begin{equation}\label{EQ:DDP Disc n2}
\cM_n\left(w, \bw; \rho_n^{\blue{k}}, \bq_n^{\blue{k}}\right)+\left(\bw, \alpha_n\bE^{\red{k-1}}\rho_n^{\blue{k}} \right)_{\Omega_S}=-\langle \bn\cdot \bw, \rho_n^e\rangle_{\Gamma_C},
\end{equation}
\begin{multline}\label{EQ:DDP Disc p1}
(w, \rho_{p}^{\blue{k}})_{\Omega_{S}} 
+ \Delta t_k \cL_p \left(w, \bw; \rho_{p}^{\blue{k}}, \bq_{p}^{\blue{k}} \right)+ \Delta t_k \langle w,  I_{ht}(\rho_{p}^{\blue{k}}, \rho_{r}^{\red{k-1}})  \rangle_{\Sigma} \\ 
= (w, \rho_{p}^{\red{k-1}})_{\Omega_{S}} +\Delta t_k \langle w, I_{ht}(\rho_{p}^{e}, \rho_{r}^{\red{k-1}}) \rangle_{\Sigma}-\Delta t_k (w, \widetilde R(\rho_n^{\red{k-1}},\rho_p^{\red{k-1}}))_{\Omega_S},
\end{multline}
\begin{equation}\label{EQ:DDP Disc p2}
\cM_p\left(w, \bw; \rho_p^{\blue{k}}, \bq_p^{\blue{k}}\right)+\left(\bw, \alpha_p\bE^{\red{k-1}}\rho_p^k \right)_{\Omega_S}=-\langle \bn\cdot \bw, \rho_p^e\rangle_{\Gamma_C}.
\end{equation}
Here the diffusion terms are treated implicitly, the drift terms are treated implicitly, while the recombination-generation terms (which is now incorporated into the function $\widetilde R$) are treated explicitly.

Now that we have obtained the electron and hole densities at step $k$, $(\rho_{n}^{\red{k}}, \rho_{p}^{\red{k}})$, we can use these values in the interface conditions~\eqref{EQ:Int Rho r} to update the reductant density:
\begin{multline}\label{EQ:DDP Disc r1}
(w, \rho_{r}^{\blue{k}})_{\Omega_{E}} 
+ \Delta t_k \cL_r \left(w, \bw; \rho_{r}^{\blue{k}}, \bq_{r}^{\blue{k}} \right) + \Delta t_k \langle w,  I_{ht}(\rho_{p}^{\red{k}}-\rho_p^e, \rho_{r}^{\blue{k}})  \rangle_{\Sigma} \\ 
= (w, \rho_{r}^{\red{k-1}})_{\Omega_{E}} +\Delta t_k \langle w, I_{et}(\rho_n^{\red{k}}-\rho_{n}^{e}, \rho_{o}^{\red{k}}) \rangle_{\Sigma},
\end{multline}
\begin{equation}\label{EQ:DDP Disc r2}
\cM_r\left(w, \bw; \rho_r^{\blue{k}}, \bq_r^{\blue{k}}\right)+\left(\bw, \alpha_r\bE^{\red{k-1}}\rho_r^{\blue{k}} \right)_{\Omega_E}=-\langle \bn\cdot \bw, \rho_r^\infty \rangle_{\Gamma_A}.
\end{equation}
With the updated values $(\rho_{n}^{\red{k}}, \rho_{p}^{\red{k}}, \rho_{r}^{\red{k}})$, we can now update the oxidant density following the equations:
\begin{multline}\label{EQ:DDP Disc o1}
(w, \rho_{o}^{\blue{k}})_{\Omega_{E}} 
+ \Delta t_k \cL_o \left(w, \bw; \rho_{o}^{\blue{k}}, \bq_{o}^{\blue{k}} \right) + \Delta t_k \langle w,  I_{et}(\rho_{n}^{\red{k}}-\rho_n^e, \rho_{o}^{\blue{k}})  \rangle_{\Sigma} \\ 
= (w, \rho_{o}^{\red{k-1}})_{\Omega_{E}} +\Delta t_k \langle w, I_{ht}(\rho_p^{\red{k}}-\rho_{p}^{e}, \rho_{r}^{\red{k}}) \rangle_{\Sigma},
\end{multline}
\begin{equation}\label{EQ:DDP Disc o2}
\cM_o\left(w, \bw; \rho_o^{\blue{k}}, \bq_o^{\blue{k}}\right)+\left(\bw, \alpha_o\bE^{\red{k-1}}\rho_o^{\blue{k}} \right)_{\Omega_E}=-\langle \bn\cdot \bw, \rho_o^\infty \rangle_{\Gamma_A}.
\end{equation}

The flow of the AS-IMIMEX algorithm is summarized in Algorithm~\ref{ALG:AS-IMIMEX}. Our experience in the numerical experiments is that the use of implicit time stepping for the diffusion terms allows us to use time steps that are a few order of magnitude greater than what an explicit scheme, for instance the adaptive forward Euler scheme, on the diffusion terms would allow. However, treating the drift term implicitly is very expensive. This is because the electric field changes at every time step. Therefore, we have to update the LDG matrices at each time step of the calculation. Assembling and factorizing the LDG matrices take up a significant portion of the run-time.

\begin{algorithm}
\caption{The AS-IMIMEX [resp. AS-IMEXEX] Algorithm}
\label{ALG:AS-IMIMEX}
\begin{algorithmic}[1]
\State Initialize the density data ($\rho_{n}^0$, $\rho_{p}^0$,  $\rho_{r}^0$, $\rho_{o}^0$); set $t_0=0$; set $k=1$;
\While{not reaching steady state} 
	\State Solve the Poisson problem~\eqref{EQ:MFEM III} for ($\Phi^{\blue{k-1}}$, $\bE^{\blue{k-1}}$) using ($\rho_{n}^{\red{k-1}}, \rho_{p}^{\red{k-1}},  \rho_{r}^{\red{k-1}}, \rho_{o}^{\red{k-1}}$) as data.
	\State Determine $\Delta t_{k}$ from the CFL condition with $\bE^{\red{k-1}}$ and ($\rho_{n}^{\red{k-1}}, \rho_{p}^{\red{k-1}},  \rho_{r}^{\red{k-1}}, \rho_{o}^{\red{k-1}}$).
	\State Update electron and hole densities from~\eqref{EQ:DDP Disc n1}-\eqref{EQ:DDP Disc p2} [resp. ~\eqref{EQ:DDP Disc n1 b}-\eqref{EQ:DDP Disc p2 b}].
	\State Update reductant density from~\eqref{EQ:DDP Disc r1}-\eqref{EQ:DDP Disc r2}  [resp. ~\eqref{EQ:DDP Disc r1 b}-\eqref{EQ:DDP Disc r2 b}] using $(\rho_n^{\red{k}}, \rho_p^{\red{k}})$.
	\State Update oxidant density from~\eqref{EQ:DDP Disc o1}-\eqref{EQ:DDP Disc o2} [resp. ~\eqref{EQ:DDP Disc o1 b}-\eqref{EQ:DDP Disc o2 b}] using $(\rho_n^{\red{k}}, \rho_p^{\red{k}}, \rho_r^{\red{k}})$.
	\State $t_{k} = t_{k-1} + \Delta t_{k}$, $k = k+1$
\EndWhile
\end{algorithmic}
\end{algorithm}

\subsection{The implicit-explicit-explicit alternating Schwarz scheme}
\label{SUBSEC:AS-IMEXEX}

The implicit-explicit-explicit alternating Schwarz (AS-IMEXEX) scheme is similar to the AS-IMIMEX scheme in the previous section. The only difference is that here we treat the drift terms in the drift-diffusion equations explicitly in time. That is, the charge densities in the drift terms are taken as there values at the previous step $k-1$, not the current $k$.

The algorithm works as follows. At time step $k$, we first solve the Poisson equation~\eqref{EQ:MFEM III} for $(\Phi^{k-1},\bE^{k-1})$. We then update the electron and hole densities as in the AS-IMIMEX scheme, but treat the drift terms explicitly to get $(\rho_n^{\blue{k}}, \rho_p^{\blue{k}})$:
\begin{multline}\label{EQ:DDP Disc n1 b}
(w, \rho_{n}^{\blue{k}})_{\Omega_{S}} 
+ \Delta t_k \cL_n \left(w, \bw; \rho_{n}^{\blue{k}}, \bq_{n}^{\blue{k}} \right)+ \Delta t_k \langle w,  I_{et}(\rho_{n}^{\blue{k}}, \rho_{o}^{\red{k-1}})  \rangle_{\Sigma} \\ 
= (w, \rho_{n}^{\red{k-1}})_{\Omega_{S}} +\Delta t_k \langle w, I_{et}(\rho_{n}^{e}, \rho_{o}^{\red{k-1}}) \rangle_{\Sigma}-\Delta t_k (w, \widetilde R(\rho_n^{\red{k-1}},\rho_p^{\red{k-1}}))_{\Omega_S},
\end{multline}
\begin{equation}\label{EQ:DDP Disc n2 b}
\cM_n\left(w, \bw; \rho_n^{\blue{k}}, \bq_n^{\blue{k}}\right)=-\left(\bw, \alpha_n\bE^{\red{k-1}}\rho_n^{\red{k-1}} \right)_{\Omega_S}-\langle \bn\cdot \bw, \rho_n^e\rangle_{\Gamma_C},
\end{equation}
\begin{multline}\label{EQ:DDP Disc p1 b}
(w, \rho_{p}^{\blue{k}})_{\Omega_{S}} 
+ \Delta t_k \cL_p \left(w, \bw; \rho_{p}^{\blue{k}}, \bq_{p}^{\blue{k}} \right)+ \Delta t_k \langle w,  I_{ht}(\rho_{p}^{\blue{k}}, \rho_{r}^{\red{k-1}})  \rangle_{\Sigma} \\ 
= (w, \rho_{p}^{\red{k-1}})_{\Omega_{S}} +\Delta t_k \langle w, I_{ht}(\rho_{p}^{e}, \rho_{r}^{\red{k-1}}) \rangle_{\Sigma}-\Delta t_k (w, \widetilde R(\rho_n^{\red{k-1}},\rho_p^{\red{k-1}}))_{\Omega_S},
\end{multline}
\begin{equation}\label{EQ:DDP Disc p2 b}
\cM_p\left(w, \bw; \rho_p^{\blue{k}}, \bq_p^{\blue{k}}\right)=-\left(\bw, \alpha_p\bE^{\red{k-1}}\rho_p^{\red{k-1}} \right)_{\Omega_S}-\langle \bn\cdot \bw, \rho_p^e\rangle_{\Gamma_C}.
\end{equation}
Now that we have obtained the electron and hole densities at step $k$, $(\rho_{n}^{\red{k}}, \rho_{p}^{\red{k}})$, we can use these values to update the reductant density:
\begin{multline}\label{EQ:DDP Disc r1 b}
(w, \rho_{r}^{\blue{k}})_{\Omega_{E}} 
+ \Delta t_k \cL_r \left(w, \bw; \rho_{r}^{\blue{k}}, \bq_{r}^{\blue{k}} \right) + \Delta t_k \langle w,  I_{ht}(\rho_{p}^{\red{k}}-\rho_p^e, \rho_{r}^{\blue{k}})  \rangle_{\Sigma} \\ 
= (w, \rho_{r}^{\red{k-1}})_{\Omega_{E}} +\Delta t_k \langle w, I_{et}(\rho_n^{\red{k}}-\rho_{n}^{e}, \rho_{o}^{\red{k}}) \rangle_{\Sigma},
\end{multline}
\begin{equation}\label{EQ:DDP Disc r2 b}
\cM_r\left(w, \bw; \rho_r^{\blue{k}}, \bq_r^{\blue{k}}\right)=-\left(\bw, \alpha_r\bE^{\red{k-1}}\rho_r^{\red{k-1}} \right)_{\Omega_E}-\langle \bn\cdot \bw, \rho_r^\infty \rangle_{\Gamma_A}.
\end{equation}
Note again that the drift term is treated explicitly here. With the updated values $(\rho_{n}^{\red{k}}, \rho_{p}^{\red{k}}, \rho_{r}^{\red{k}})$, we can now update the oxidant density following the equations:
\begin{multline}\label{EQ:DDP Disc o1 b}
(w, \rho_{o}^{\blue{k}})_{\Omega_{E}} 
+ \Delta t_k \cL_o \left(w, \bw; \rho_{o}^{\blue{k}}, \bq_{o}^{\blue{k}} \right) + \Delta t_k \langle w,  I_{et}(\rho_{n}^{\red{k}}-\rho_n^e, \rho_{o}^{\blue{k}})  \rangle_{\Sigma} \\ 
= (w, \rho_{o}^{\red{k-1}})_{\Omega_{E}} +\Delta t_k \langle w, I_{ht}(\rho_p^{\red{k}}-\rho_{p}^{e}, \rho_{r}^{\red{k}}) \rangle_{\Sigma},
\end{multline}
\begin{equation}\label{EQ:DDP Disc o2 b}
\cM_o\left(w, \bw; \rho_o^{\blue{k}}, \bq_o^{\blue{k}}\right)=-\left(\bw, \alpha_o\bE^{\red{k-1}}\rho_o^{\red{k-1}} \right)_{\Omega_E}-\langle \bn\cdot \bw, \rho_o^\infty \rangle_{\Gamma_A}.
\end{equation}

The overall flow of the AS-IMEXEX algorithm is identical to that of the AS-IMIMEX algorithm, and is summarized in Algorithm~\ref{ALG:AS-IMIMEX}, in the brackets. As discussed in Section~\ref{SEC:Num}, numerical experiments showed that treating the drift term explicitly saves significant computational cost overall.

\subsection{The implicit-explicit-explicit parallel Schwarz scheme}
\label{SUBSEC:PS-IMEXEX}

In the implicit-explicit-explicit parallel Schwarz scheme, we treat the diffusion terms implicitly, the drift terms explicitly, and the recombination-generation term explicitly, in the same way as in the AS-IMEXEX scheme. What is different here is that we now treat the reaction interface conditions explicitly in all drift-diffusion equations. This way, all four drift-diffusion equations are decoupled from each other and can be solved simultaneously in parallel. To be precise, we use the values of all the densities at time step $k-1$ in the functions $I_{et}$ and $I_{ht}$ in interface conditions~\eqref{EQ:Int Rho n}-\eqref{EQ:Int Rho o} to update the values of the densities simultaneously to the current time step:
\begin{multline}\label{EQ:DDP Disc n1 c}
(w, \rho_{n}^{\blue{k}})_{\Omega_{S}} 
+ \Delta t \cL_n \left(w, \bw; \rho_{n}^{\blue{k}}, \bq_{n}^{\blue{k}} \right)\\ 
= (w, \rho_{n}^{\red{k-1}})_{\Omega_{S}} -\Delta t \langle w, I_{et}(\rho_n^{\red{k-1}}-\rho_{n}^{e}, \rho_{o}^{\red{k-1}}) \rangle_{\Sigma}-\Delta t (w, \widetilde R(\rho_n^{\red{k-1}},\rho_p^{\red{k-1}}))_{\Omega_S},
\end{multline}
\begin{equation}\label{EQ:DDP Disc n2 c}
\cM_n\left(w, \bw; \rho_n^{\blue{k}}, \bq_n^{\blue{k}}\right)=-\left(\bw, \alpha_n\bE^{\red{k-1}}\rho_n^{\red{k-1}} \right)_{\Omega_S}-\langle \bn\cdot \bw, \rho_n^e\rangle_{\Gamma_C},
\end{equation}
\begin{multline}\label{EQ:DDP Disc p1 c}
(w, \rho_{p}^{\blue{k}})_{\Omega_{S}} 
+ \Delta t \cL_p \left(w, \bw; \rho_{p}^{\blue{k}}, \bq_{p}^{\blue{k}} \right)\\ 
= (w, \rho_{p}^{\red{k-1}})_{\Omega_{S}} -\Delta t \langle w, I_{ht}(\rho_p^{\red{k-1}}-\rho_{p}^{e}, \rho_{r}^{\red{k-1}}) \rangle_{\Sigma}-\Delta t (w, \widetilde R(\rho_n^{\red{k-1}},\rho_p^{\red{k-1}}))_{\Omega_S},
\end{multline}
\begin{equation}\label{EQ:DDP Disc p2 c}
\cM_p\left(w, \bw; \rho_p^{\blue{k}}, \bq_p^{\blue{k}}\right)=-\left(\bw, \alpha_p\bE^{\red{k-1}}\rho_p^{\red{k-1}} \right)_{\Omega_S}-\langle \bn\cdot \bw, \rho_p^e\rangle_{\Gamma_C},
\end{equation}
\begin{multline}\label{EQ:DDP Disc r1 c}
(w, \rho_{r}^{\blue{k}})_{\Omega_{E}} 
+ \Delta t \cL_r \left(w, \bw; \rho_{r}^{\blue{k}}, \bq_{r}^{\blue{k}} \right)\\ 
= (w, \rho_{r}^{\red{k-1}})_{\Omega_{E}} - \Delta t \langle w,  I_{ht}(\rho_{p}^{\red{k-1}}-\rho_p^e, \rho_{r}^{\red{k-1}}) - I_{et}(\rho_n^{\red{k-1}}-\rho_{n}^{e}, \rho_{o}^{\red{k-1}}) \rangle_{\Sigma},
\end{multline}
\begin{equation}\label{EQ:DDP Disc r2 c}
\cM_r\left(w, \bw; \rho_r^{\blue{k}}, \bq_r^{\blue{k}}\right)=-\left(\bw, \alpha_r\bE^{\red{k-1}}\rho_r^{\red{k-1}} \right)_{\Omega_E}-\langle \bn\cdot \bw, \rho_r^\infty \rangle_{\Gamma_A},
\end{equation}
\begin{multline}\label{EQ:DDP Disc o1 c}
(w, \rho_{o}^{\blue{k}})_{\Omega_{E}} 
+ \Delta t \cL_o \left(w, \bw; \rho_{o}^{\blue{k}}, \bq_{o}^{\blue{k}} \right) \\ 
= (w, \rho_{o}^{\red{k-1}})_{\Omega_{E}}- \Delta t \langle w,  I_{et}(\rho_{n}^{\red{k-1}}-\rho_n^e, \rho_{o}^{\red{k-1}}) - I_{ht}(\rho_p^{\red{k-1}}-\rho_{p}^{e}, \rho_{r}^{\red{k-1}}) \rangle_{\Sigma},
\end{multline}
\begin{equation}\label{EQ:DDP Disc o2 c}
\cM_o\left(w, \bw; \rho_o^{\blue{k}}, \bq_o^{\blue{k}}\right)=-\left(\bw, \alpha_o\bE^{\red{k-1}}\rho_o^{\red{k-1}} \right)_{\Omega_E}-\langle \bn\cdot \bw, \rho_o^\infty \rangle_{\Gamma_A}.
\end{equation}

The PS-IMEXEX algorithm is summarized in Algorithm~\ref{ALG:PS-IMEXEX}. Note that due to the fact that all the nonlinear terms are treated explicitly, we do not need to adjust the size of the time step $\Delta t$ during the time evolution. Moreover, we only need to factorize the system matrix once at the beginning of this algorithm.

\begin{algorithm}
\caption{The PS-IMEXEX Algorithm}
\label{ALG:PS-IMEXEX}
\begin{algorithmic}[1]
\State Initialize the density data ($\rho_{n}^0$, $\rho_{p}^0$,  $\rho_{r}^0$, $\rho_{o}^0$); set $t_0=0$; set $k=1$;
\State Determine $\Delta t$ from the CFL condition;
\While{not reaching steady state} 
	\State Solve the Poisson problem~\eqref{EQ:MFEM III} for $(\Phi^{\blue{k-1}}, \bE^{\blue{k-1}})$ using $(\rho_{n}^{\red{k-1}}, \rho_{p}^{\red{k-1}},  \rho_{r}^{\red{k-1}}, \rho_{o}^{\red{k-1}})$ as data;
	\State Update the densities from~\eqref{EQ:DDP Disc n1 c}-\eqref{EQ:DDP Disc o2 c} using $\bE^{\red{k-1}}$ and $(\rho_{n}^{\red{k-1}}, \rho_{p}^{\red{k-1}},  \rho_{r}^{\red{k-1}}, \rho_{o}^{\red{k-1}})$ as data;
	\State $t_{k} = t_{k-1} + \Delta t$, $k = k+1$;
\EndWhile
\end{algorithmic}
\end{algorithm}	

\subsection{The two-scale PS-IMEXEX scheme}
\label{SUBSEC:TsPS-IMEXEX}

In the PS-IMEXEX algorithm, we choose a fixed time step $\Delta t$ even though it can be calculated in an adaptive manner. In our numerical experiments, we observe that we can take a $\Delta t$ that is much larger than the one that we calculated using the CFL condition for the drift-diffusion equations in the semiconductor system (mainly based on the operator pairs $(\cL_n, \cM_n)$ and $(\cL_p, \cM_p)$. Intuitively, this is due to the fact that the characteristic time scale for the semiconductor drift-diffusion system and that for the electrolyte drift-diffusion system are quite different. In general, the electrolyte system evolves in a much slower pace than the semiconductor system. To take advantage of this scale separation, we implemented the two-scale PS-IMEXEX (TsPU-IMEXEX) algorithm. The main idea here is to update quantities in the semiconductor domain more frequently than the quantities in the electrolyte domain. 

Let $\Delta t_s$ and $\Delta t_e$ be the time steps given by the CFL conditions in the semiconductor and electrolyte domains respectively. Let $K$ be the positive integer such that $\Delta t_e=K\Delta t_s$. At time step $k$, we first solve the Poisson equation~\eqref{EQ:MFEM III} to get $(\Phi^{k-1}, \bE^{k-1})$. We then perform $K$ time steps for the semiconductor system, starting with the initial condition $(\tilde\rho_n^0, \tilde\rho_p^0)=(\rho_n^{\red{k-1}}, \rho_p^{\red{k-1}})$:
\begin{multline}\label{EQ:DDP Disc n1 d}
(w, \tilde\rho_{n}^{\blue{j}})_{\Omega_{S}} 
+ \Delta t_s \cL_n \left(w, \bw; \tilde\rho_{n}^{\blue{j}}, \tilde \bq_{n}^{\blue{j}} \right)\\ 
= (w, \rho_{n}^{{\red{j-1}}})_{\Omega_{S}} -\Delta t_s \langle w, I_{et}(\tilde\rho_n^{\red{j-1}}-\rho_{n}^{e}, \rho_{o}^{\red{k-1}}) \rangle_{\Sigma}-\Delta t_s (w, \widetilde R(\tilde\rho_n^{\red{j-1}},\tilde\rho_p^{\red{j-1}}))_{\Omega_S},
\end{multline}
\begin{equation}\label{EQ:DDP Disc n2 d}
\cM_n\left(w, \bw; \tilde\rho_n^{\blue{j}}, \tilde\bq_n^{\blue{j}}\right)=-\left(\bw, \alpha_n\bE^{\red{k-1}}\tilde\rho_n^{j-1} \right)_{\Omega_S}-\langle \bn\cdot \bw, \rho_n^e\rangle_{\Gamma_C},
\end{equation}
\begin{multline}\label{EQ:DDP Disc p1 d}
(w, \tilde\rho_{p}^{\blue{j}})_{\Omega_{S}} 
+ \Delta t_s \cL_p \left(w, \bw; \tilde\rho_{p}^{\blue{j}}, \tilde\bq_{p}^{\blue{j}} \right)\\ 
= (w, \tilde\rho_{p}^{\red{j-1}})_{\Omega_{S}} -\Delta t_s \langle w, I_{ht}(\tilde\rho_p^{\red{j-1}}-\rho_{p}^{e}, \rho_{r}^{\red{k-1}}) \rangle_{\Sigma}-\Delta t_s (w, \widetilde R(\tilde\rho_n^{\red{j-1}},\tilde\rho_p^{\red{j-1}}))_{\Omega_S},
\end{multline}
\begin{equation}\label{EQ:DDP Disc p2 d}
\cM_p\left(w, \bw; \tilde\rho_p^{\blue{j}}, \tilde\bq_p^{\blue{j}}\right)=-\left(\bw, \alpha_p\bE^{\red{k-1}}\tilde\rho_p^{\red{j-1}} \right)_{\Omega_S}-\langle \bn\cdot \bw, \rho_p^e\rangle_{\Gamma_C},
\end{equation}
Note that the values for the densities of the reductants and the oxidants are kept unchanged during this iteration.

We then set $(\rho_n^{\blue{k-1}}, \rho_p^{\blue{k-1}})=(\tilde\rho_n^{\red{K}}, \tilde \rho_p^{\red{K}})$, and update the density values of the reductants and the oxidants using time step $\Delta t_e$ following~\eqref{EQ:DDP Disc r1 b}-\eqref{EQ:DDP Disc o2 b}, that is:
\begin{multline}\label{EQ:DDP Disc r1 d}
(w, \rho_{r}^{\blue{k}})_{\Omega_{E}} 
+ \Delta t_e \cL_r \left(w, \bw; \rho_{r}^{\blue{k}}, \bq_{r}^{\blue{k}} \right)\\ 
= (w, \rho_{r}^{\red{k-1}})_{\Omega_{E}} - \Delta t_e \langle w,  I_{ht}(\rho_{p}^{\red{k-1}}-\rho_p^e, \rho_{r}^{\red{k-1}}) - I_{et}(\rho_n^{\red{k-1}}-\rho_{n}^{e}, \rho_{o}^{\red{k-1}}) \rangle_{\Sigma},
\end{multline}
\begin{equation}\label{EQ:DDP Disc r2 d}
\cM_r\left(w, \bw; \rho_r^{\blue{k}}, \bq_r^{\blue{k}}\right)=-\left(\bw, \alpha_r\bE^{\red{k-1}}\rho_r^{\red{k-1}} \right)_{\Omega_E}-\langle \bn\cdot \bw, \rho_r^\infty \rangle_{\Gamma_A},
\end{equation}
\begin{multline}\label{EQ:DDP Disc o1 d}
(w, \rho_{o}^{\blue{k}})_{\Omega_{E}} 
+ \Delta t_e \cL_o \left(w, \bw; \rho_{o}^{\blue{k}}, \bq_{o}^{\blue{k}} \right) \\ 
= (w, \rho_{o}^{\red{k-1}})_{\Omega_{E}}- \Delta t_e \langle w,  I_{et}(\rho_{n}^{\red{k-1}}-\rho_n^e, \rho_{o}^{\red{k-1}}) - I_{ht}(\rho_p^{\red{k-1}}-\rho_{p}^{e}, \rho_{r}^{\red{k-1}}) \rangle_{\Sigma},
\end{multline}
\begin{equation}\label{EQ:DDP Disc o2 d}
\cM_o\left(w, \bw; \rho_o^{\blue{k}}, \bq_o^{\blue{k}}\right)=-\left(\bw, \alpha_o\bE^{\red{k-1}}\rho_o^{\red{k-1}} \right)_{\Omega_E}-\langle \bn\cdot \bw, \rho_o^\infty \rangle_{\Gamma_A}.
\end{equation}

The algorithm is summarized in Algorithm~\ref{ALG:TsPS-IMEXEX}. In the numerical simulations, we are relatively conservative on the selection of $K$. We take a $K$ that is smaller than the one determined by $\Delta t_s$ and $\Delta t_e$. We are able to significantly achieve faster simulations with $K$ between $5$ and $10$, without sacrificing the accuracy and stability of the algorithm.

\begin{algorithm}
\caption{The TsPS-IMEXEX Algorithm}
\label{ALG:TsPS-IMEXEX}
\begin{algorithmic}[1]
\State Initialize the density data ($\rho_{n}^0$, $\rho_{p}^0$,  $\rho_{r}^0$, $\rho_{o}^0$); set $t_0=0$; set $k=1$;
\State Determine $\Delta t_s$ from $(\cL_n, \cM_n)$; determine $\Delta t_e$ from $(\cL_r, \cM_r)$; set $K=\Delta t_e/\Delta t_s$;
\While{not reaching steady state}
	\State Solve the Poisson problem~\eqref{EQ:MFEM III} for $(\Phi^{\red{k-1}}, \bE^{\red{k-1}})$ using $(\rho_{n}^{\red{k-1}}, \rho_{p}^{\red{k-1}},  \rho_{r}^{\red{k-1}}, \rho_{o}^{\red{k-1}})$ as data;
	\State Set $(\tilde\rho_n^0, \tilde\rho_p^0)=(\rho_n^{k-1}, \rho_p^{k-1})$;
	\For{j=0, $\cdots$, K}
		\State Update the densities $(\tilde\rho_n, \tilde\rho_p)$ according to~\eqref{EQ:DDP Disc n1 d}-\eqref{EQ:DDP Disc p2 d};
	\EndFor
	\State Set $(\rho_n^{k-1}, \rho_p^{k-1})=(\tilde\rho_n^K, \tilde \rho_p^K)$;
	\State Update the reductant and oxidant densities according to~\eqref{EQ:DDP Disc r1 d}-\eqref{EQ:DDP Disc o2 d};
	\State $t_{k} = t_{k-1} + \Delta t_e$, $k = k+1$;
\EndWhile
\end{algorithmic}
\end{algorithm}	

\section{Numerical experiments}
\label{SEC:Num}

We now present some numerical simulations in one and two dimensions to demonstrate the performance of the algorithms we developed in this work. We also intend to investigate the characteristics of PEC solar cells under various conditions. In our implementation of the algorithms in one dimension, we use the SparseLU~\cite{DeEiGiLiLi-SIAM99} direct solver provided by the Eigen Library~\cite{Eigen} for the inversion of the MFEM and LDG matrices. In two dimensions we use the deal.ii finite element library \cite{BaHaKan-ACM07,BaHeHeKaKrMaTu-Preprint15} and the UMFPACK~\cite{DaDu-SIAM97,DaDu-ACM99,Davis-ACM04-1,Davis-ACM04-2} direct solver for the inversion of the MFEM and LDG matrices.

\subsection{Benchmark on spatial discritizations}

\begin{figure}[!htb]
\centering
\subfloat[$L^{2}$ error plots for steady-state electron densities in 1D.]{
\includegraphics[scale=0.3]{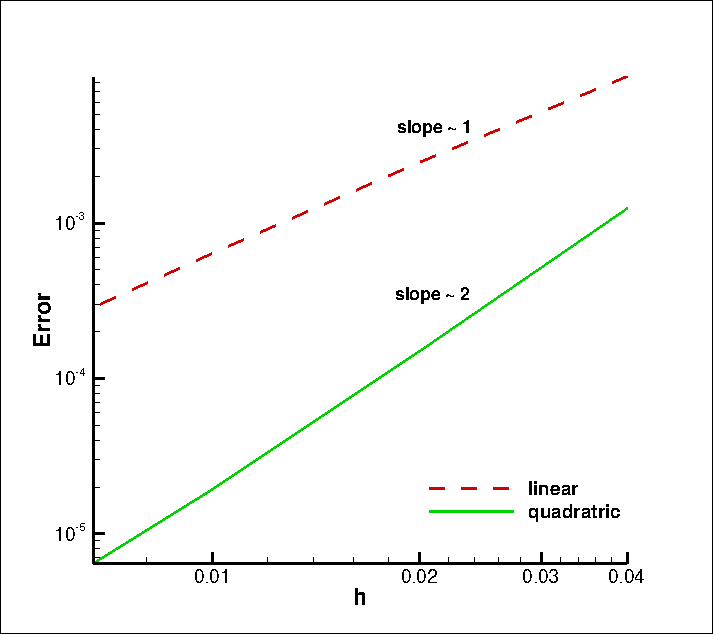}
}
\hspace{3mm}
\subfloat[$L^{2}$ error plots for $u$ in problem \eqref{eq:Interface_Errors_1} and \eqref{eq:Interface_Errors_2}.]{
\includegraphics[scale=0.3]{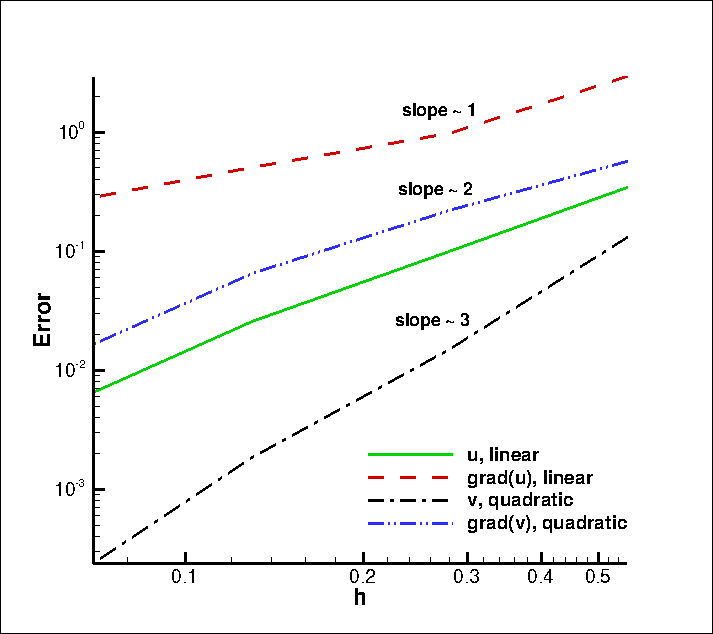}
}
\label{FIG:Schwarz Error}
\end{figure}

To benchmark the algorithms for the semiconductor-electrolyte interface simulation, we plot $L^{2}$ errors in the steady-state electron densities in Fig.~\ref{FIG:Schwarz Error} (a) in one dimension.  In our numerical experiments we did not notice any difference between the accuracy of the steady state solutions computed by any of the four stepping algorithms.  The results in Fig.~\ref{FIG:Schwarz Error} (a) we obtained using the TsPS-IMEXEX algorithm.  Note that since we do not have a manufactured analytical solution here, we use the numerical solution on an extremely fine mesh as the true solution. In Fig.~\ref{FIG:Schwarz Error} (a) we show the results for linear and quadratic LDG approximations. The results show almost optimal convergence rates. The reason for plotting only the $L^{2}$ error of the steady-state electron density and not the steady-state current distributions is that the steady-state currents in one dimension are constants and therefore the $L^{2}$ error plots would not reveal the convergence rate.

In the two dimensional code we tested all the individual solvers and coupling procedures that will be used in simulations of the semiconductor-electrolyte interface.  We present here a benchmark on the underlying LDG discretization for reactive-interfaces on a simplified problem.  We take $u: \, (0,T] \times \Omega_{S} \rightarrow \mathbb{R}$ and $v: \, (0,T] \times \Omega_{E}  \rightarrow \mathbb{R}$ and solve the coupled parabolic problems,
\begin{equation}
\begin{aligned}
u_{t} - \Delta u  \; &= \; f_{1} (\textbf{x},t) && \text{in} \quad  (0,T] \times \Omega_{S}, \\
\textbf{n}_{\Sigma^{S}} \cdot \left( -\nabla u \right) \cdot  \; &= \; u(\textbf{x},t)v(\textbf{x},t) - I(\textbf{x},t) && \text{on} \quad  (0,T] \times \Sigma,  \\
u \; &= \; g_{1,D}(\textbf{x},t) && \text{on} \quad (0,T] \times \Gamma_{C},  \\
u \; &= \; u_{0}(\textbf{x}) && \text{in} \quad \{t=0\} \times \Omega_{S}, 
\end{aligned}
\label{eq:Interface_Errors_1}
\end{equation}

\noindent
and,
\begin{equation}
\begin{aligned}
v_{t} - \Delta v \; &= \; f_{2} (\textbf{x},t) && \text{in} \quad   (0,T] \times \Omega_{E} ,\\
 \textbf{n}_{\Sigma^{E}} \cdot \left(-\nabla v \right)  \; &= \; u(\textbf{x},t)v(\textbf{x},t) - I(\textbf{x},t) && \text{on} \quad  (0,T] \times \Sigma , \\
v \; &= \; g_{2,D}(\textbf{x},t) && \text{on} \quad  (0,T] \times\Gamma_{A}, \\
v \; &= \; v_{0}(\textbf{x}) && \text{in} \quad  \{t=0\} \times\Omega_{E}. 
\end{aligned}
\label{eq:Interface_Errors_2}
\end{equation}

\noindent
We take the domain to be $\Omega = [0,1] \times [0,1]$ with $\Omega_{S} = [0,1/2] \times [0,1]$ and $\Omega_{E} = [1/2,1] \times [0,1]$. Let $\bx=(x,y)$. The interface is $\Sigma = \{x  = 1/2 \} \times [0,1]$ and the boundaries are $\Gamma_{C} = \partial \Omega_{S}  \setminus \Sigma $ and $\Gamma_{A} = \partial \Omega_{E} \setminus \Sigma$. The manufactured solutions for this problem are,
\begin{equation}
u(x,y,t) 
\; = \; v(x,y,t) 
\; = \; e^{-t} + \cos(2 \pi x) + \cos(2  \pi y).
\end{equation}

\noindent
The corresponding right hand side functions are,
\begin{equation}
f_{1} (x,y,t)
 \; = \; f_{2} (x,y,t)
 \; = \; -e^{-t} + 4\pi^{2} \cos(2 \pi x) + 4\pi^{2} \cos(2 \pi y) + 2 \pi \sin(2\pi x) .
\end{equation}

\noindent
The Dirichlet boundary conditions $g_{1,D}(x,y,t)$ and $g_{2,D}(x,y,t)$ are taken to be $u(x,y,t)$ on $\Gamma_{C}$ and $v(x,y,t)$ on $\Gamma_{A}$ respectively.  The interface function is,
\begin{equation}
{I}(x,y,t) \; = \; \left(e^{-t} + \cos(\pi y) - 1\right)^{2}.
\end{equation}

\noindent
The initial conditions $u_{0}(x,y)$ and $v_{0}(x,y)$ are taken to be the $L^{2}$ projection of the solutions $u(x,y,0)$ and $v(x,y,0)$ onto the DG basis.  To perform time stepping we use a first-order PS-IMEXEX method and an end time of $T=1$. In order to obtain the underlying errors of the LDG method we take time step to be $\Delta t = h^{k+1}$ when using basis functions of order $k$. The results can be seen in Fig.~\ref{FIG:Schwarz Error} (b) and show that we obtain optimal convergence rates for the LDG method.  We note that $h$ will be the same value for both of the triangulations of $\Omega_{S}$ and $\Omega_{E}$, that is $h_{S} = h_{E} = h$.

\subsection{Performance of time stepping algorithms}

\begin{figure}
\centering
\includegraphics[scale=0.7]{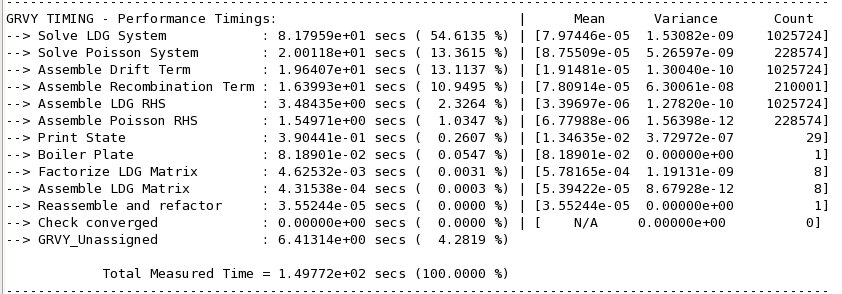}
\caption{\small{Run time performance results of key subroutines provided by the GRVY library.}}
\label{GRVY}
\end{figure}

We now look at the performance of the four time stepping algorithms. Our focus will be on the comparison between the four algorithms. We take the domain $\Omega=(-1, +1)$ where the interface $\Sigma$ is located at $x=0$. The semiconductor domain is $\Omega_S=(-1, 0)$ with $\Gamma_C$ being at $x=-1$.  The electrolyte domain is $\Omega_E=(0,1)$ with $\Gamma_A$ being at $x = +1$. We used 100 linear elements on a Dell Precision T1700 Workstation (i5-4590 Processor, Quad Core 3.30GHz). We use the GRVY library to monitor the computational time that is spent on each major subroutines in the algorithms. In Fig.~\ref{GRVY} we show a typical example of the information provided by the GRVY library.

\begin{table}[!ht]
\centering
\begin{tabular}{l| c c c c } 
\hline
\hline
& \small{AS-IMIMEX} & \small{AS-IMEXEX} & \small{PS-IMEXEX} & \small{TsPS-IMEXEX} \\ 
\hline
\small{Fact. LDG} & 2,322  & 2,421 &   $<1$  &  $<1$ \\
\small{Drift Term} & 9,133 & 48  & 92 &  51\\
\small{Recom. Term} & 307  & 326 &  293 & 286   \\
\small{Sol. LDG} &  265 & 265 &  258 & 141 \\
\small{Sol. MFEM} & 73 & 73 &  71 & 7  \\
\hline
Run Time & 12,498 & 3,529 &  766 &  518 \\
\hline
\end{tabular}
\caption{Time in seconds spent in the subroutines of the time stepping algorithms. Note: Total time includes time in other subroutines that were not recorded in this table.}
\label{TAB:Runtimes}
\end{table}

In Tab.~\ref{TAB:Runtimes}, we summarize the performance information of the four algorithms provided by GRVY. Besides the total run time of each algorithm, we showed, for each algorithm, the portion of the computational time spent on (i) assembling the drift term (first row), (ii) assembling the recombination term (second row), (iii) factorizing the LDG matrices (third row), (iv) inverting the LDG matrices (fourth row), and (v) inverting the MFEM matrices (fifth row).

The difference between the AS-IMIMEX algorithm and the AS-IMEXEX algorithm is that AS-IMIMEX uses implicit density values to assemble the drift term while AS-IMEXEX uses explicit density values to assemble the drift term. Using implicit density values to assemble the drift term is prohibitively expensive. In this case assembling the matrix corresponding to the drift term alone takes up nearly $73\%$ of the run time. We see that the AS-IMEXEX method is about $3.5$ times faster than the AS-IMIMEX method (in total CPU time). This is because using explicit density values to assemble the drift term we never compute the matrix that is used in AS-IMIMEX. Instead we assemble a right hand side vector corresponding to the drift term by performing local quadrature over each cell. It should be noted that the time step size needed to maintain stability for the AS-IMEXEX method was about one half of that needed for the AS-IMIMEX algorithm. This loss of stability is consistent with the choice of an explicit method instead of an implicit method.

We see from Tab.~\ref{TAB:Runtimes} that the main bottle neck in the AS-IMEXEX algorithm  is due to the factorization of the LDG matrices. If we treat all the densities on the interface using explicit density values then the LDG matrices will remain constant in time. All the matrices can then be factorized in the beginning of the simulation and only a linear solve must be performed at every time step. This is the reason that the PU-IMEXEX algorithm is about $4.5$ times faster (in total CPU time) than the AS-IMEXEX algorithm.

Finally, we observe from Tab.~\ref{TAB:Runtimes} that using a TsPU-IMEXEX method over PU-IMEXEX method results in a total speed up of approximately $1.5$. The faster run time is a result of the reduction in the total number of linear solves performed in TsPU-IMEX.  As noted previously, this method uses two time steps, one for the semiconductor systems and one for the electrolyte system.  In our simulations we found that time step for the electrolyte system, $\Delta t_e$, is almost $100$ bigger than the time step for the semiconductor system, $\Delta t_s$. Therefore the choice of $K\leq10$ mentioned in Section \ref{SUBSEC:TsPS-IMEXEX} could be relaxed and one might obtain even more savings from TsPU-IMEXEX.

\subsection{Device simulations under various model parameters}

We now present some numerical studies on the behavior of the model solutions under different model parameters. We are mainly interested in studying the impact of (i) device size, (ii) charge transfer rates on the interface,  (iii) variations in the doping profile  and (iv) interface geometry on the performance of the device. We use silicon as our semiconductor since it is the most commonly used semiconductor material for terrestrial photovoltaic devices. The choice of redox system is rather arbitrary as most simulation studies neglect the semiconductor electrode's interaction with it all together. Therefore we choose values that are computationally convenient, but are still representative of realistic electrolytes.  The material parameters values used in all simulations are recorded in Tab.~\ref{TAB:Sim_Parameters} and come from \cite{GrKe-PP95,Memming-Book15,Selberherr-Book01,Sze-Book81}. Other parameters such as the size~\cite{FoPrFeMa-EES12}, doping profile~\cite{FoPrFeMa-EES12} and charge transfer rates~\cite{Lewis-JPCB98} will be introduced when necessary. 
\begin{table}
\centering
\begin{tabular}{cl|cl}
\hline
\hline
Parameter & Value  & Parameter & Value  \\
\hline 
$\mu_{n}$ & $3.4911\times 10^{-3}$ & $\mu_{r}$ & $5.172 \times 10^{-4}$  \\
$\mu_{p}$ & $1.24128\times 10^{-3}$ & $\mu_{o}$ & $5.172 \times 10^{-4}$  \\
$\Phi_{\text{bi}}$ & $15.85$ & $\Phi^{\infty}$ & $0.0$  \\
$\lambda^{2}_{S} $ & $1.70215 \times 10^{-3}$   & $\lambda^{2}_{E} $ & $1.43038 \times 10^{-1}$   \\
$\tau_{n}$ & $5\times 10^{7}$ & - & -  \\
$\tau_{p}$ & $5\times 10^{7}$ & - & -\\
$\rho_{i}$ & $2.564 \times 10^{-7}$   & - & -\\
$\sigma_a$ & 17.4974  & - & -\\
$G_{0}$ & $1.2 \times 10^{-11}$  & - & -\\
\hline
\end{tabular}
\caption{Non-dimensionalized material parameters used in numerical simulations.}
\label{TAB:Sim_Parameters}
\end{table}

Besides the accuracy of the solutions to the mathematical model, we introduce two parameters to measure device performances: efficiency and fill factor. In order to define the efficiency and fill factor we introduce some necessary notations.  The short circuit current, labeled $J_{\text{SC}} $, is defined as the absolute value of current when $\Phi_{\text{app}} = 0$ and the device is illuminated. The open circuit potential, labeled $\Phi_{\text{OC}}$, is the value of the applied potential that produces zero current when the device is illuminated. We label the applied potential and absolute value of current that maximize the power output of the cell as $\Phi_{m}$ and $J_{m}$ respectively.  The efficiency $\eta_{\text{eff.}}$ and fill factor $\mathit{ff}$ of a solar cell are then defined as,
\begin{align}
\eta_{\text{eff.}} = \frac{\Phi_{m} J_{m}}{P_{\text{sun}}}, 
&&
\textit{ff} = \frac{\Phi_{m} J_{m}}{\Phi_{\text{OC}}  J_{\text{SC}}}
\end{align}
\noindent
where $P_{\text{sun}}$ is the power provided by the sun. We remark that $0 < \Phi_{m} < \Phi_{\text{OC}}$ and $0 < J_{m} <  J_{\text{SC}}$, so therefore $0 <  \mathit{ff} < 1$.

In the next four subsections, namely Sections~\ref{sec:size}-~\ref{sec:Schottky}, we focus on one-dimensional simulations with different model parameters. The last subsection, Section~\ref{sec:2D}, is devoted to two-dimensional simulations that focus on the impact of change of interface geometry on the performance of the device.

\subsubsection{Effects of device size}
\label{sec:size}

\begin{figure}[!ht]
\centering
\subfloat[Dark densities]{
\includegraphics[scale=0.16]{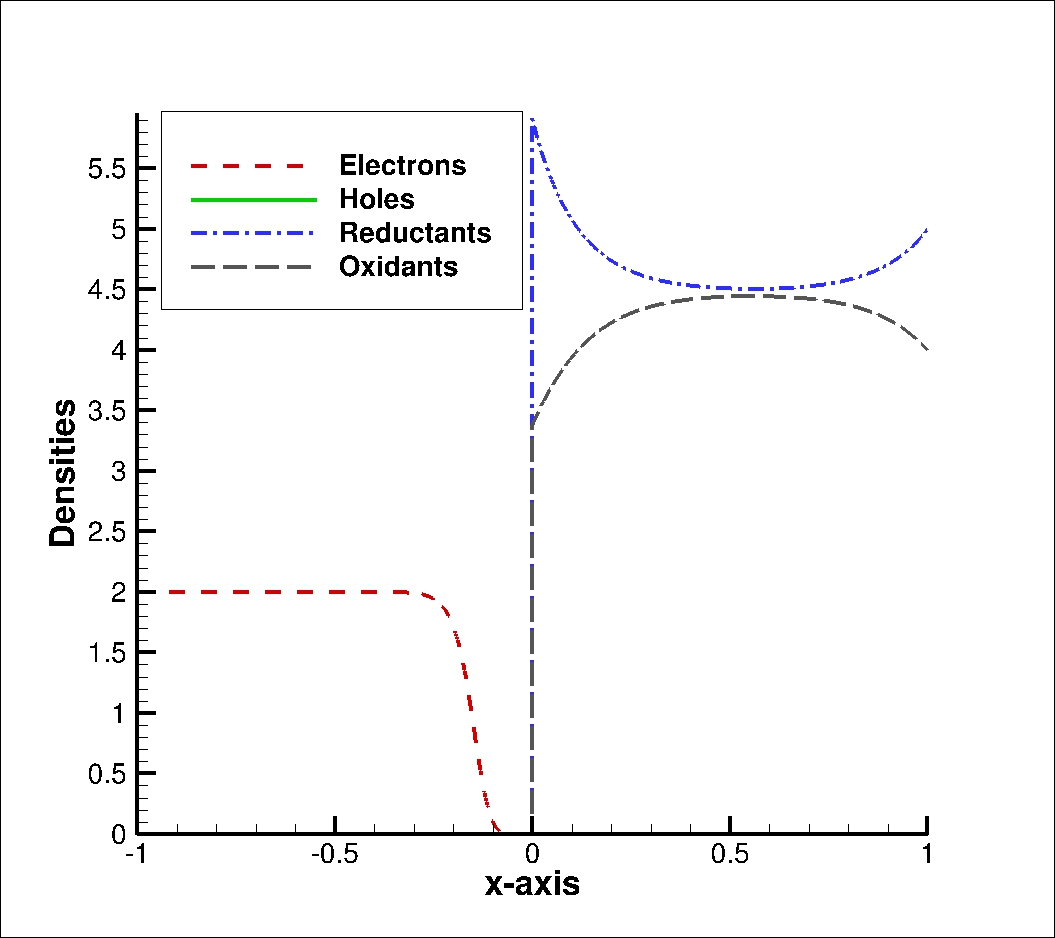} } \hskip 0.5cm
\subfloat[Illuminated densities]{
\includegraphics[scale=0.16]{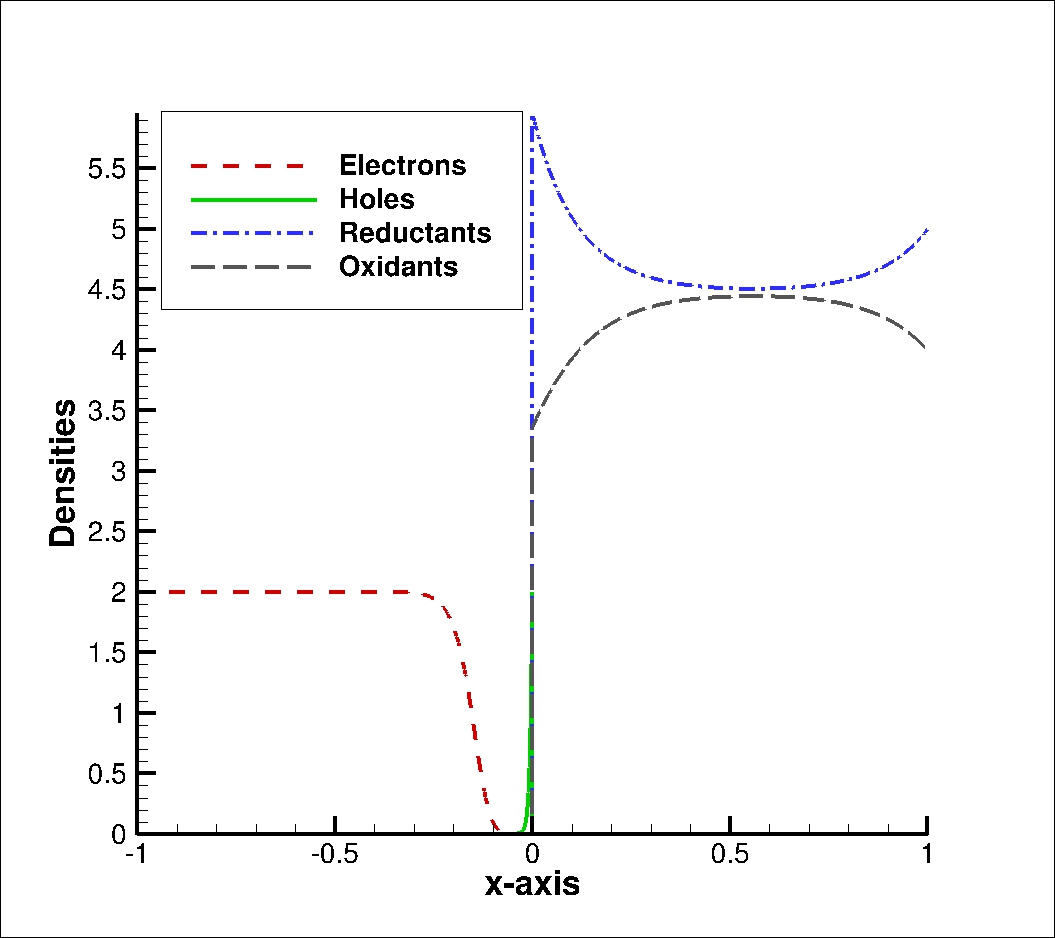} } \\
\subfloat[Dark current]{
\includegraphics[scale=0.16]{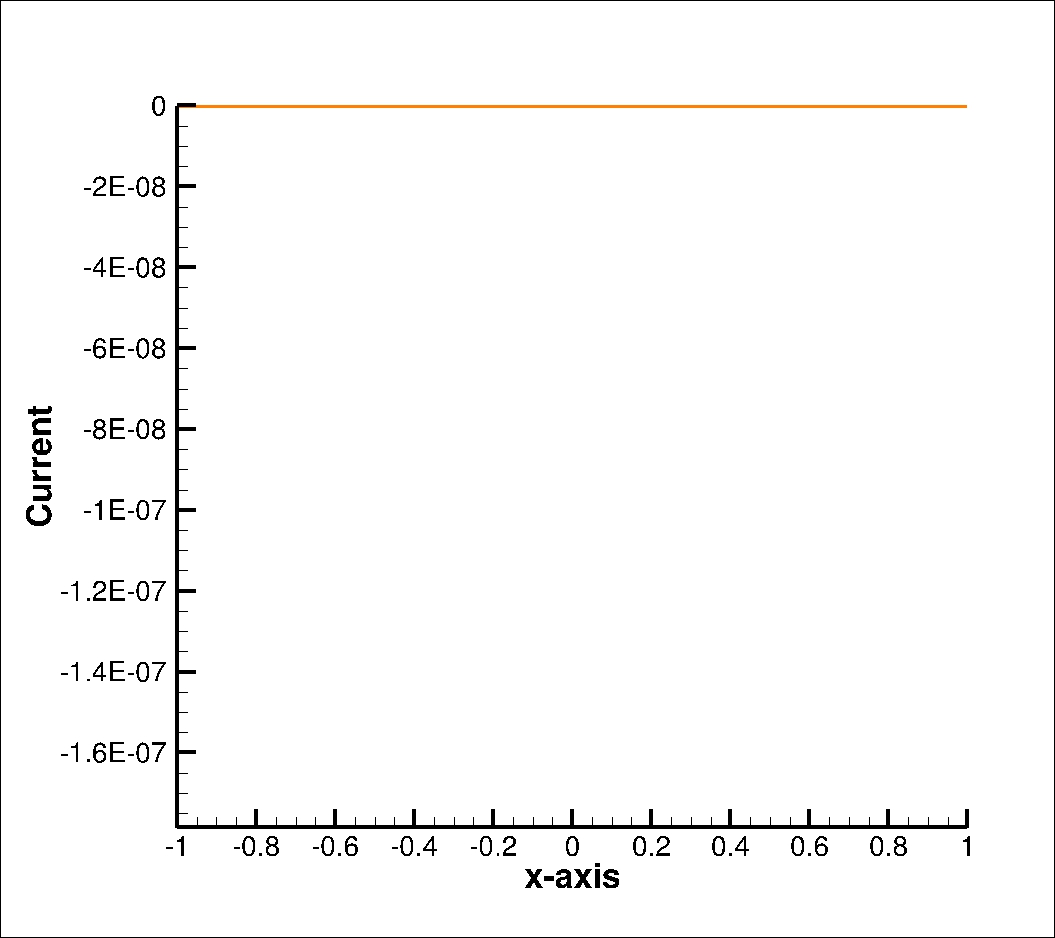} } \hskip 0.5cm
\subfloat[Illuminated current]{
\includegraphics[scale=0.16]{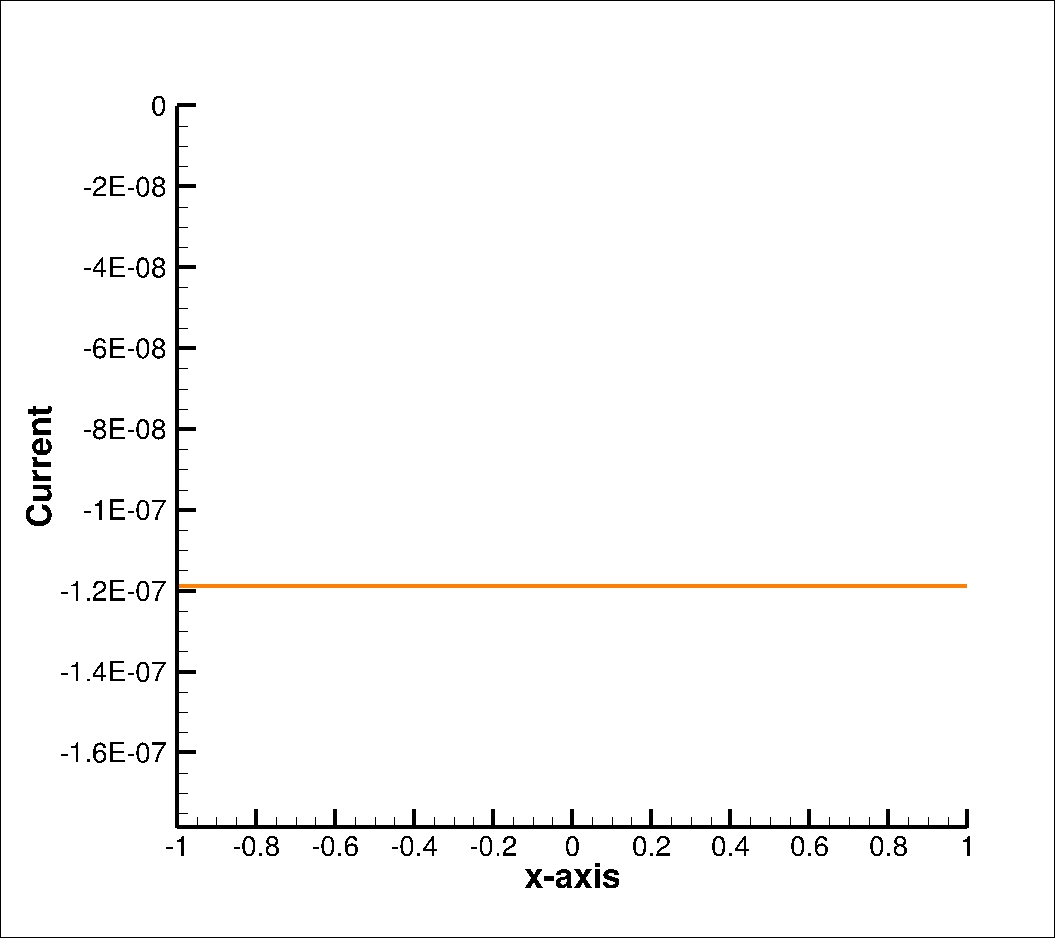} }
\caption{Steady state characteristics of device D-I.}
\label{FIG:MicroDensities}
\end{figure}

\begin{figure}[!ht]
\centering
\subfloat[Dark densities]{
\includegraphics[scale=0.16]{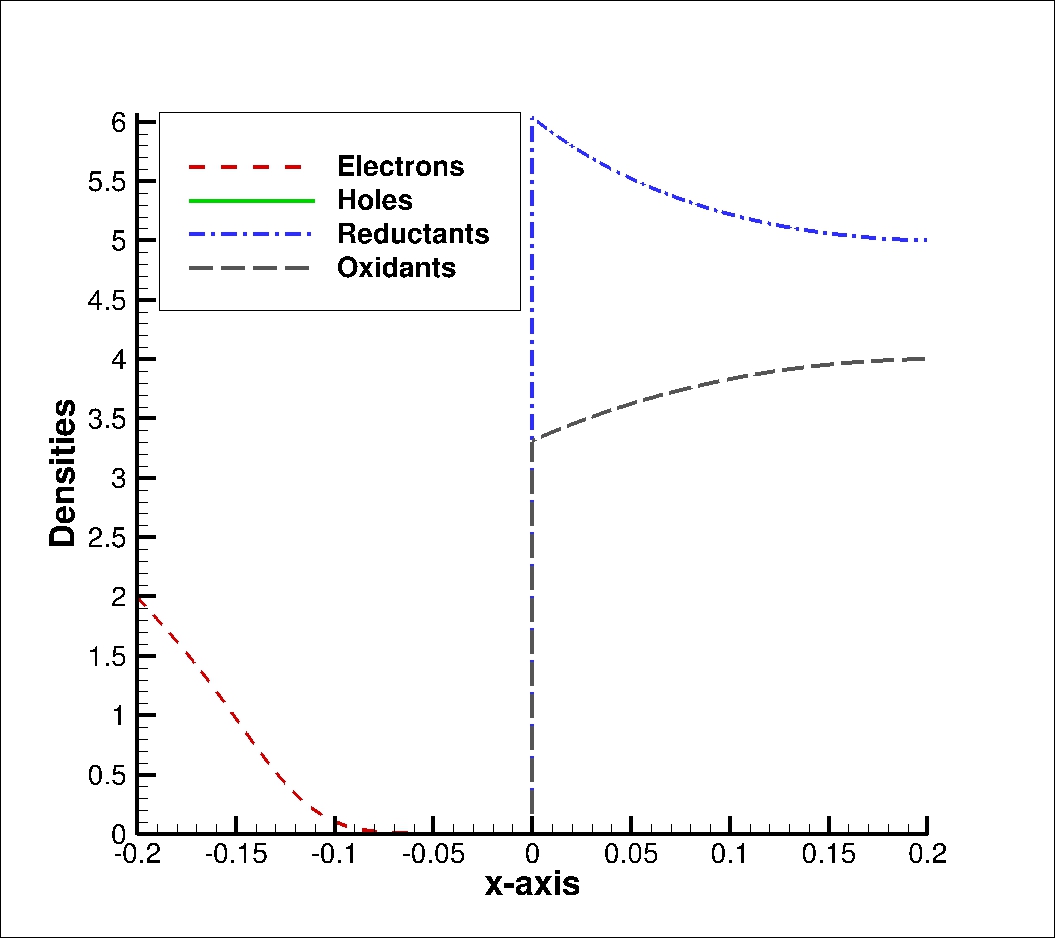} }\hskip 0.5cm
\subfloat[Illuminated densities]{
\includegraphics[scale=0.16]{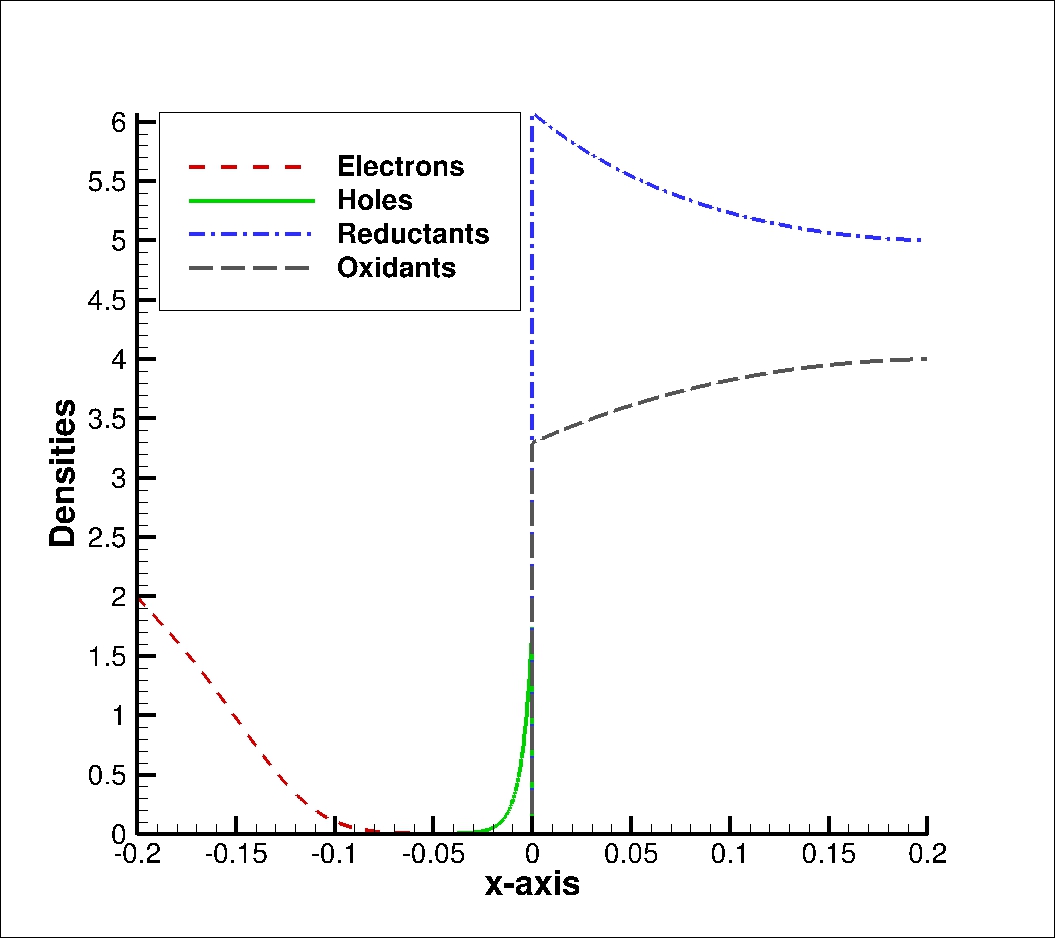} } \\
\subfloat[Dark current]{
\includegraphics[scale=0.16]{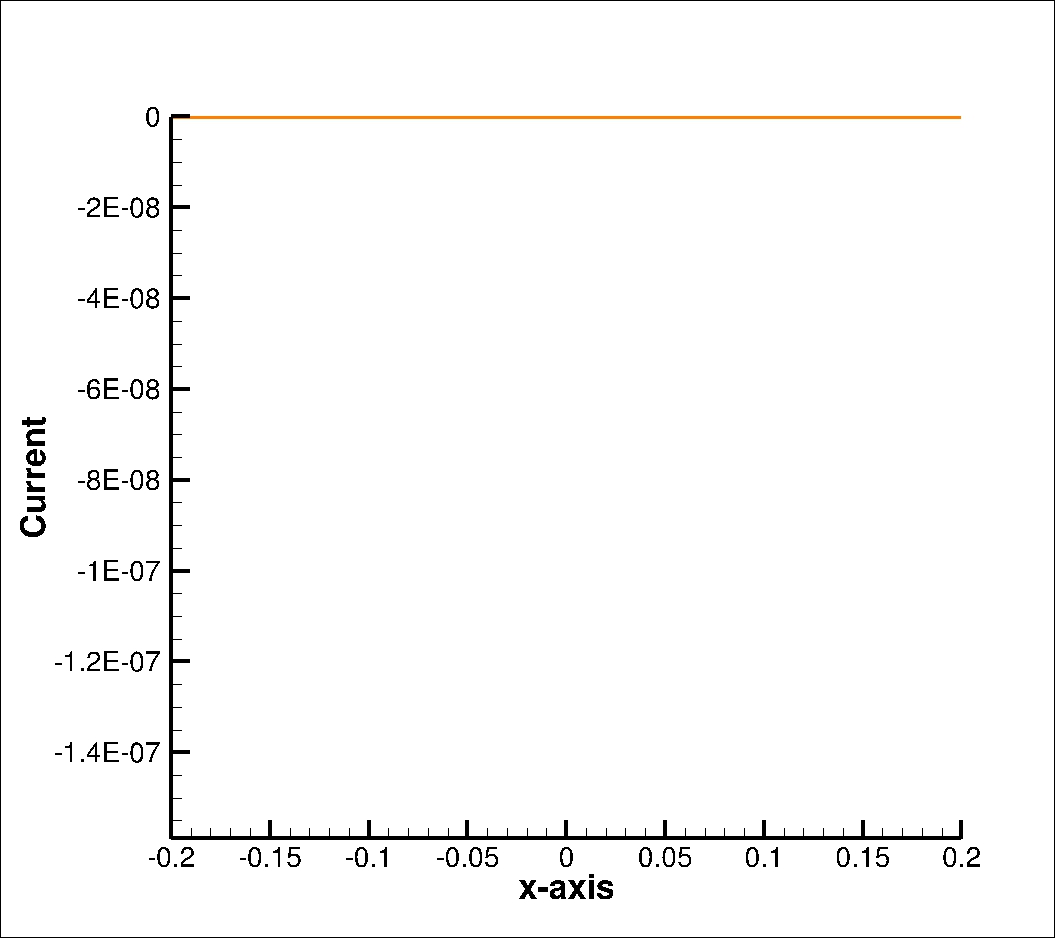} } \hskip 0.5cm
\subfloat[Illuminated current]{
\includegraphics[scale=0.16]{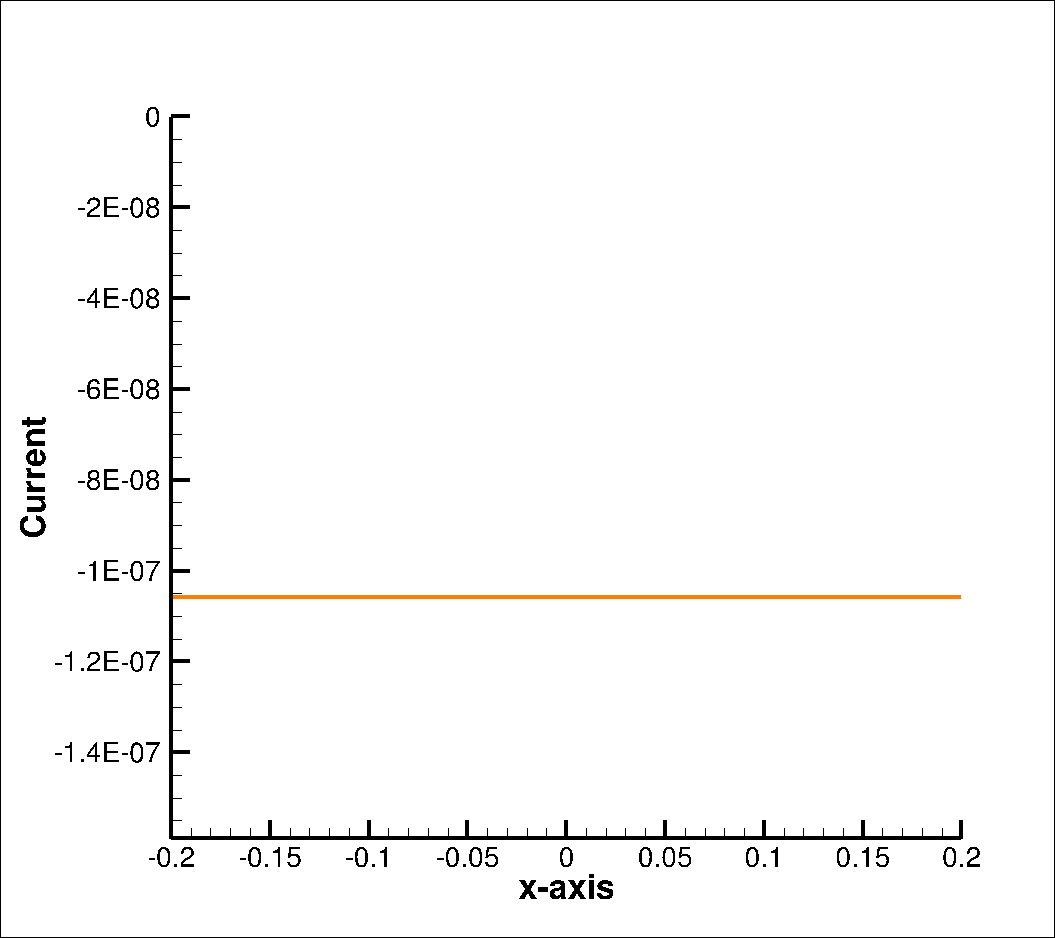} }
\caption{Steady state characteristics of device D-II.}
\label{FIG:NanoDensities}
\end{figure}

We consider two devices here: (i) D-I, $\Omega=(-1, 1)$ and (ii) D-II, $\Omega=(-0.2, 0.2)$. Besides their sizes, the other model parameters of the two devices are identical and are recorded in Tab.~\ref{TAB:VisualizationMicron} (left table). The parameters are chosen to better visualize the effects of illumination. The steady state characteristics of the devices under dark and illuminated conditions are presented in Fig.~\ref{FIG:MicroDensities} (for D-I) and Fig.~\ref{FIG:NanoDensities} (for D-II).

\begin{table}[!ht]
\centering
\begin{tabular}{cl|cl}
\hline
\hline
\multicolumn{2}{c|}{In $\Omega_S$}  & \multicolumn{2}{|c}{In $\Omega_E$ and On $\Sigma$} \\
\hline
$\rho_{n}^{e}$ & $ 2 \times 10^{16} $ & $\rho_{r}^{\infty}$ &  $ 5 \times 10^{16} $  \\
$\rho_{p}^{e}$ & $0 $ &  $\rho_{o}^{\infty}$ &  $ 4 \times 10^{16} $  \\
$\Phi_{\text{app.}}$ & 0 & $k_{et} $ & $ 10^{-11}$  \\
  & - & $k_{ht} $ & $ 10^{-8}$ \\
\hline
\end{tabular}
\hskip 1cm
\begin{tabular}{cl|cr}
\hline
\hline
\multicolumn{2}{c|}{In $\Omega_S$}  & \multicolumn{2}{|c}{In $\Omega_E$ and On $\Sigma$} \\
\hline
$\rho_{n}^{e}$ & $ 2$ & $\rho_{r}^{\infty}$ &  $30$ \\ 
$\rho_{p}^{e}$ & $0 $ & $\rho_{o}^{\infty}$ &  $ 29 $ \\
-& - & $k_{et} $ & $ 10^{-11}$ \\
-& - & $k_{ht} $ & $ 10^{-6}$ \\
\hline
\end{tabular}
\caption{Left: Parameter values for devices D-I and D-II; Right: Parameter values for devices D-III and D-IV.}
\label{TAB:VisualizationMicron}
\end{table}

\begin{figure}[!ht]
\centering
\parbox[t]{0.48\textwidth}{\null
\centering
\includegraphics[scale=0.25]{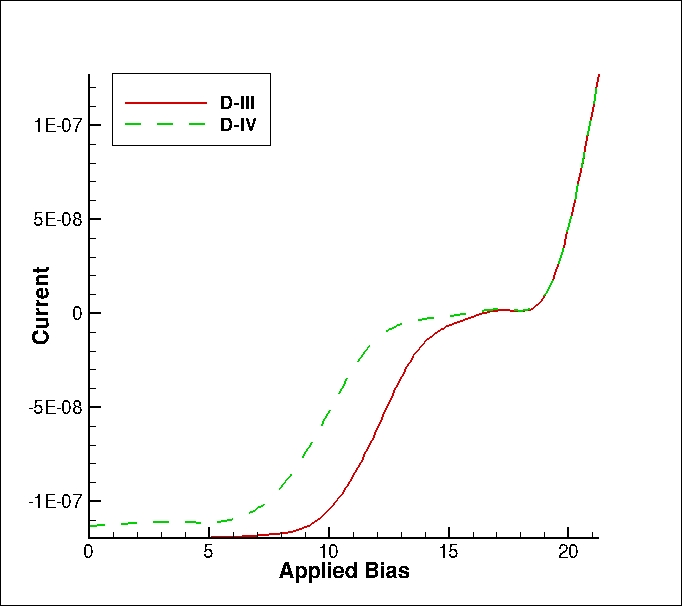}
\captionof{figure}{Illuminated current-voltage curves for devices D-III and D-IV.}
\label{FIG:Micro_Vs_Nano_IV} 
}
\hfill
\parbox[t]{0.48\textwidth}{\null
\centering
\vskip 3.9cm
\begin{tabular}{ccc}
\hline
\hline
Device & Efficiency & Fill Factor \\
\hline
D-III & $4.3\%$ & $0.524$ \\
D-IV & $3.1\%$ & $0.411$ \\
\hline
\end{tabular}
\captionof{table}[t]{Performance characteristics of devices D-III and D-IV.}
\label{TAB:Performance III IV}
}
\end{figure}

Comparing the results of the two devices it is evident that device characteristics vary greatly depending on its size. Indeed, the size of the device can have significant impact on the device's energy conversion rates. To see that, we perform simulations on two more devices: (iii) D-III and (iv) D-IV. D-III and D-IV have the same dimensions as D-I and D-II respectively. The model parameters, however, are very different; see the values listed in Tab.~\ref{TAB:VisualizationMicron}. Devices D-III and D-IV use higher concentrations of reductants and oxidants and also have faster minority transfer rates ($k_{ht}$). We apply a range of applied biases to the devices and record their steady state currents under illumination.  The results of these simulations are displayed in Fig.~\ref{FIG:Micro_Vs_Nano_IV} and the performance of the devices are summarized in Tab.~\ref{TAB:Performance III IV}. From Tab.~\ref{TAB:Performance III IV} it is observed that device D-IV has a lower efficiency and fill factor than D-III. Fig.~\ref{FIG:Micro_Vs_Nano_IV} shows that this is because both the short circuit current and the open circuit voltage have been reduced when using a smaller devices.

\subsubsection{Impact of interfacial charge transfer rates}

\begin{figure}[!ht]
\centering
\includegraphics[scale=0.18]{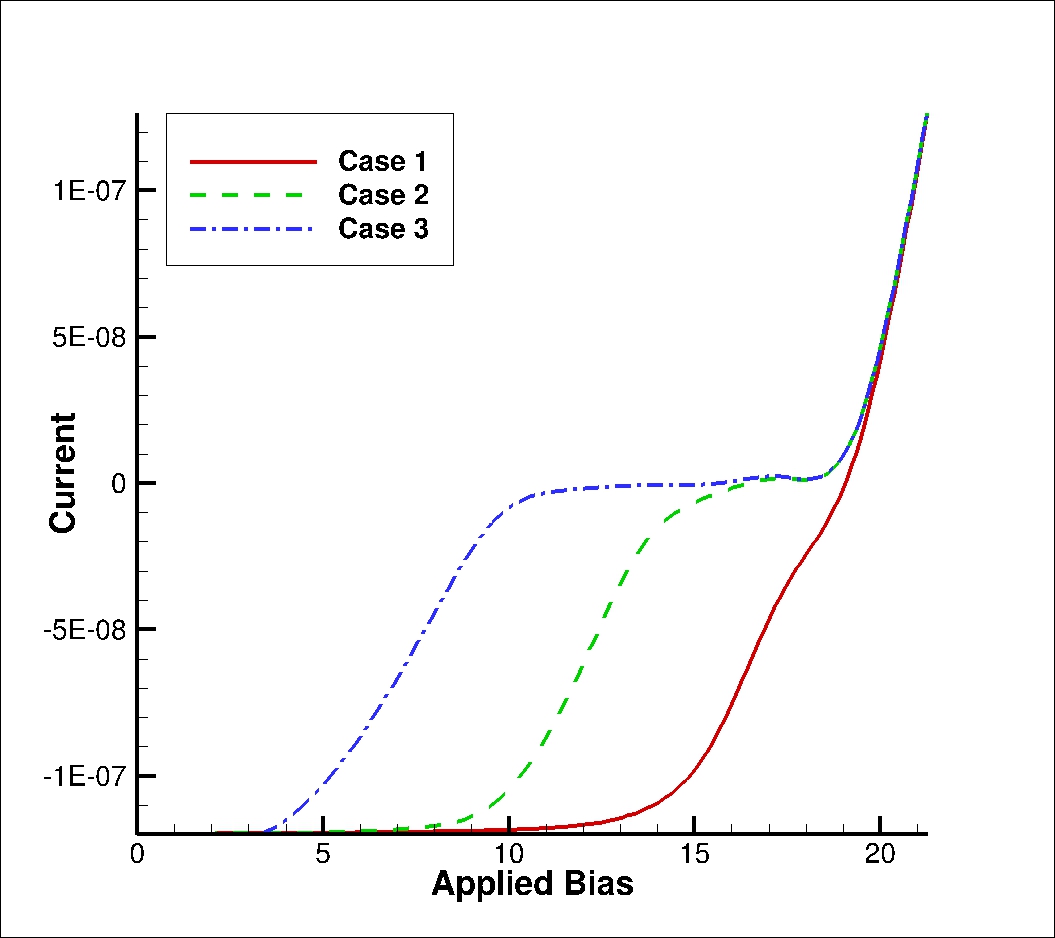} \hskip 0.5cm
\includegraphics[scale=0.278]{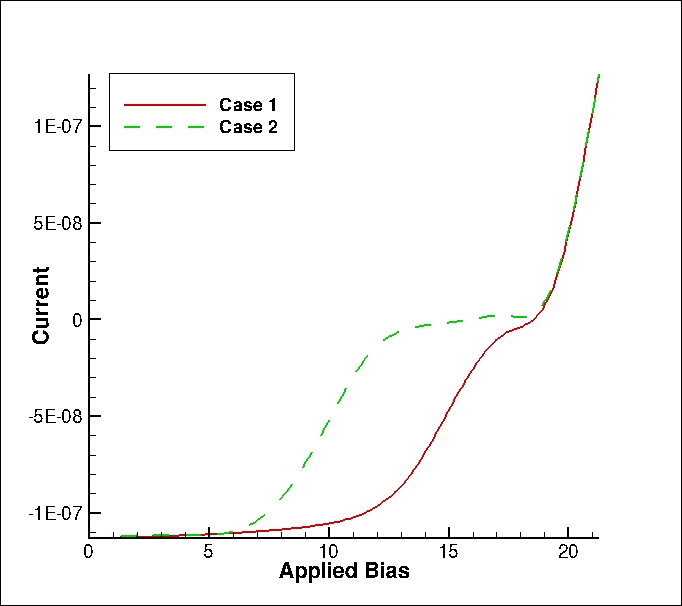}  
\caption{Illuminated current-voltage curves for devices D-III (left) and D-IV (right) different interfacial minority charge transfer rates.}
\label{FIG:Transfer-Rates}
\end{figure}

The charge transfer dynamics across the semiconductor-electrolyte interface is controlled mainly by the electron and hole transfer rates $k_{et}$ and $k_{ht}$. There has been tremendous progress in recent years on the theoretical calculation and experimental measurement of these rate constants~\cite{Lewis-JPCB98}. Here we study numerically the impact of these parameters on the performance of the PEC cell. 

In Fig.~\ref{FIG:Transfer-Rates}, we show the current-voltage curves of devices D-III and D-IV under different rates constants. For simplicity, we fix the electron transfer rate $k_{et}$, varying only the hole transfer rate $k_{ht}$ to see its effect. The values of $k_{ht}$ and the corresponding performance indicators are summarized in Tab.~\ref{TAB:Transfer-Rates}. It is clear from the data that the value of the minority transfer rate greatly effects the efficiency and fill factor of a PEC solar cell. Specifically it is observed that faster minority rates increase both the solar cell efficiency and fill factor. While the values of $k_{ht}$ that we considered here might be relatively high, they are still within the upper limit of its value given in ~\cite{Lewis-JPCB98}. In fact, the large amount of efforts in practical studies focus on how to select the combination of the semiconductor material and the electrolyte so that the transfer rate is high.
\begin{table}[!ht]
\centering
\begin{tabular}{cccc}
\hline
\hline
Case & $k_{ht}$ &  Efficiency & Fill Factor \\
\hline
1 & $10^{-4}$ & $6.3\%$ & $0.671$ \\
2 & $10^{-6}$ & $4.3\%$ & $0.524$ \\
3 & $10^{-8} $ & $2.2\%$ & $0.280$ \\
\hline
\end{tabular}
\hskip 1cm
\begin{tabular}{cccc}
\hline
\hline
Case & $k_{ht}$ & Efficiency & Fill Factor \\
\hline
1 &$10^{-4} $ & $4.8\%$ & $0.561$ \\
2 & $10^{-6} $ & $3.1\%$ & $0.411$ \\
\hline
\end{tabular}
\caption{Performance of devices D-III (left) and D-IV (right) under different minority transfer rates.}
\label{TAB:Transfer-Rates}
\end{table}

\subsubsection{Impact of the doping profile}

We now study the impact of the doping profile on the performance of PEC cells. We focus on the doping profile of the majority carrier. We perform simulations on device D-V which has the same size as device D-I, that is, $\Omega=(-0.2, 0.2)$, but different model parameters that are listed in Tab.~\ref{TAB:D-V} (left). We consider four different majority doping profiles that are given as:
\[
\cN_D^1 =2.0,\ x\in[-0.2,0]; \qquad 
\cN_D^2 = \left\{ 
\begin{array}{ll}
10.0 &  x\in[-0.2,-0.07) \\
2.0  &  x\in[-0.07 , 0] 
\end{array} \right.
\]
\[
\cN_D^3 = \left\{ 
\begin{array}{ll}
10.0 &  x\in[-0.20, -0.13)  \\
2.0 & x\in[-0.13, 0]
\end{array} \right. ;
\qquad
\cN_D^4=\left\{ 
\begin{array}{ll}
20.0 &  x\in[-0.20, -0.13) \\
2.0  &  x\in[-0.13, 0]
\end{array} \right.
\]
The values of $\rho_n^e$ is then given as $\rho_n^e=\cN_D$ since the minority carrier doping profile is set as $\cN_A=0$.
\begin{table}[!ht]
\centering
\begin{tabular}{cl}
\hline
\hline
Parameter & Value\\
\hline
$\rho_{p}^{e}$ & $0 $  \\
$\rho_{r}^{\infty}$ &  $ 30 $  \\
$\rho_{o}^{\infty}$ &  $ 29$ \\
$k_{et} $ & $ 10^{-11}$  \\
$k_{ht} $ & $ 10^{-4}$  \\
\hline
\end{tabular}
\hskip 1cm
\begin{tabular}{ccc}
\hline
\hline
Majority Doping & Efficiency & Fill Factor \\
\hline
$\cN_D^1$ & $4.8\%$ & $0.561$ \\
$\cN_D^2$ & $4.9\%$ & $0.626$ \\
$\cN_D^3$ & $4.9\%$ & $0.617$ \\
$\cN_D^4$ & $4.9\%$ & $0.634$ \\
\hline
\end{tabular}
\caption{Left: model parameters for device D-V; Right: device performance of D-V with different majority carrier doping profiles.}
\label{TAB:D-V}
\end{table}

\begin{figure}[ht!]
\centering
\subfloat[I-V curves]{
\includegraphics[scale=0.18]{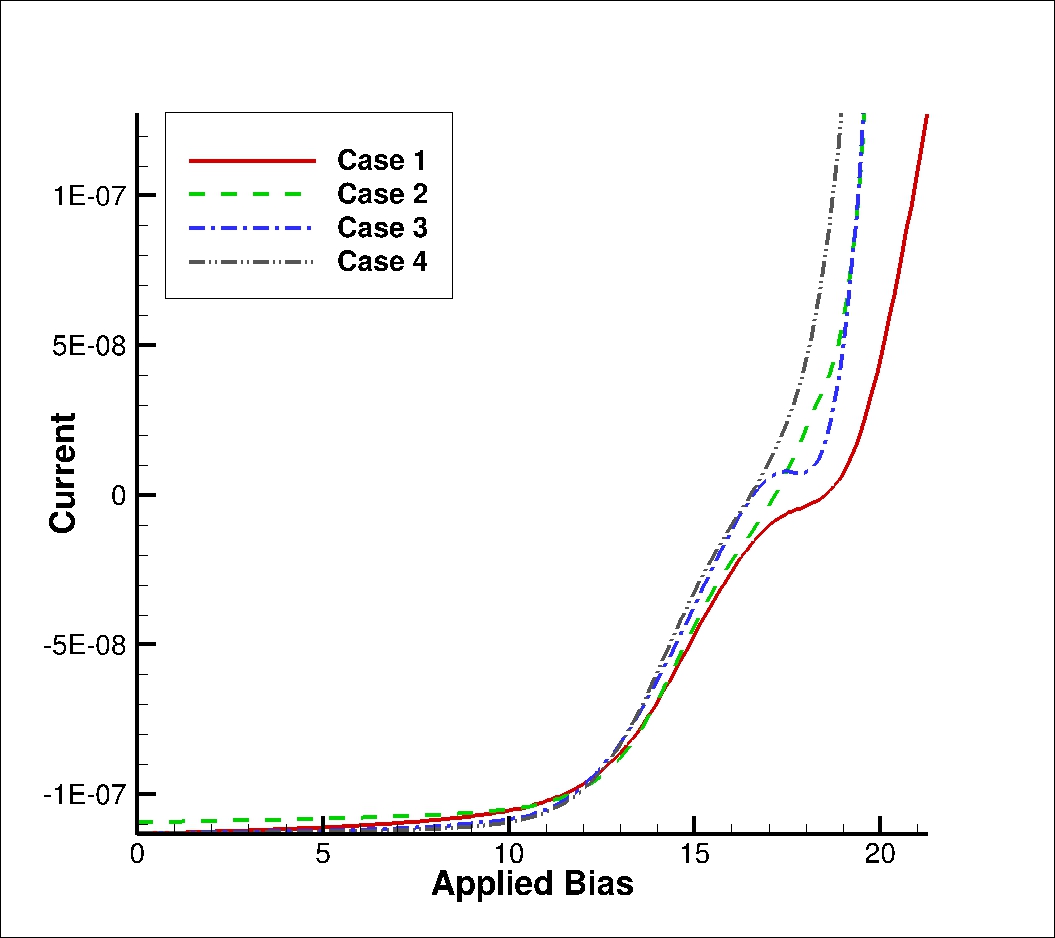}  }
\subfloat[$\Phi_{\text{app}} = 12$]{
\includegraphics[scale=0.18]{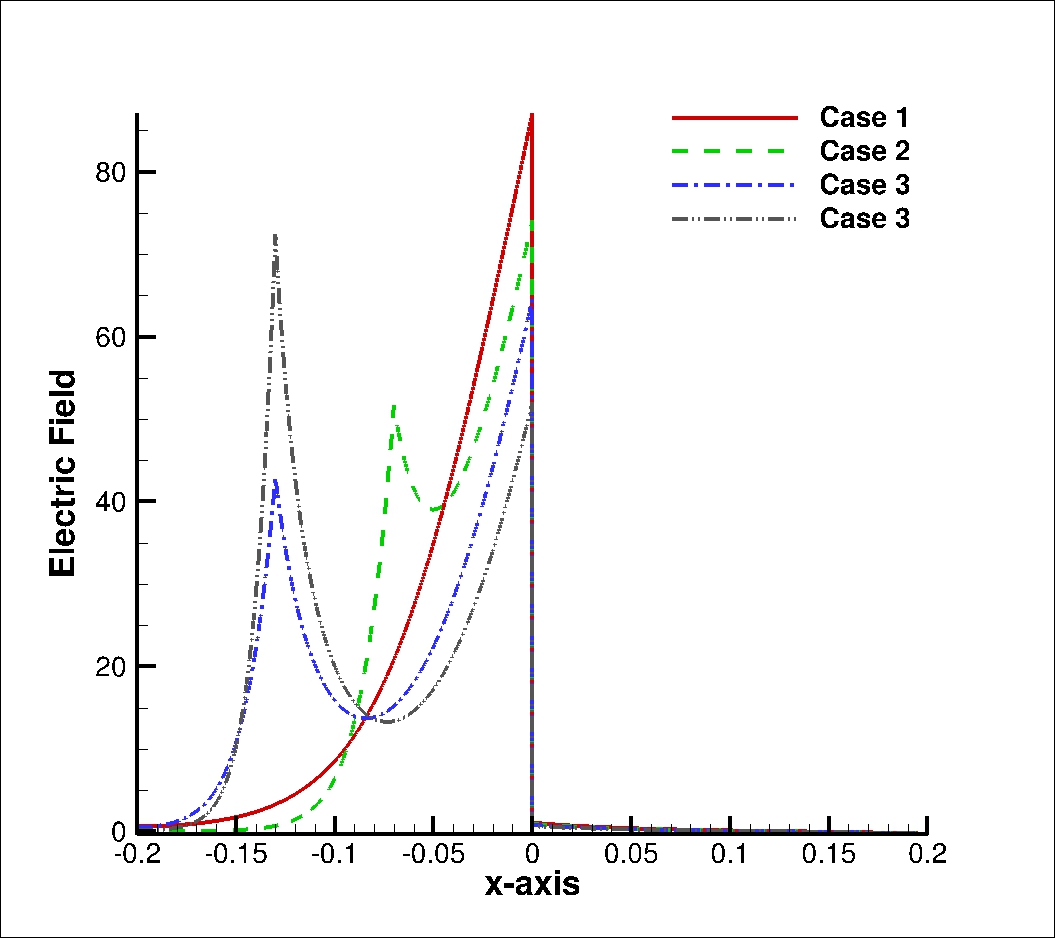} }
\\
\subfloat[$\Phi_{\text{app}} = 15.47$]{
\includegraphics[scale=0.18]{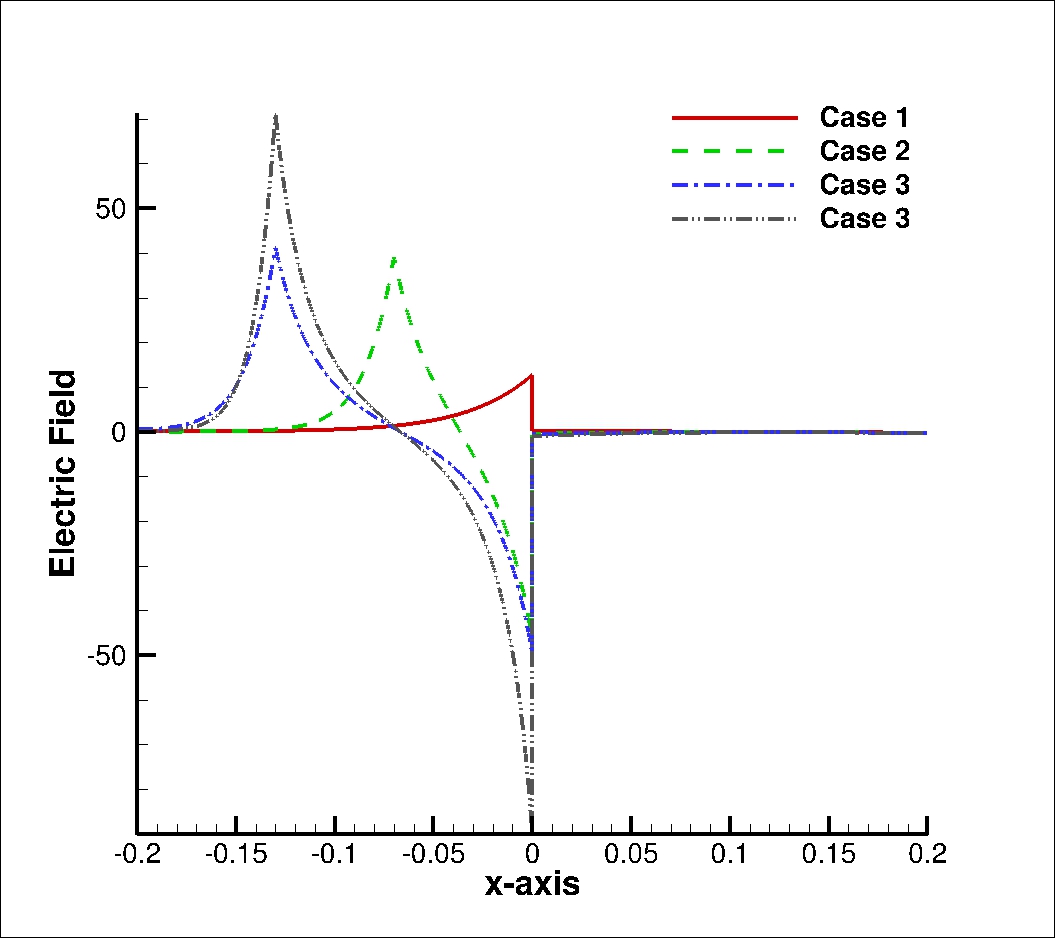} }
\subfloat[$\Phi_{\text{app}} = 19.34$]{
\includegraphics[scale=0.18]{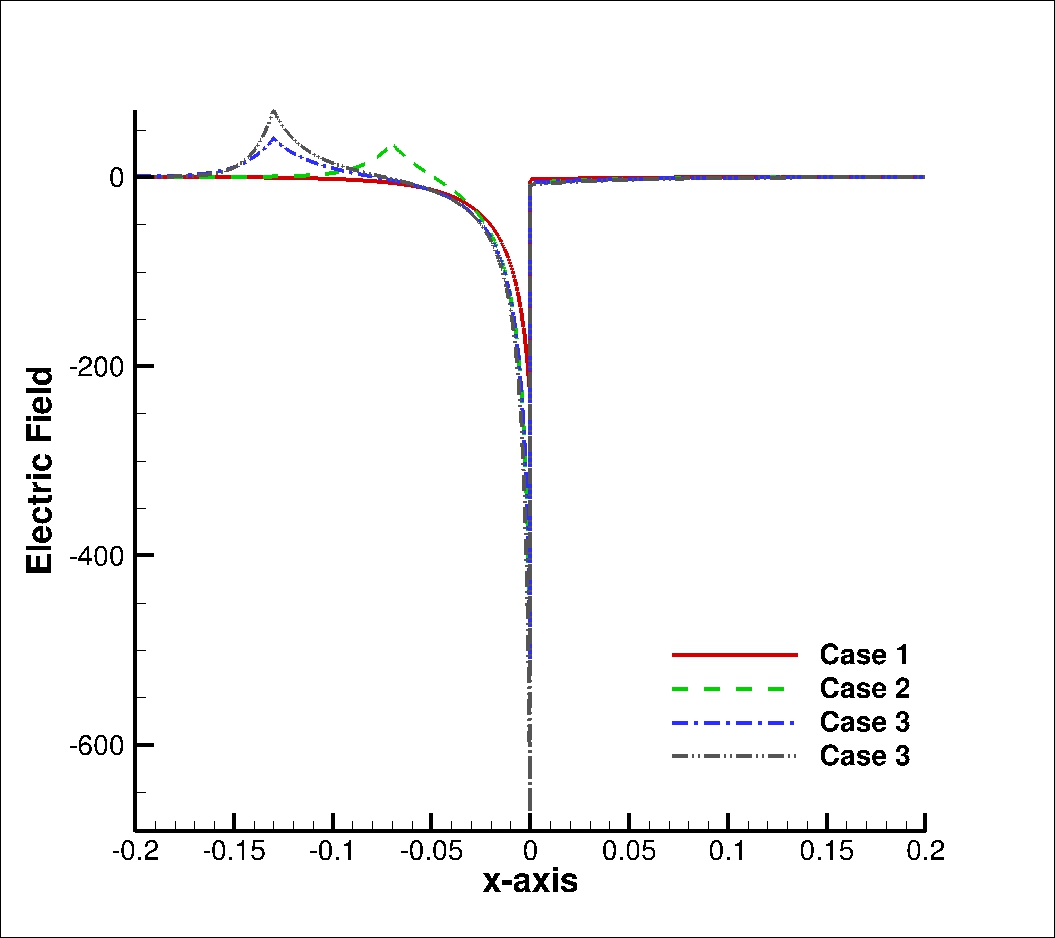} }
\caption{Device characteristics under different majority carrier doping profiles for D-V.}
\label{FIG:Nano_Step_Doping}
\end{figure}

The simulation results are displayed in Fig.~\ref{FIG:Nano_Step_Doping} and Tab.~\ref{TAB:D-V} (right). It seems that by adding variations to the doping profile we greatly improve the fill factor. It is evident in Fig.~\ref{FIG:Nano_Step_Doping} that the improvement of the fill factor is caused by a decrease in open circuit voltage. We see that a doping profile that has a thin, highly doped layer near the Ohmic contact has the best fill factor of all four cases. Interestingly, in Fig.~\ref{FIG:Nano_Step_Doping} the introduction of non-constant doping profiles causes the electric field to switch signs within the semiconductor domain for certain applied biases. This could be the reason for the reduction in the cells' open circuit voltage and the subsequent increase in the cell's fill factor. The phenomenon we observe here are consistent with what is known in the literature; see, for instance, the discussions in~\cite{FoPrFeMa-EES12,Nelson-Book03}.

\subsubsection{Comparison with the Schottky approximation}
\label{sec:Schottky}

In most previous simulation-based studies of PECs, the impact of the electrolyte solutions on the performance of the system has been modeled in a simpler manner. Essentially, the simulations are only performed on the semiconductor component. The impact of the electrolyte component comes in from the so-called ``Schottky boundary conditions'' on the interface~\cite{FoPrFeMa-EES12,MaRiSc-Book90,Nelson-Book03}. The Schottky boundary conditions on the interface are:
\begin{equation}
\begin{aligned}
\boldsymbol \bn_{\Sigma^{S}} \cdot (-\alpha_{n} \mu_{n} \rho_{n} \nabla \Phi 
\, - \, \mu_n \nabla \rho_{n} ) 
\; &= \;
v_{n} ( \rho_{n} \, - \, \rho_{n}^{e} ), \\
\bn_{\Sigma^{S}} \cdot (-\alpha_{p} \mu_{p} \rho_{p} \nabla \Phi 
\, - \, \mu_p \nabla \rho_{p} ) 
\; &= \;
v_{p} ( \rho_{p} \, - \, \rho_{p}^{e} ) ,
\end{aligned}
\end{equation}
where $v_{n}$ and $v_{p}$ are respectively the electron and hole recombination velocities. It is clear that when the densities of the reductants and oxidants do not change overall applied biases, these Schottky boundary conditions are simply the interface conditions~\eqref{EQ:Int Rho n} and~\eqref{EQ:Int Rho p} with $v_n=k_{et}\rho_o$ and $v_p=k_{ht}\rho_r$. If the densities of the reductants and oxidants change with applied biases, then the Schottky approximation fails to faithfully reflect the impact of the electrolyte system on the performance of the device.
\begin{table}[!ht]
\centering
\begin{tabular}{cl|cl}
\hline
\hline
\multicolumn{2}{c|}{In $\Omega_S$  and On $\Sigma$} & \multicolumn{2}{c}{In $\Omega_E$} \\
\hline
$\rho_{p}^{n}$ & $2 $ & $\rho_{r}^{\infty}$ &  $ 30$ \\
$\rho_{p}^{e}$ & $0 $ & $\rho_{o}^{\infty}$ &  $ 29$  \\
$k_{et} $ & $ 10^{-11}$ & - & - \\
$k_{ht} $ & $ 10^{-4}$ & - & - \\
$v_{n} $ & $ 3\times 10^{-9}$ &  - & -\\
$v_{p} $ & $ 2.9 \times 10^{-2}$ &  - & - \\
$\Phi_{\text{app.}} $ & $ 0$ & - & - \\
\hline
\end{tabular}
\hskip 1.5cm
\begin{tabular}{cl|cl}
\hline
\hline
\multicolumn{2}{c|}{In $\Omega_S$ and On $\Sigma$ } & \multicolumn{2}{c}{In $\Omega_E$} \\
\hline
$\rho_{p}^{n}$ & $2 $ & $\rho_{r}^{\infty}$ &  $ 5$ \\
$\rho_{p}^{e}$ & $0 $ & $\rho_{o}^{\infty}$ &  $ 4 $ \\
$k_{et} $ & $ 10^{-11}$ & - & - \\
$k_{ht} $ & $ 10^{-4}$ & - & - \\
$v_{n} $ & $ 3\times 10^{-9}$ &  - & -\\
$v_{p} $ & $ 2.9 \times 10^{-2}$ &  - & - \\
$\Phi_{\text{app.}} $ & $ 0$ & - & - \\
\hline
\end{tabular}
\caption{Model parameters for devices D-VI (left) and D-VII (right).}
\label{TAB:D-VI-VII}
\end{table}

\begin{figure}[!ht]
\centering
\subfloat[D-VI]{
\includegraphics[scale=0.18]{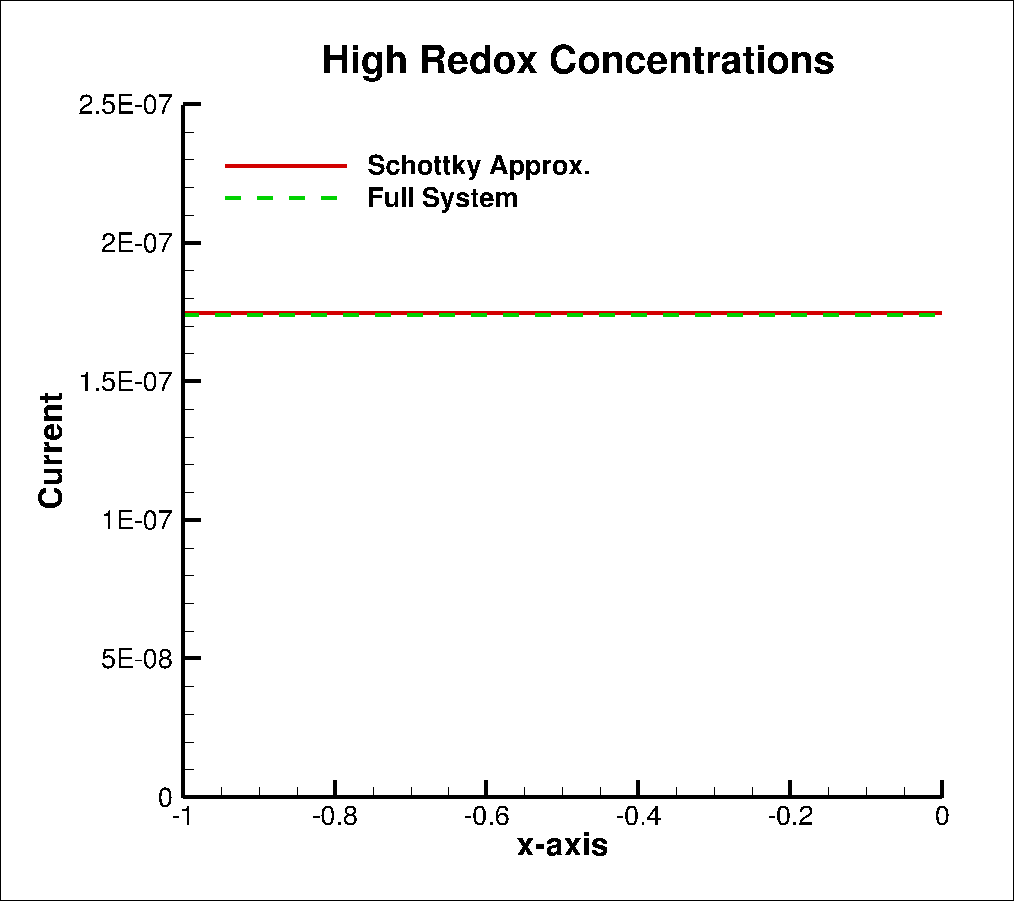}  }
\subfloat[D-VII]{
\includegraphics[scale=0.18]{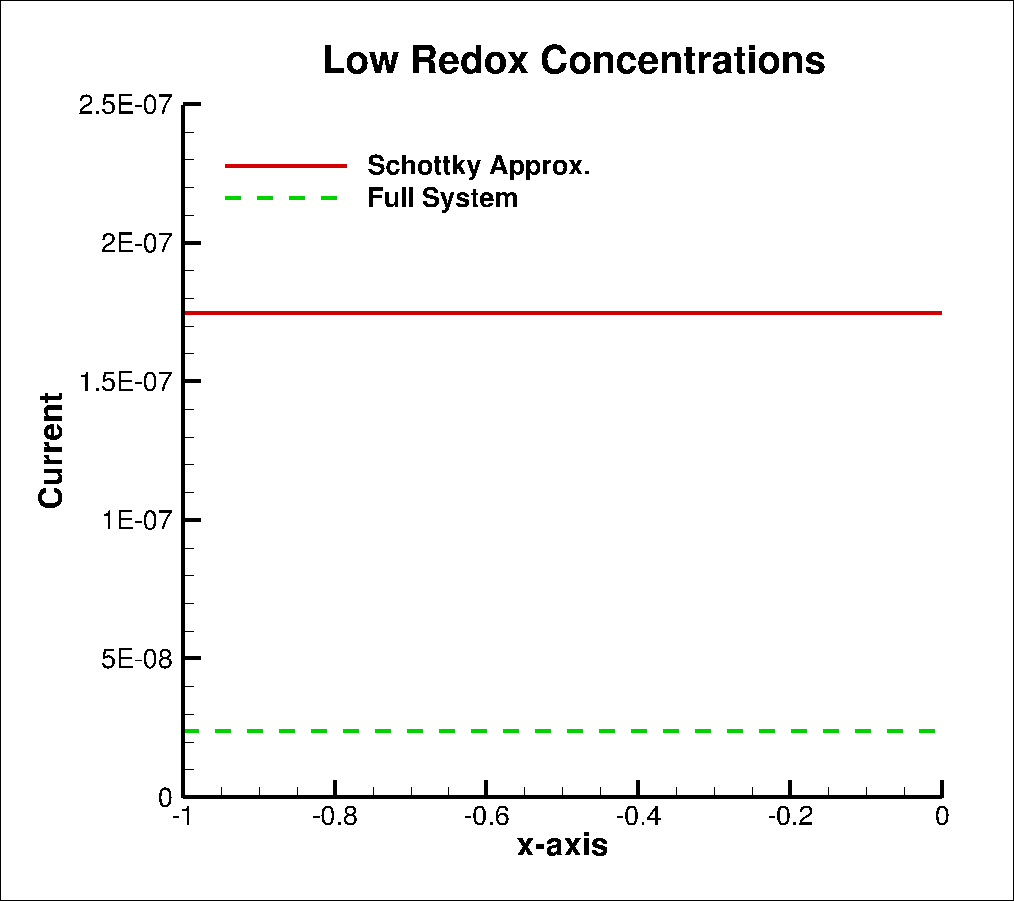} } \\
\caption{Comparison of of an illuminated device's current densities resulting from simulations of the full semiconductor-electrolyte system and its Schottky approximation.}
\label{FIG:Comparison_Carriers}
\end{figure}

We compare here the simulation results using our model in this paper with these using the Schottky approximation. We consider the comparison for a device with $\Omega=(-0.1, 0.1)$ with two different sets of parameters: D-VI, listed in Tab.~\ref{TAB:D-VI-VII} (left) and D-VII, listed in Tab.~\ref{TAB:D-VI-VII} (right). The main difference between the two devices are their redox pair concentrations. The comparison of current densities between simulations of the full system and the Schottky approximation for the two devices are displayed in Fig.~\ref{FIG:Comparison_Carriers}.

\begin{figure}[!ht]
\centering
\subfloat[I-V curves]{
\includegraphics[scale=0.32]{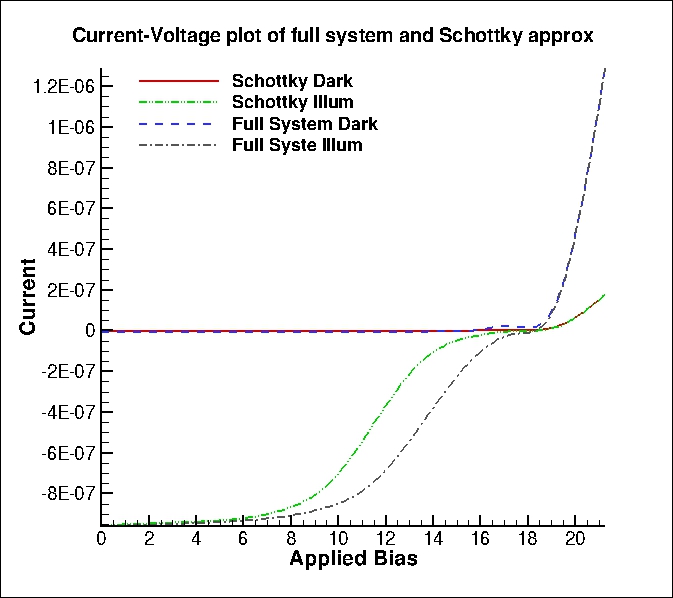} }
\subfloat[Zoomed in I-V curves]{
\includegraphics[scale=0.32]{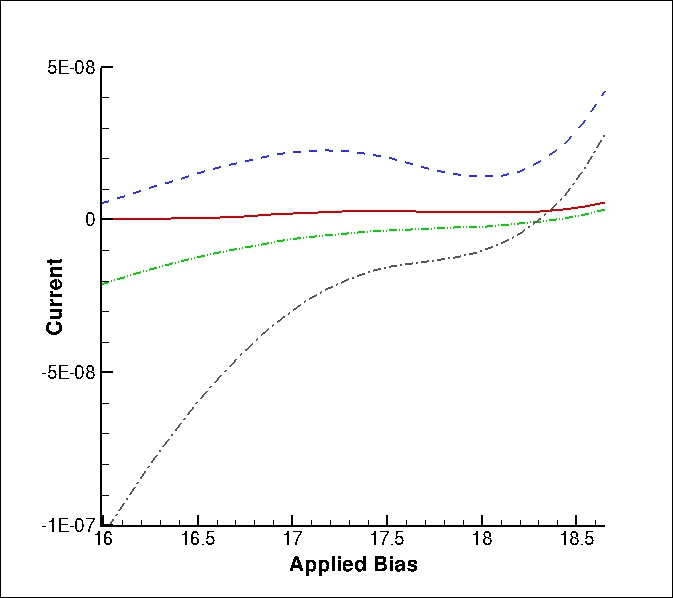} }
\caption{Comparison of current-voltage curves resulting from simulations of the full semiconductor-electrolyte system and the Schottky approximation with low redox concentrations.}
\label{FIG:Comparison_Schottky}
\end{figure}

It can be observed from Fig.~\ref{FIG:Comparison_Carriers} that the Schottky approximation produces results which are acceptable when the concentration of redox species are high compared to the densities of electrons and holes. However, when the redox concentrations are comparable with the density of electrons and holes the Schottky approximation produces results that deviate appreciably from those computed with the full system.  Specifically, in device D-VII, the current density computed with the full system is are much smaller than the current density computed with the Schottky approximation.  Deviations in current calculations can yield erroneous estimates of solar cell efficiency and fill factors. We demonstrate this as by plotting the current-voltage curves for both devices in Fig.~\ref{FIG:Comparison_Schottky}. The Schottky approximation over estimated the efficiency $(3.58\%$) and fill factor ($0.5$) compared to the simulation with the full reactive-interface conditions (efficiency $= 3.01\%$ and fill factor $= 0.42$) by over $0.5\%$ and 0.08 respectively. We can see that it would be much more accurate to use the full semiconductor-electrolyte systems in such instances.

\subsubsection{Two-dimensional simulations: impact of interface geometry}
\label{sec:2D}

\begin{figure}[ht!]
\centering
\subfloat[Nanowire shapes~\cite{FoPrFeMa-EES12}]
{\includegraphics[scale=0.65]{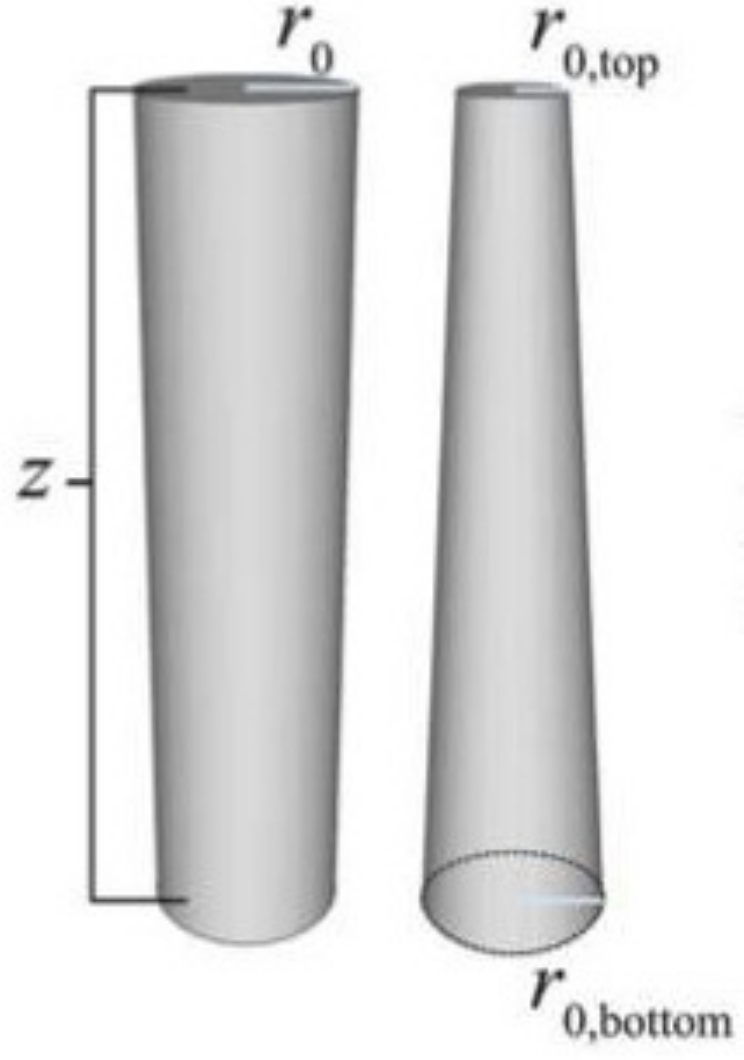}}
\hspace{5mm}
\subfloat[2D domain for nanowires]{
\includegraphics[scale=0.35]{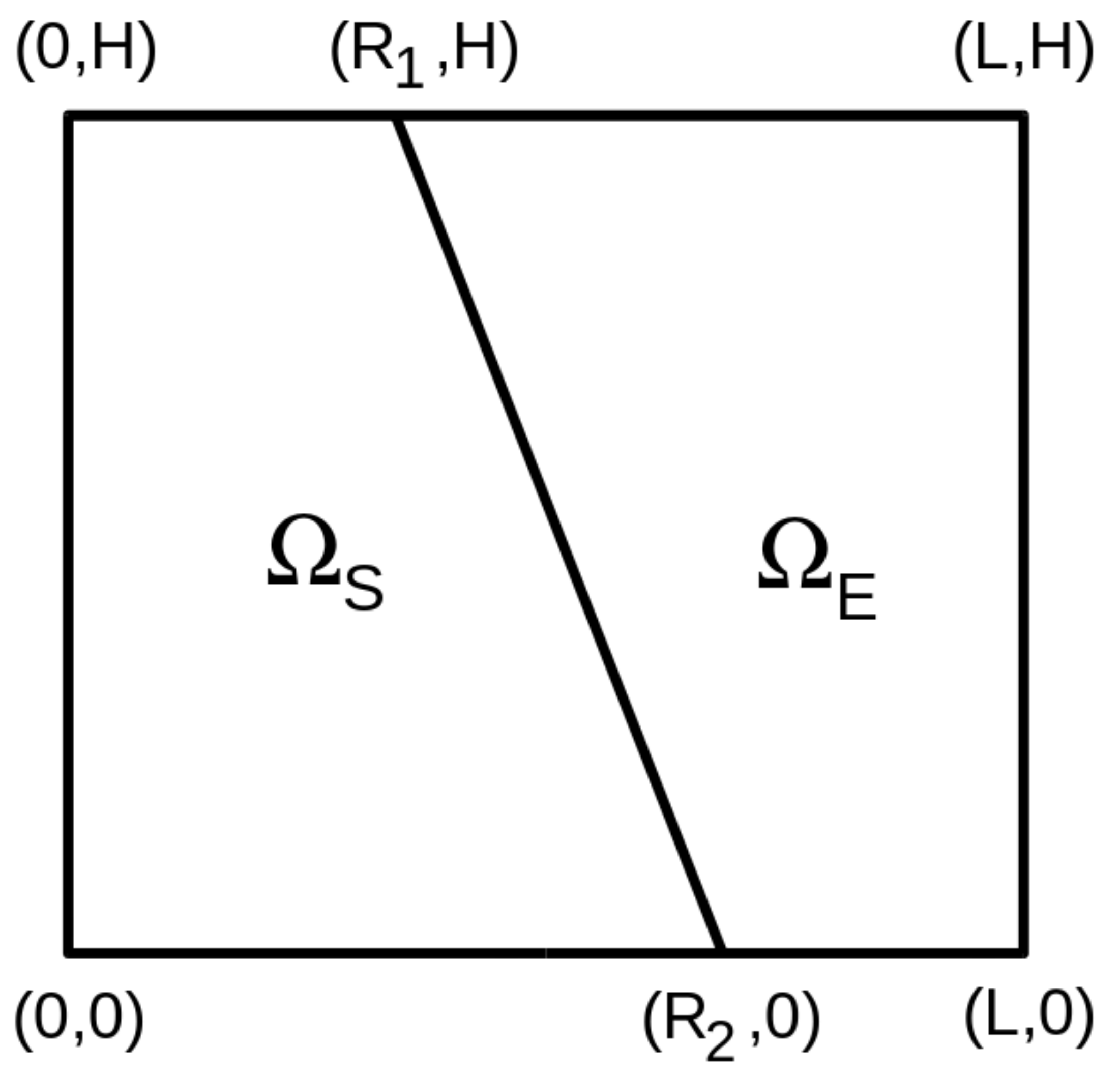}}
\caption{The domain for  two dimensional nanowire simulations.}
\label{fig:Domain}
\end{figure}

We now turn to the two-dimensional simulations to investigate the dependence of the numerical solution on geometry of the nanowire. The simulations for a two-dimensional semiconductor-electrolyte interface use the domain in Fig.~\ref{fig:Domain}  (b) to represent a cross section of the nanowires depicted in Fig.~\ref{fig:Domain} (a). Specifically we aim to investigate the effect the radii lengths $R_{1}$ and $R_{2}$ have on a PEC solar cell's characteristics.   The parameter values for the two-dimensional simulations are the same as those for device D-I and recorded in Tab.~\ref{TAB:VisualizationMicron} (left table).  These values were chosen to better visualize the dynamics of the charge densities under illumination.  

The top ($y = H$) and bottom ($y = 0$) of the electrolyte domain will always be insulated and the right boundary of the electrolye domain ($y=L$) will have bulk redox and bulk potential values.  That is we define the boundaries of the electrolyte to be,
\begin{equation}
\begin{aligned}
\Gamma_{A} \,  &= \,  \{ x \, = \, \text{L} \} \times [0,\text{H}],  \\
\Gamma_{N}^{e} \, &= \,  \left( [\text{R}_{2},\text{L}] \times \{y \, = \, 0\} \right) \cup \left( [\text{R}_{1},\text{L}] \times \{y \, = \, \text{H}\} \right) .
\end{aligned}
\end{equation}
\noindent
The boundary conditions of the semiconductor are the top and bottom of the domain are insulating and the left boundary of the domain is an Ohmic contact. That is,
\begin{equation}
\begin{aligned}
\Gamma_{C} \,  &= \,  \{ x \, = \, 0 \} \times [0,\text{H} ],  \\
\Gamma_{N}^{s} \,  &= \,  \left( [0,\text{R}_{2}] \times \{y \, = \, 0\} \right) \cup\left(  [0,\text{R}_{1}] \times \{y \, = \, \text{H}\} \right).
\end{aligned}
\end{equation}

\noindent
This choice of boundary conditions for semiconductor domain corresponds to the top of the device ($y = H$) being covered by glass. Light enters the device through the glass and travels straight downwards into the semiconductor.  Light will not be able to enter the bottom ($y=0$) or the left boundary ($x=0$) of the semiconductor domain.  We use the design of device D-VIII to simulate a cylindrical nanowire with dimensions that are summarized in Tab.~\ref{tab:2D_Flat} (left).  The other device we investigate, D-IX, is used to simulate a conic nanowire design and its parameter values are summarized in Tab.~\ref{tab:2D_Flat} (right).

\begin{figure}[!ht]
\centering
\subfloat[$\Phi_{\text{app}} \, = \, 0.0 $ ]{
\includegraphics[scale=0.3]{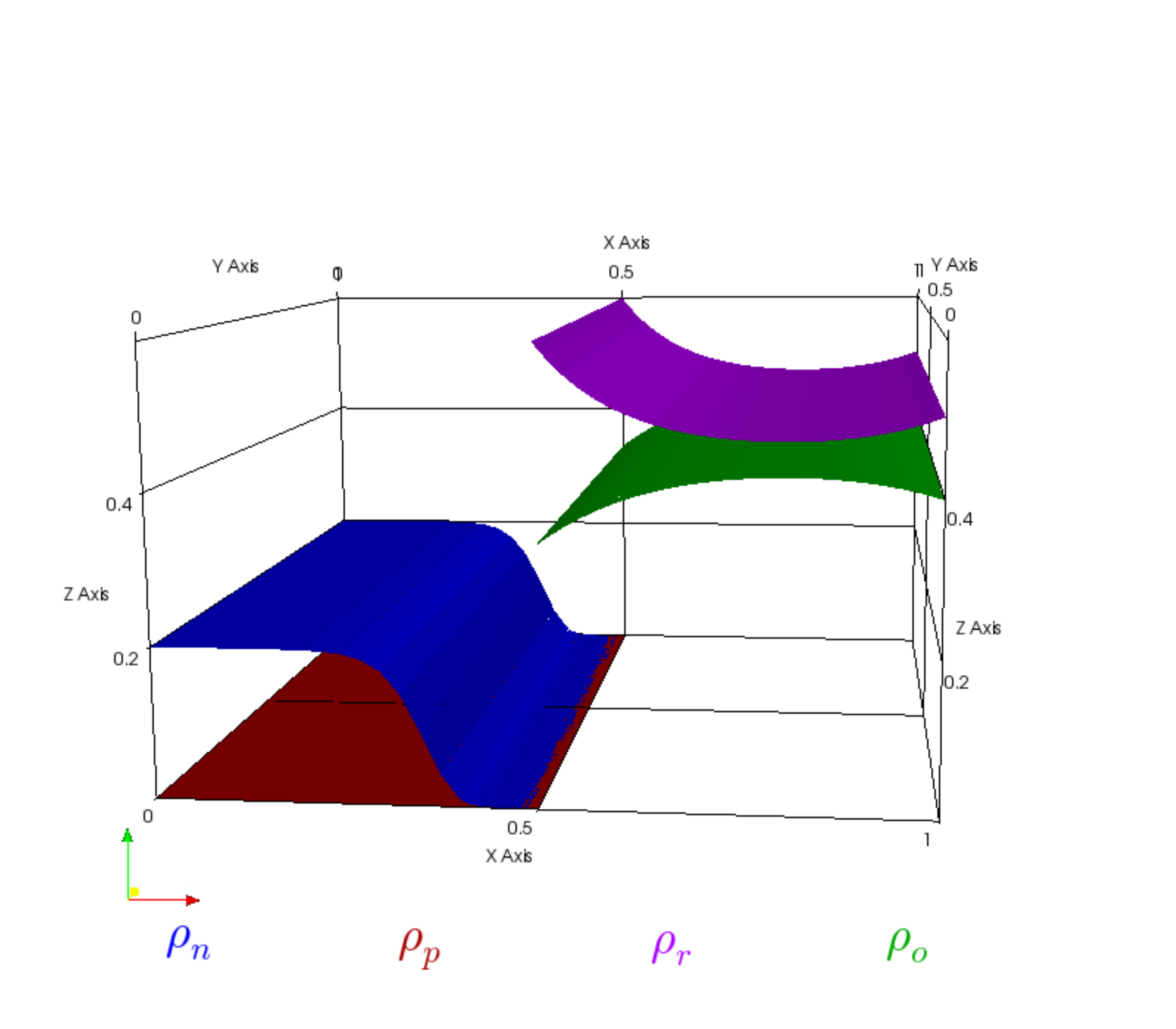} } 
\subfloat[$\Phi_{\text{app}} \, = \, 0.5 $]{
\includegraphics[scale=0.3]{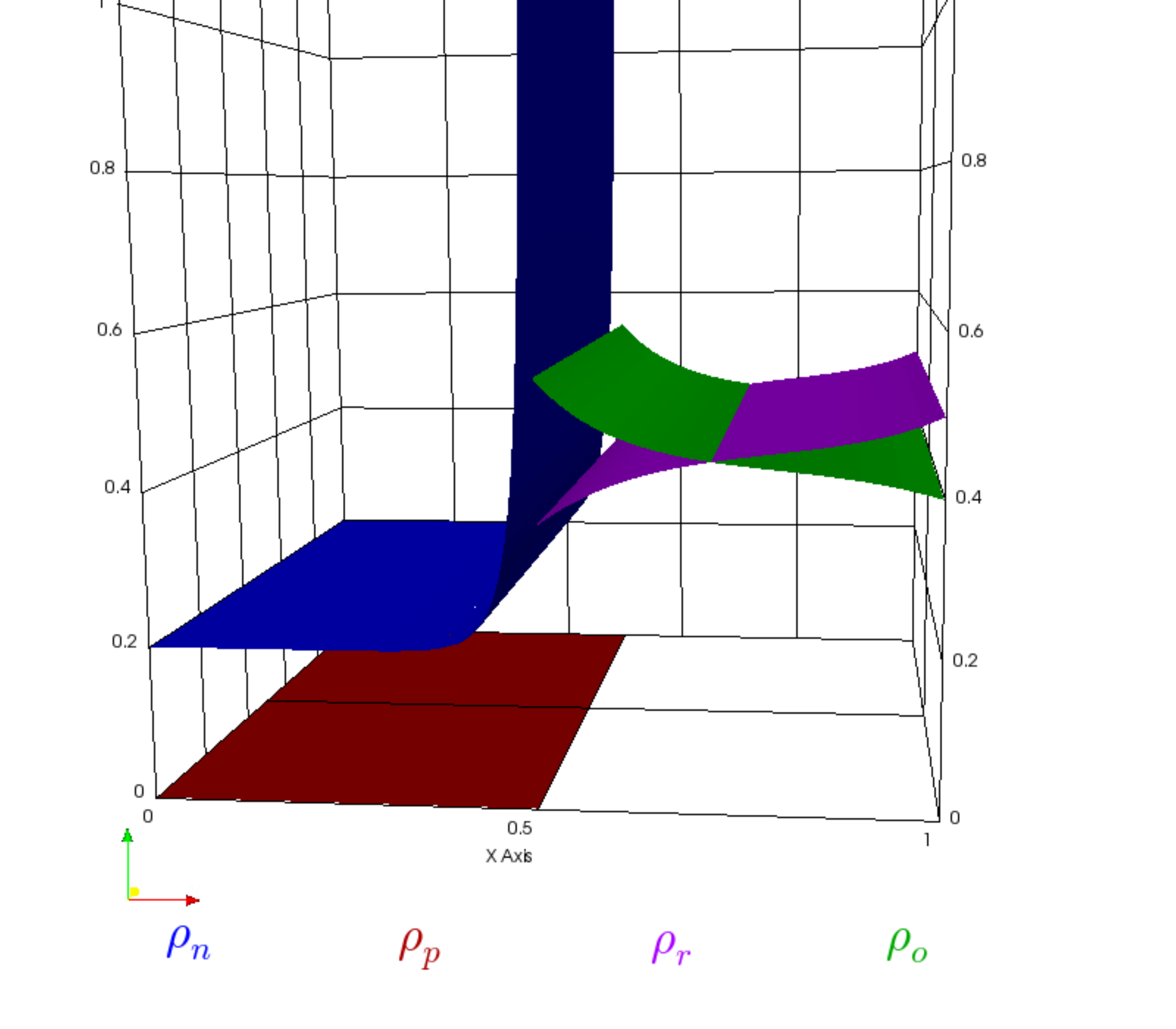} }
\caption{Densities of charge carriers in D-VIII. Note: Densities are in $10^{1} $.}
\label{fig:2D_flat_densities}
\end{figure}

\begin{figure}[!ht]
\centering
\subfloat[$\textbf{J}_{n}$]{
\includegraphics[scale=0.2]{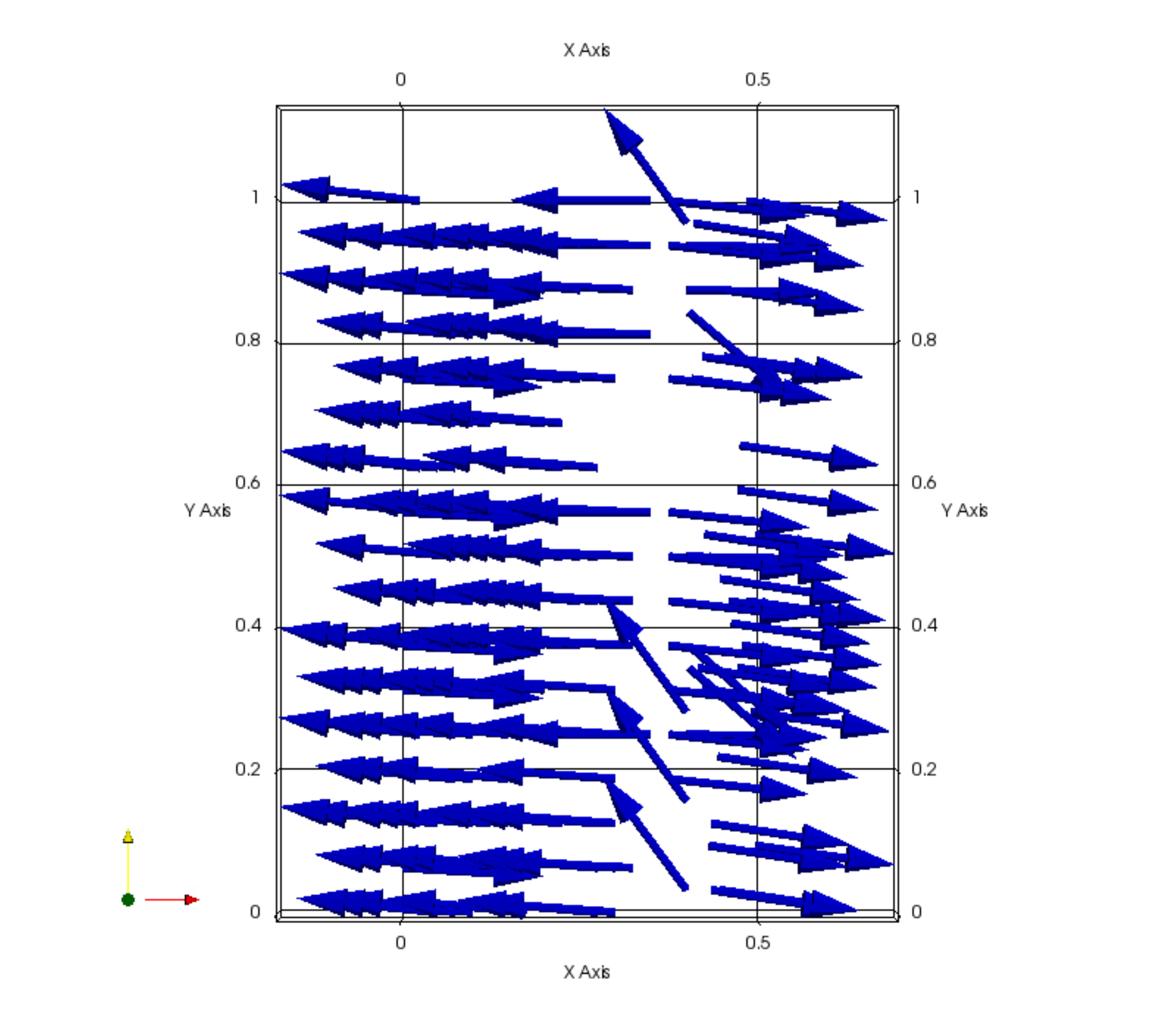} }
\subfloat[$\textbf{J}_{r}$]{
\includegraphics[scale=0.2]{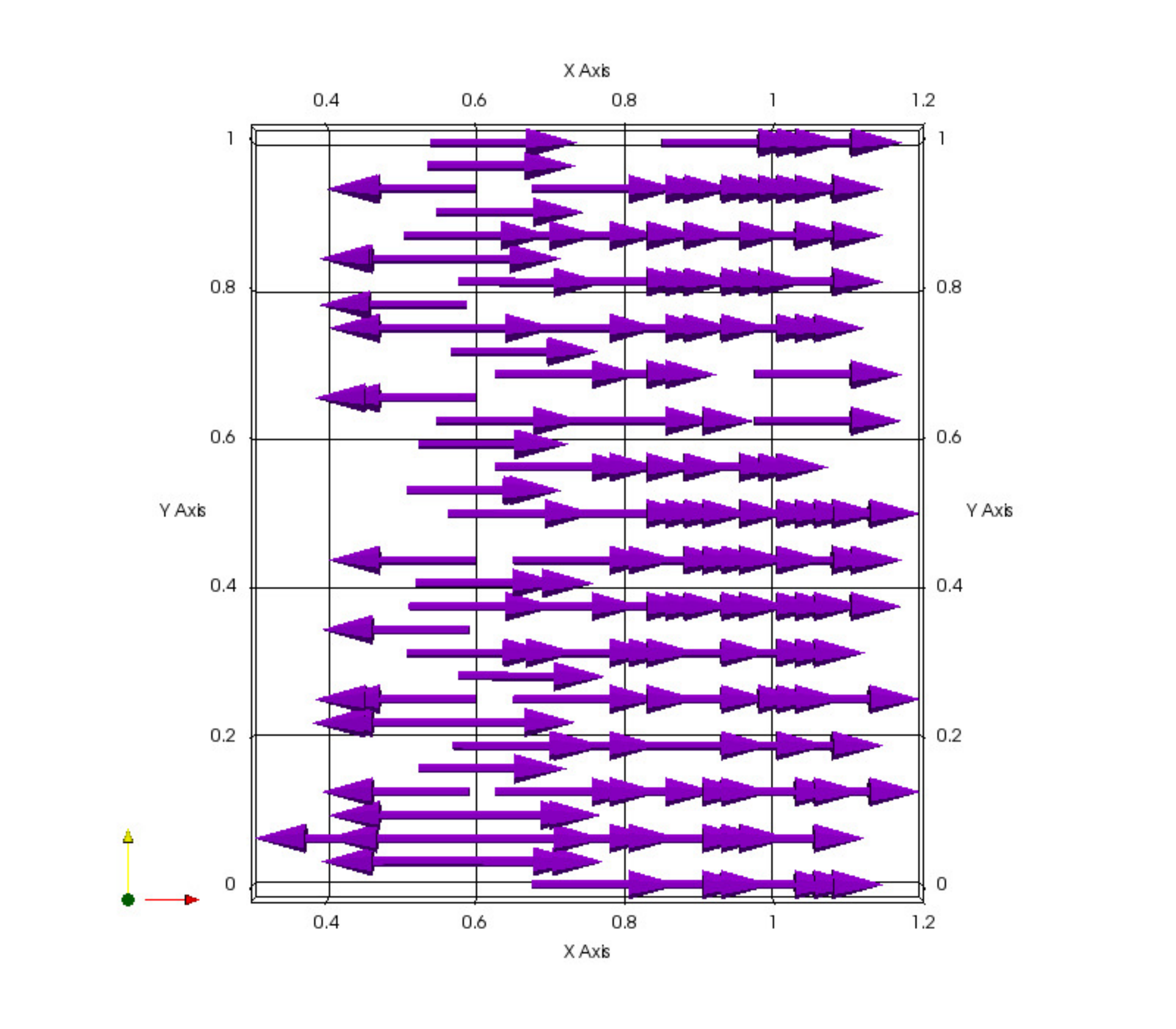} } \\
\subfloat[$\textbf{J}_{p}$]{
\includegraphics[scale=0.2]{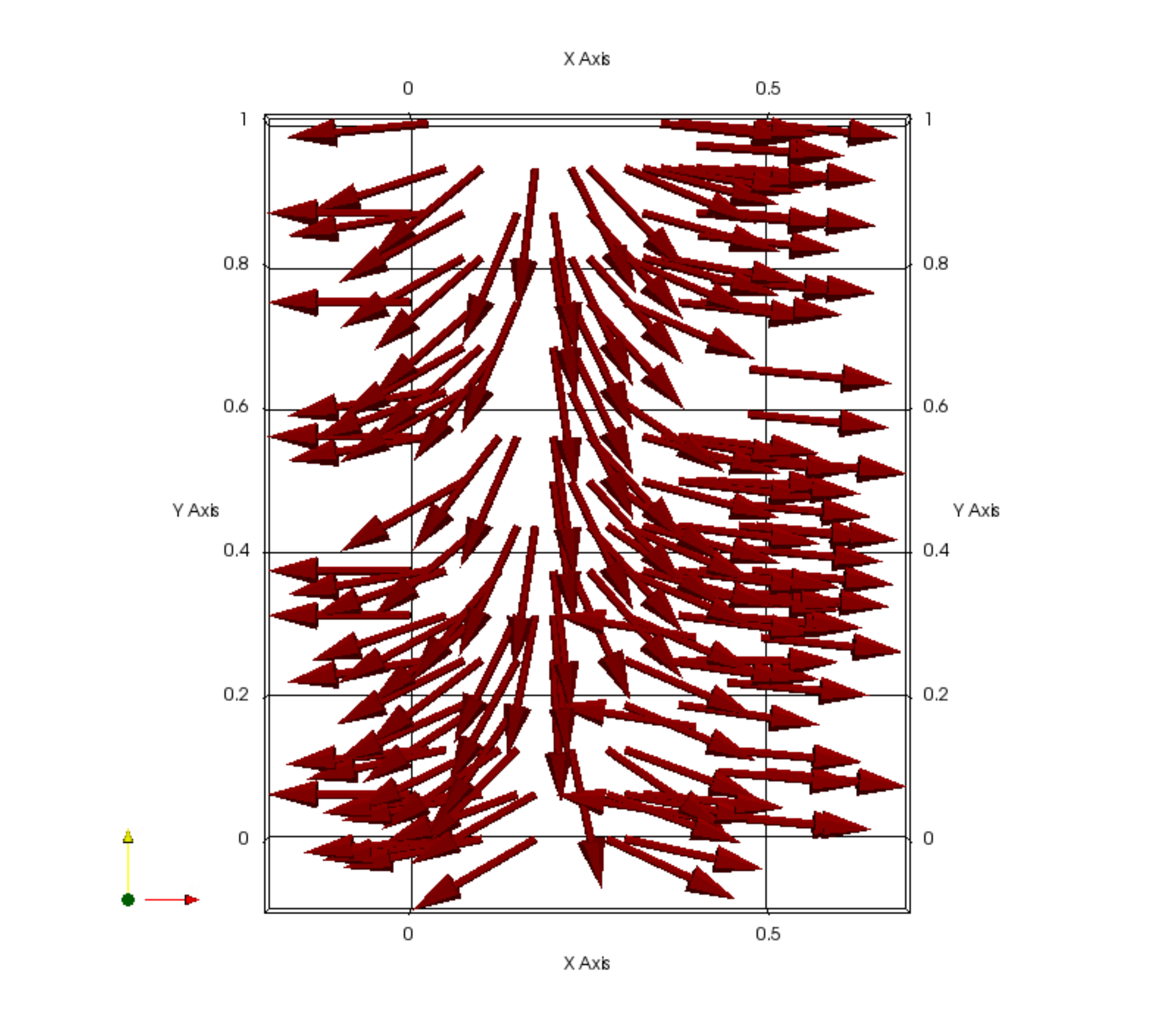} } 
\subfloat[$\textbf{J}_{o}$]{
\includegraphics[scale=0.2]{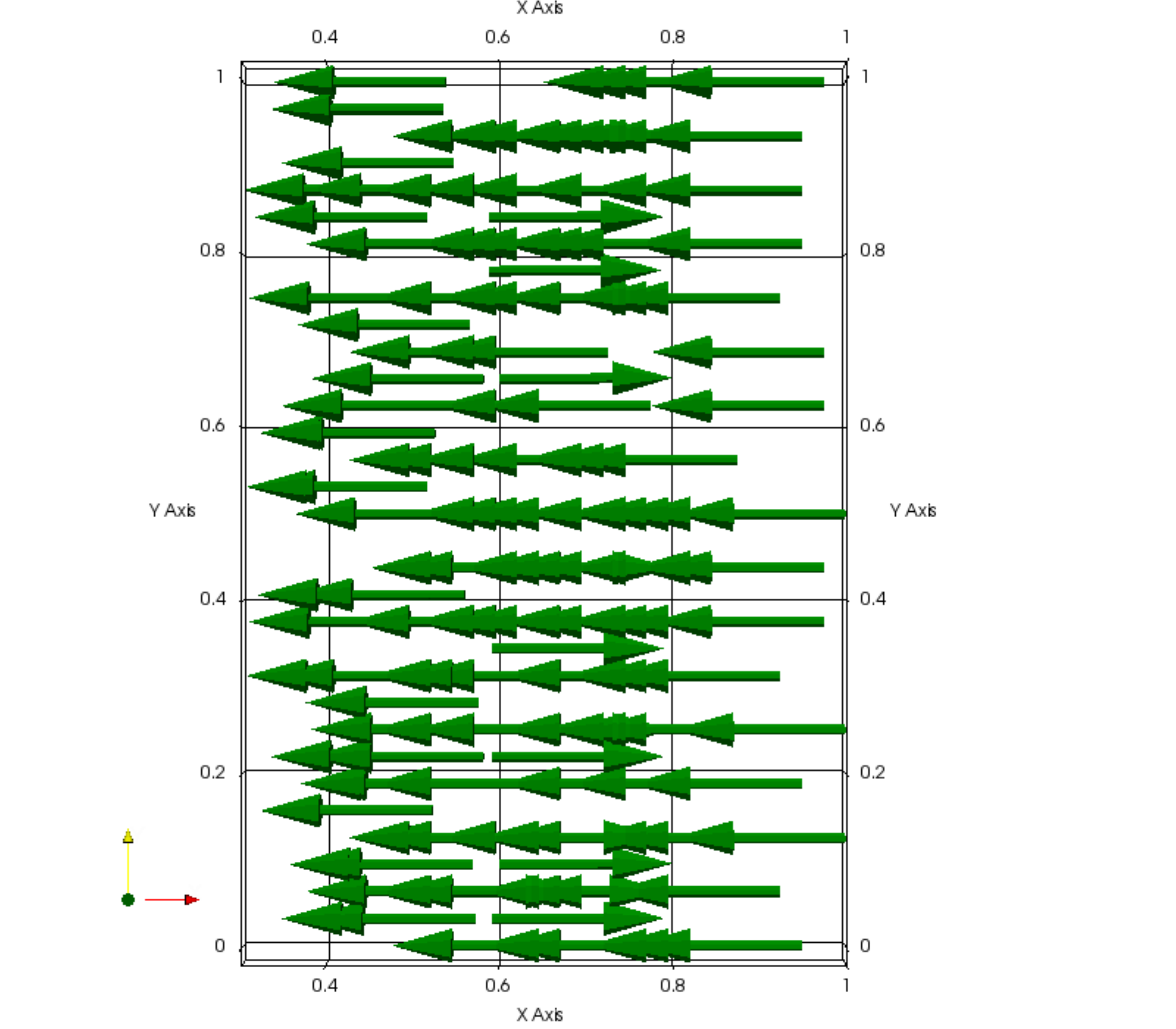} } 
\caption{Carrier currrents in device D-VIII with $\Phi_{\text{app}} \, = \, 0.0 $ for D-VIII.}
\label{fig:2D_flat_current_0}
\end{figure}

\begin{figure}[!ht]
\centering
\subfloat[$\textbf{J}_{n}$]{
\includegraphics[scale=0.2]{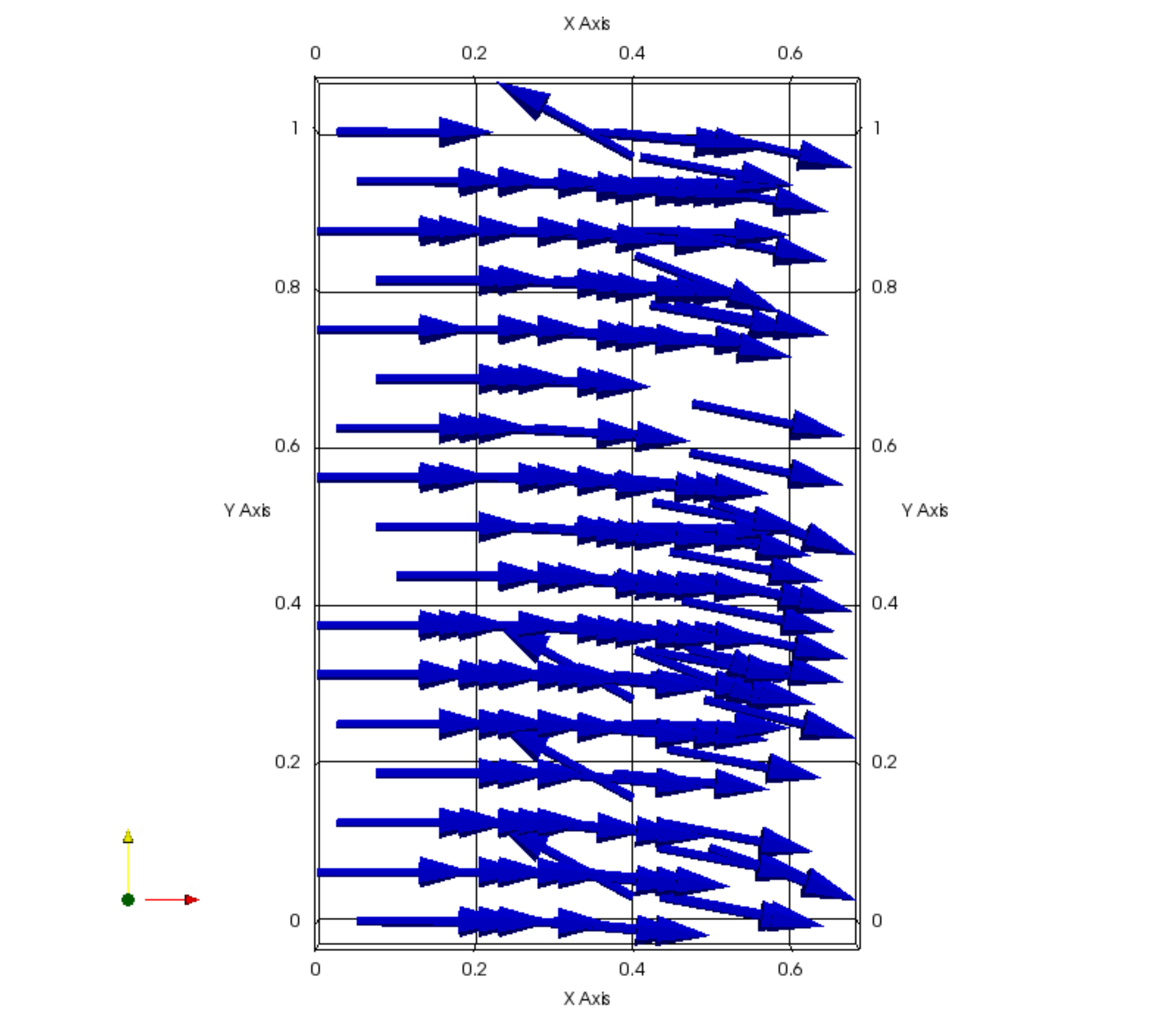} }
\subfloat[$\textbf{J}_{r}$]{
\includegraphics[scale=0.2]{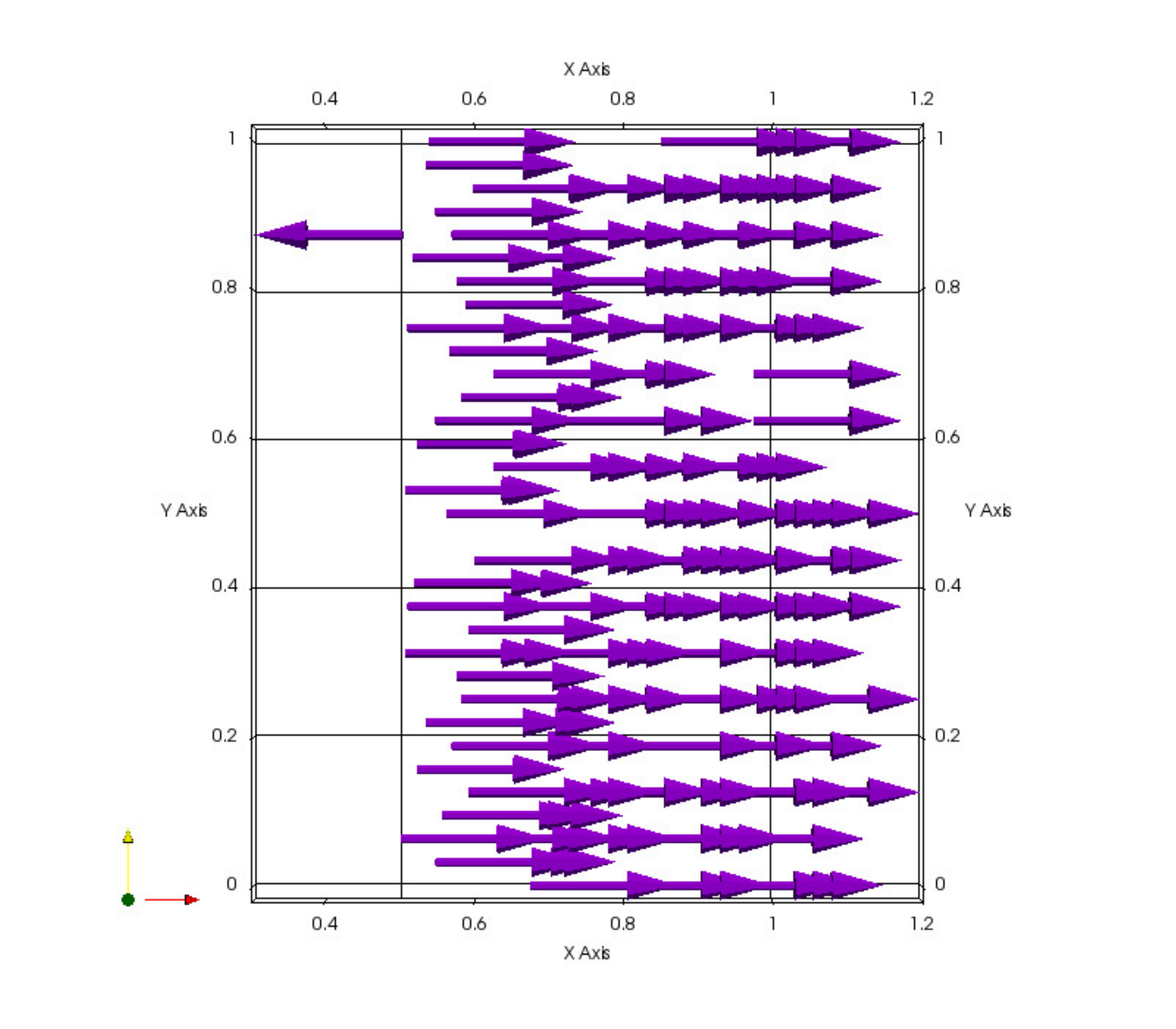} } \\
\subfloat[$\textbf{J}_{p}$]{
\includegraphics[scale=0.2]{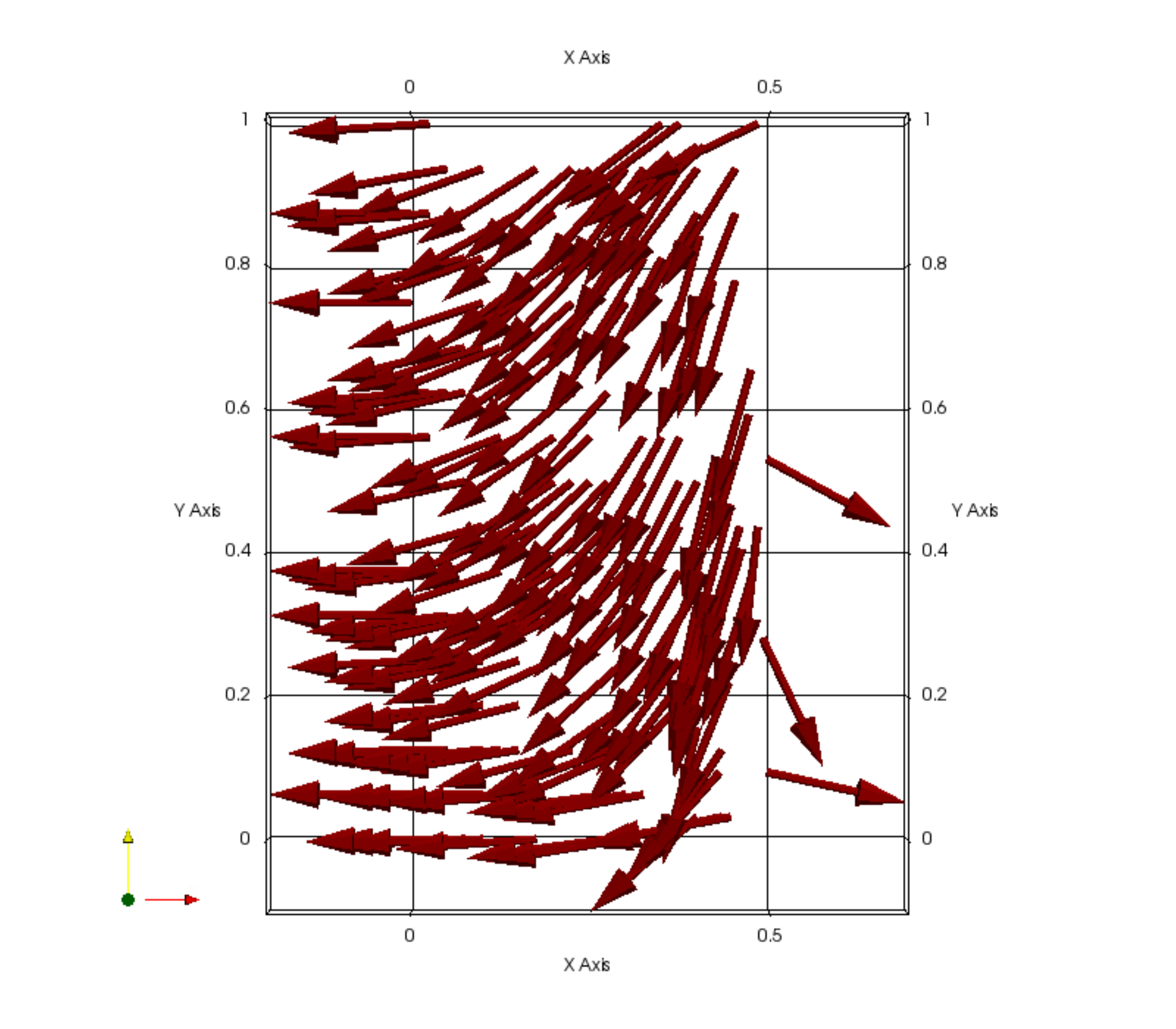} } 
\subfloat[$\textbf{J}_{o}$]{
\includegraphics[scale=0.2]{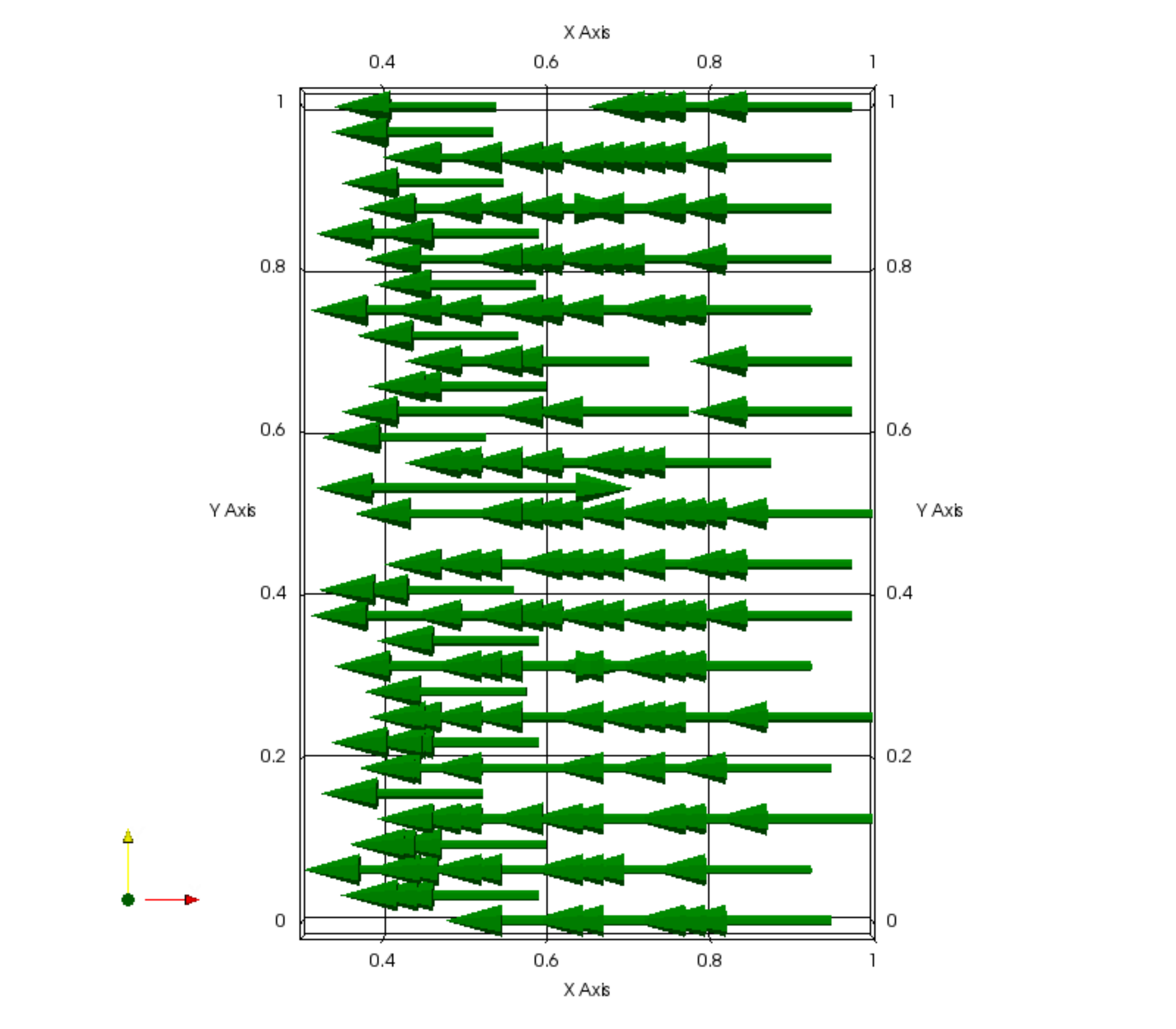} } 
\caption{Carrier currrents in device D-VIII with $\Phi_{\text{app}} \, = \, 0.5 $.}
\label{fig:2D_flat_current_50}
\end{figure}

The densities of electron, holes, reductants and oxidants under an applied bias $\Phi_{\text{app.}} = 0.0 ]$  and an applied bias of $\Phi_{\text{app.}} = 0.5$ are displayed in Fig.~\ref{fig:2D_flat_densities}.  We can see that under zero applied bias electrons are mostly forced away from the interface and under a large applied bias they are forced toward the interface.  This can also be verified from Fig.~\ref{fig:2D_flat_current_0} (a) and Fig.~\ref{fig:2D_flat_current_50} (a) which show the electron current density, $\textbf{J}_{n}$, for $\Phi_{\text{app.}} = 0.0 $ and $\Phi_{\text{app.}} = 0.5$ respectively.  From Fig.~\ref{fig:2D_flat_current_0} (b) and Fig.~\ref{fig:2D_flat_current_50} (b), we can see that the photo-generated holes are mostly forced towards the interface for the case of $\Phi_{\text{app.}} = 0.0 $ and mostly forced away from the interface for $\Phi_{\text{app.}} = 0.5 $. 

\begin{table}[!ht]
\centering
\caption{D-VIII (left) \& D-IX (right) domains.}
\vspace{1mm}
\begin{tabular}{cc|cc}
\hline
\hline
Parameter & Value & Parameter & Value  \\
\hline
R$_{1}$ & $0.5$ & R$_{1}$ & $0.3 $ \\
R$_{2}$ & $0.5 $ &  R$_{2}$ & $0.7 $ \\
H & $1 $ & H & $1 $ \\
L & $1 $ & L & $1 $ \\
\hline
\end{tabular}
\label{tab:2D_Flat}
\end{table}

Similarly, we can observe from Fig.~\ref{fig:2D_flat_densities}, Fig.~\ref{fig:2D_flat_current_0} and Fig.~\ref{fig:2D_flat_current_50} that reductants and oxidants are forced towards or away from interface depending on the value of the applied bias.  It is important to note that in Fig.~\ref{fig:2D_flat_densities} (b) the reductant and oxidant densities at the interface deviate appreciably from their bulk values.  Indeed, the bulk density of the oxidants is lower than the bulk density of the reductants, however, at the interface the density of oxidants is higher than the density of reductants.  In such instances, the validity of the assumption used for the Schottky approximations that the density or redox species remains constants is clearly not true. 

\begin{figure}[!ht]
\subfloat[$\Phi_{\text{app}} \, = \, 0.0$]{
\includegraphics[scale=0.32]{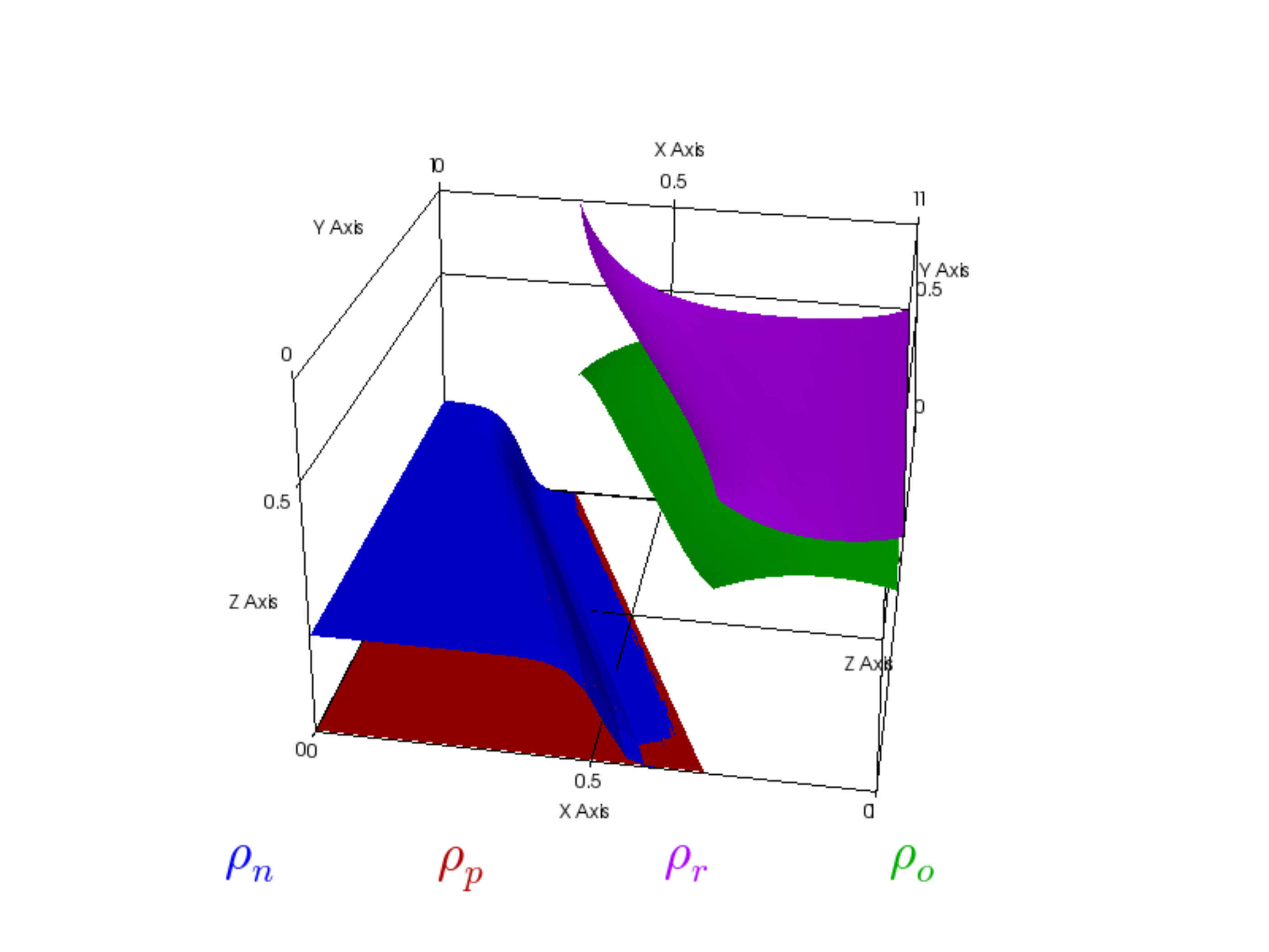} } 
\subfloat[$\Phi_{\text{app}} \, = \, 0.0 $]{
\includegraphics[scale=0.27]{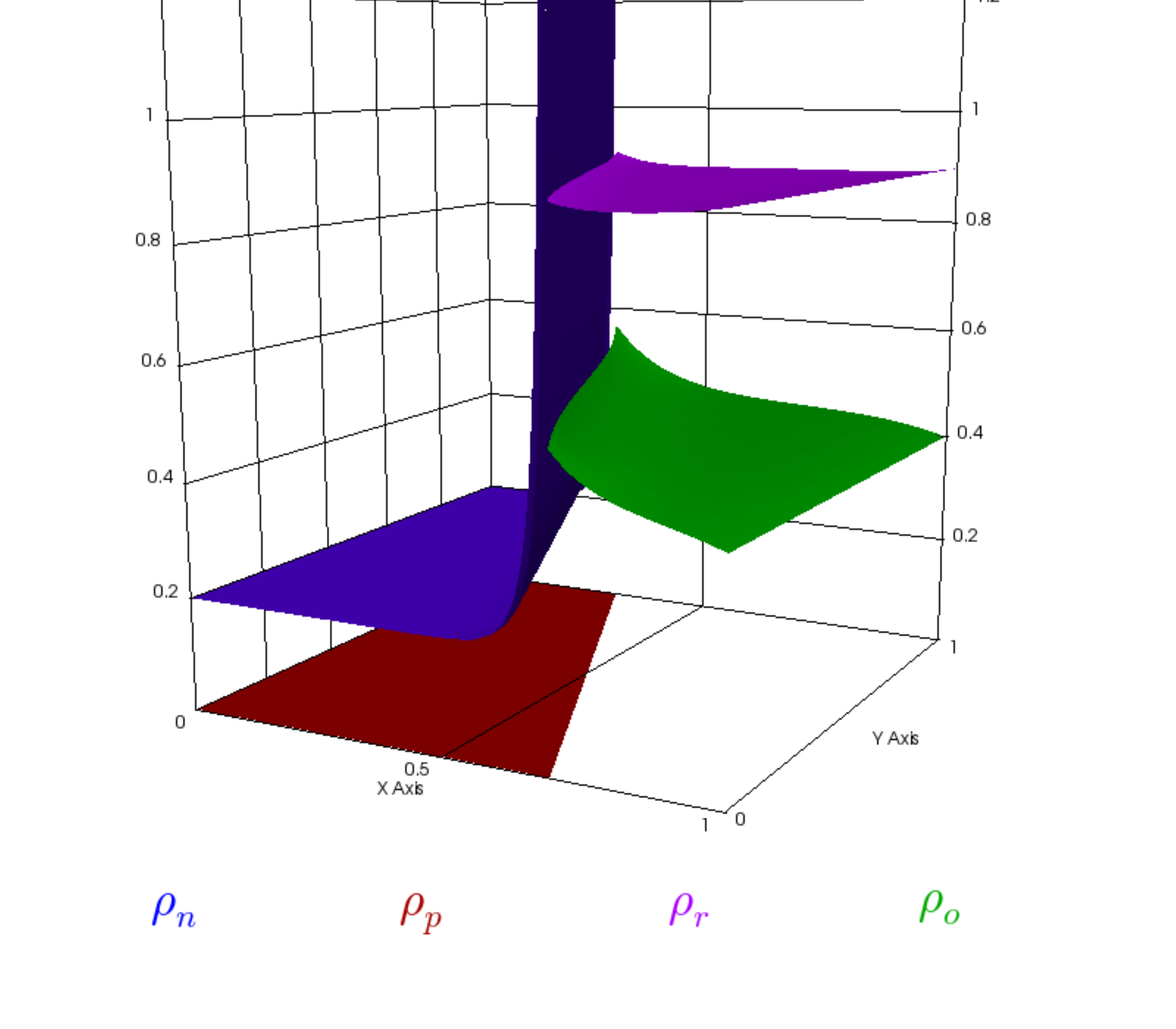}  }
\caption{Densities of charge carriers in device D-IX.  Note: Densities are in $10^{1} $.}
\label{fig:2D_slanted_densities}
\end{figure}

The densities of electrons, holes, reductants and oxidants of D-IX under illumination with an applied bias $\Phi_{\text{app.}} = 0.0$ and $\Phi_{\text{app.}} = 0.5 $ are displayed in Fig.~\ref{fig:2D_slanted_densities}.  We can immediately see the behavior of device D-IX is very different from the behavoir of device D-VIII.  The densities values of carriers in D-VIII were uniform y-axis, however, in device D-IX, the densities are not even uniform along the interface.  Indeed, in Fig.~\ref{fig:2D_slanted_densities} we can observe local increases and decreases in the densities of carriers; specifically at the corners ends of the interface, that is the points $(R_{1},H)$ and $(R_{2},0)$.

\begin{figure}[!ht]
\centering
\subfloat[$\textbf{J}_{n}$]{
\includegraphics[scale=0.2]{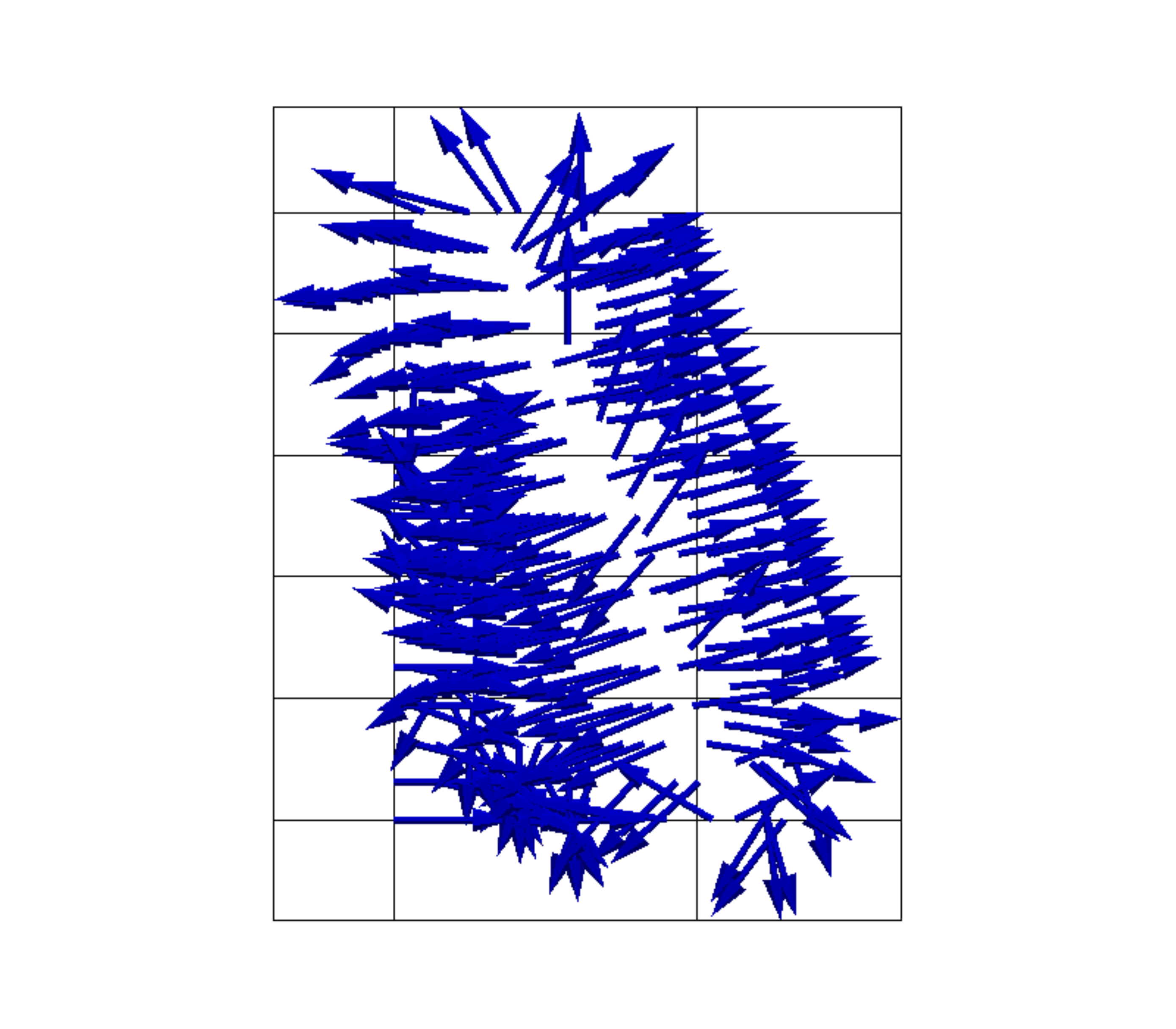} }
\subfloat[$\textbf{J}_{r}$]{
\includegraphics[scale=0.2]{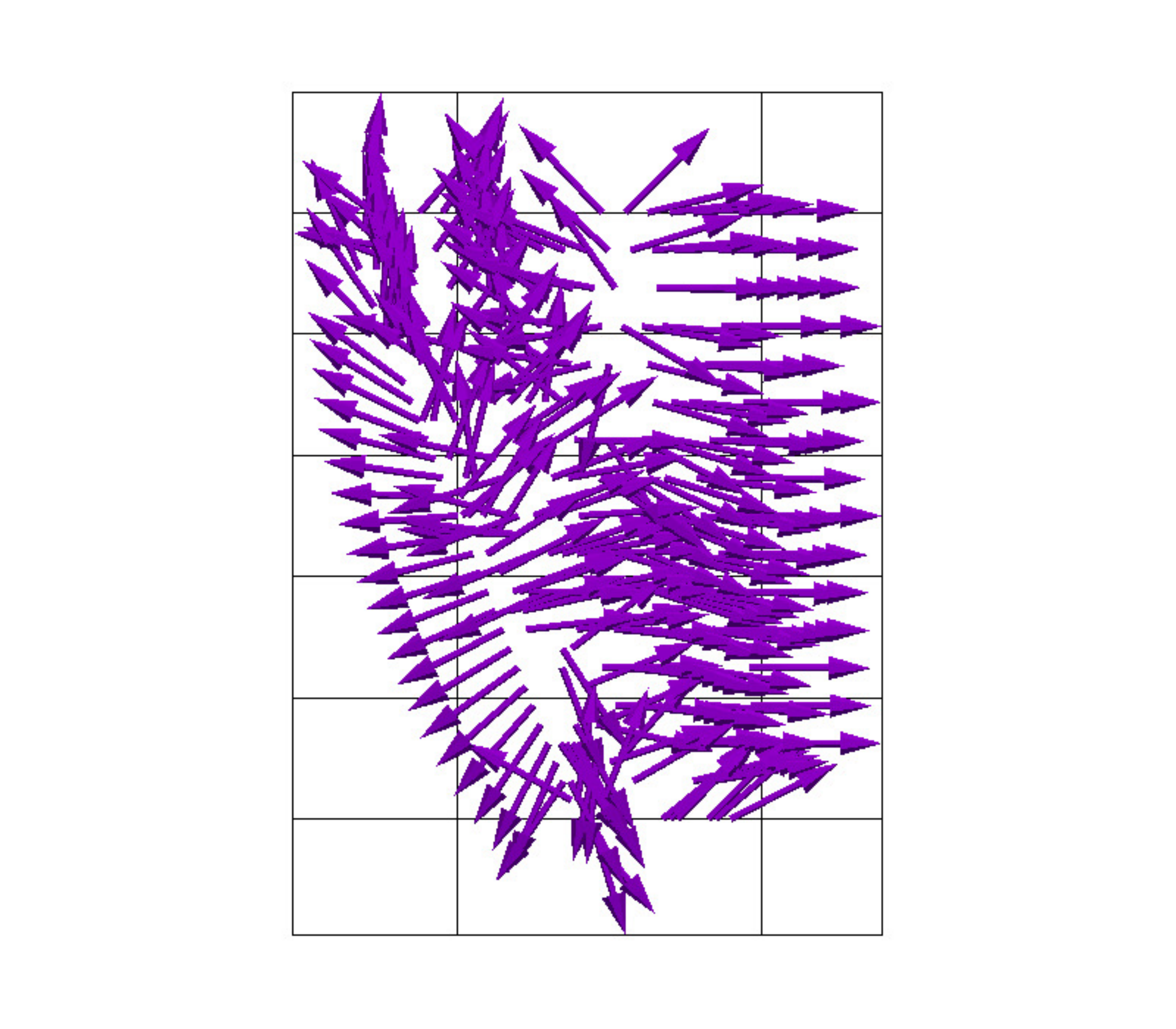} } \\
\subfloat[$\textbf{J}_{p}$]{
\includegraphics[scale=0.2]{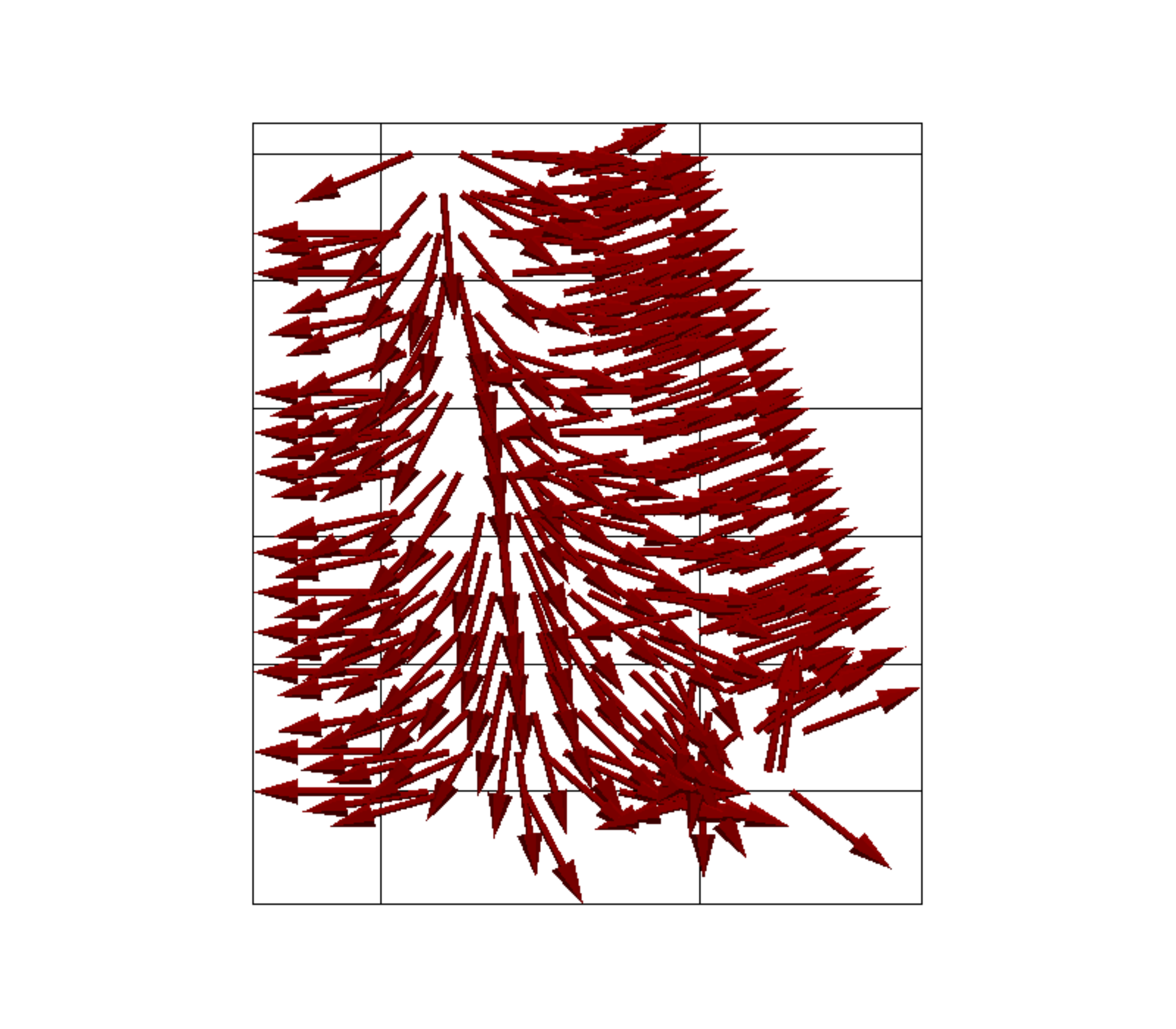} } 
\subfloat[$\textbf{J}_{o}$]{
\includegraphics[scale=0.2]{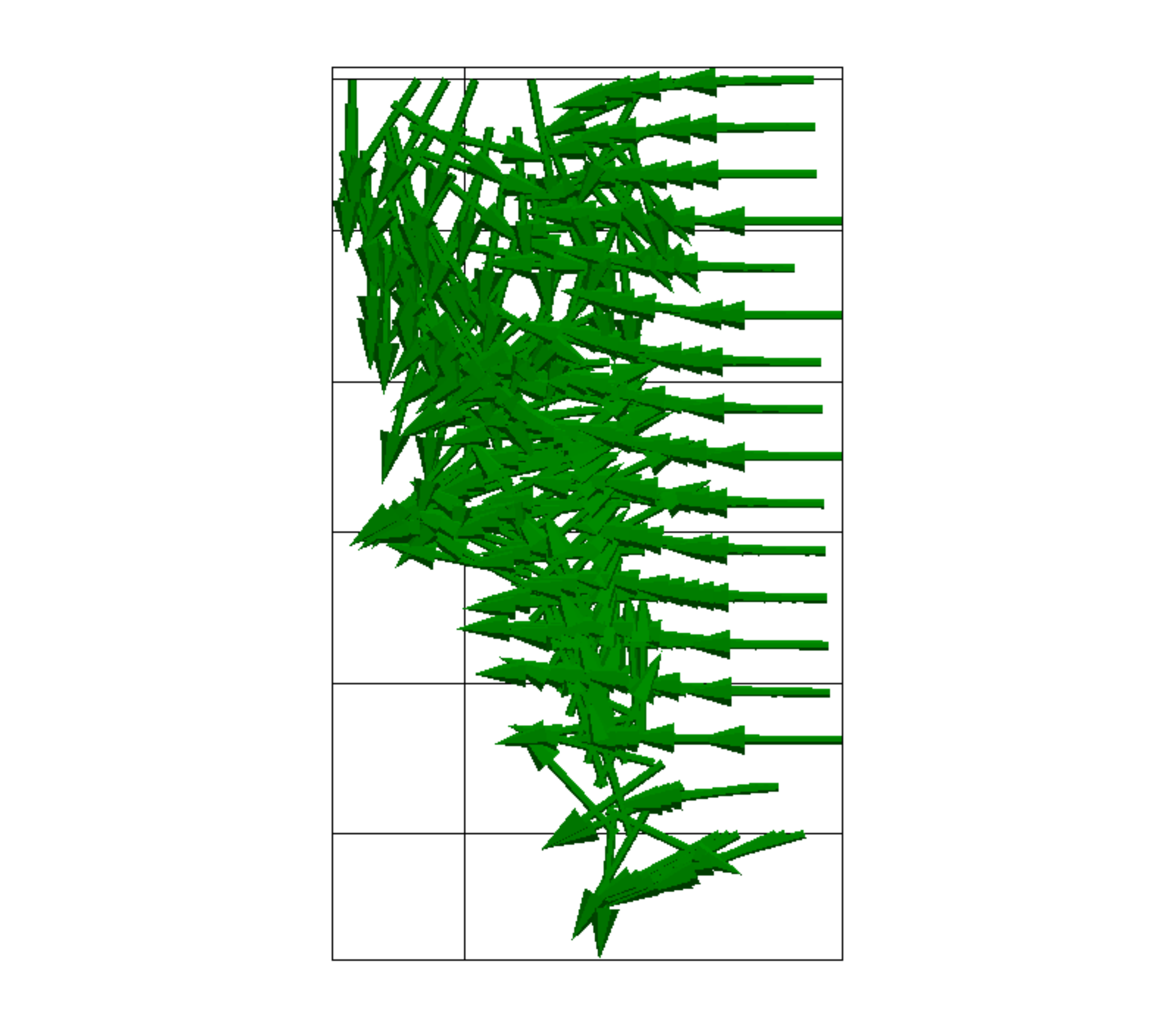} } 
\caption{Carrier currrents for device D-IX with $\Phi_{\text{app}} \, = \, 0$ .}
\label{fig:2D_slanted_current_0}
\end{figure}

\begin{figure}[!ht]
\centering
\subfloat[$\textbf{J}_{n}$]{
\includegraphics[scale=0.2]{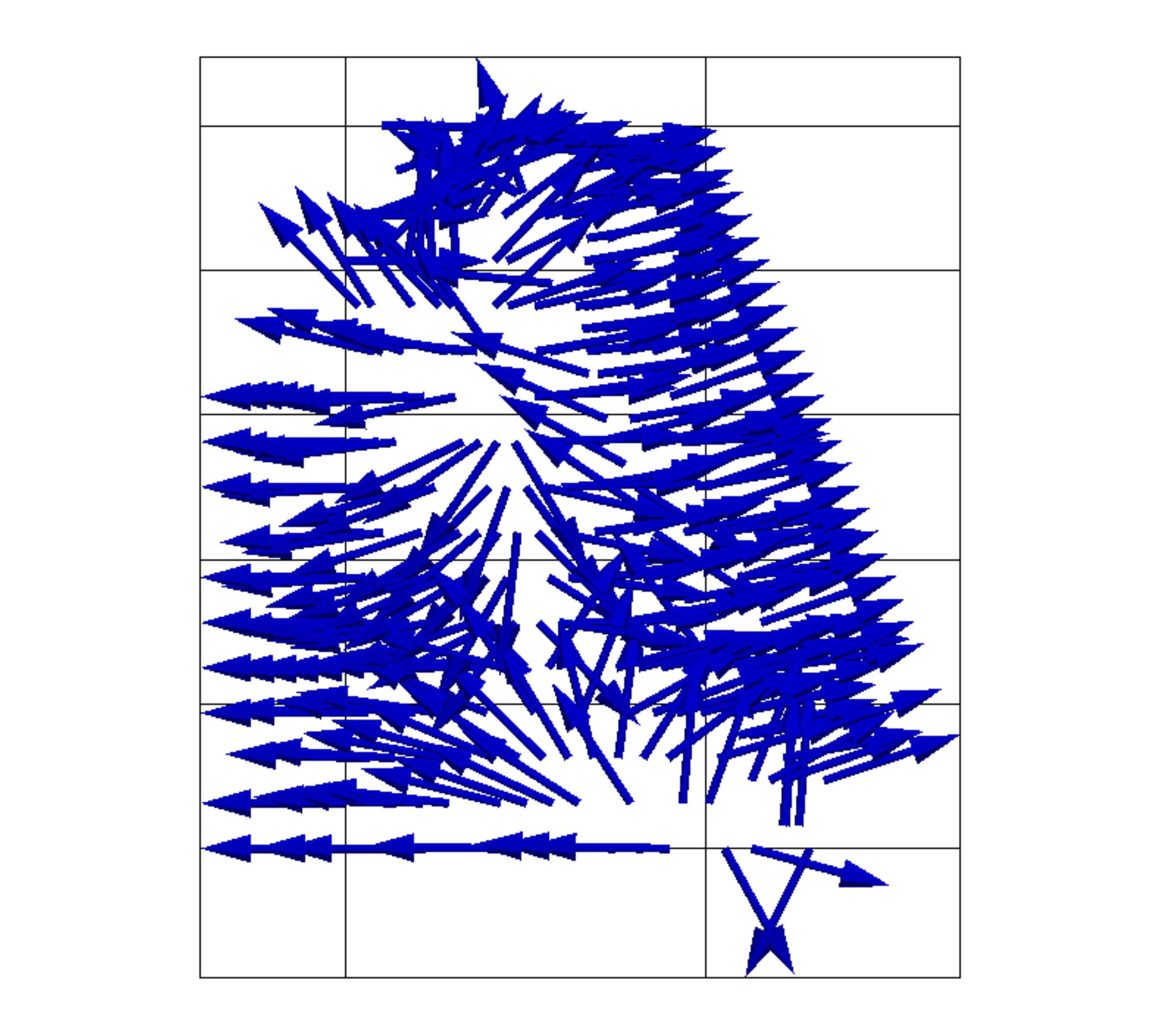} }
\subfloat[$\textbf{J}_{r}$]{
\includegraphics[scale=0.2]{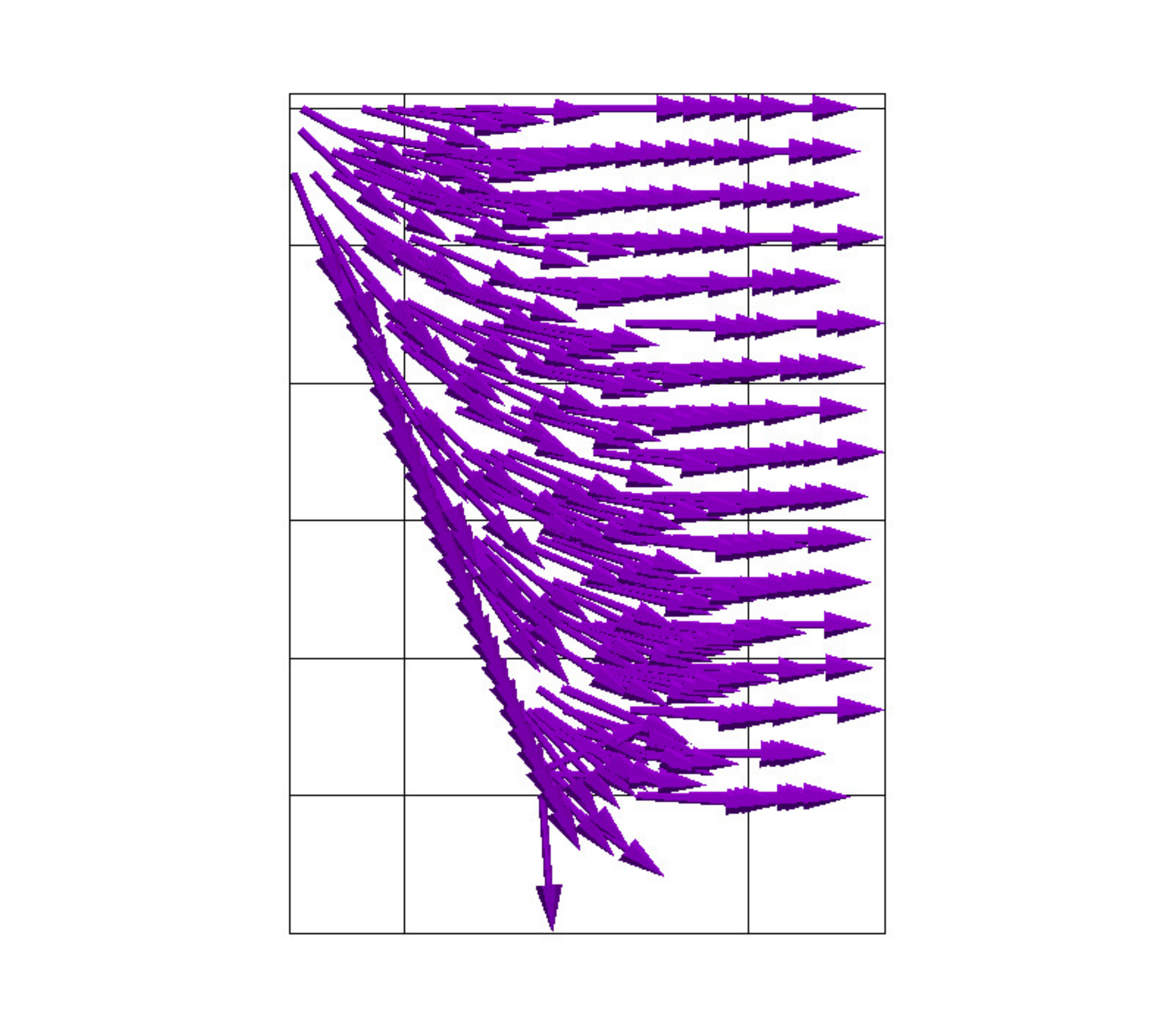} } \\
\subfloat[$\textbf{J}_{p}$]{
\includegraphics[scale=0.2]{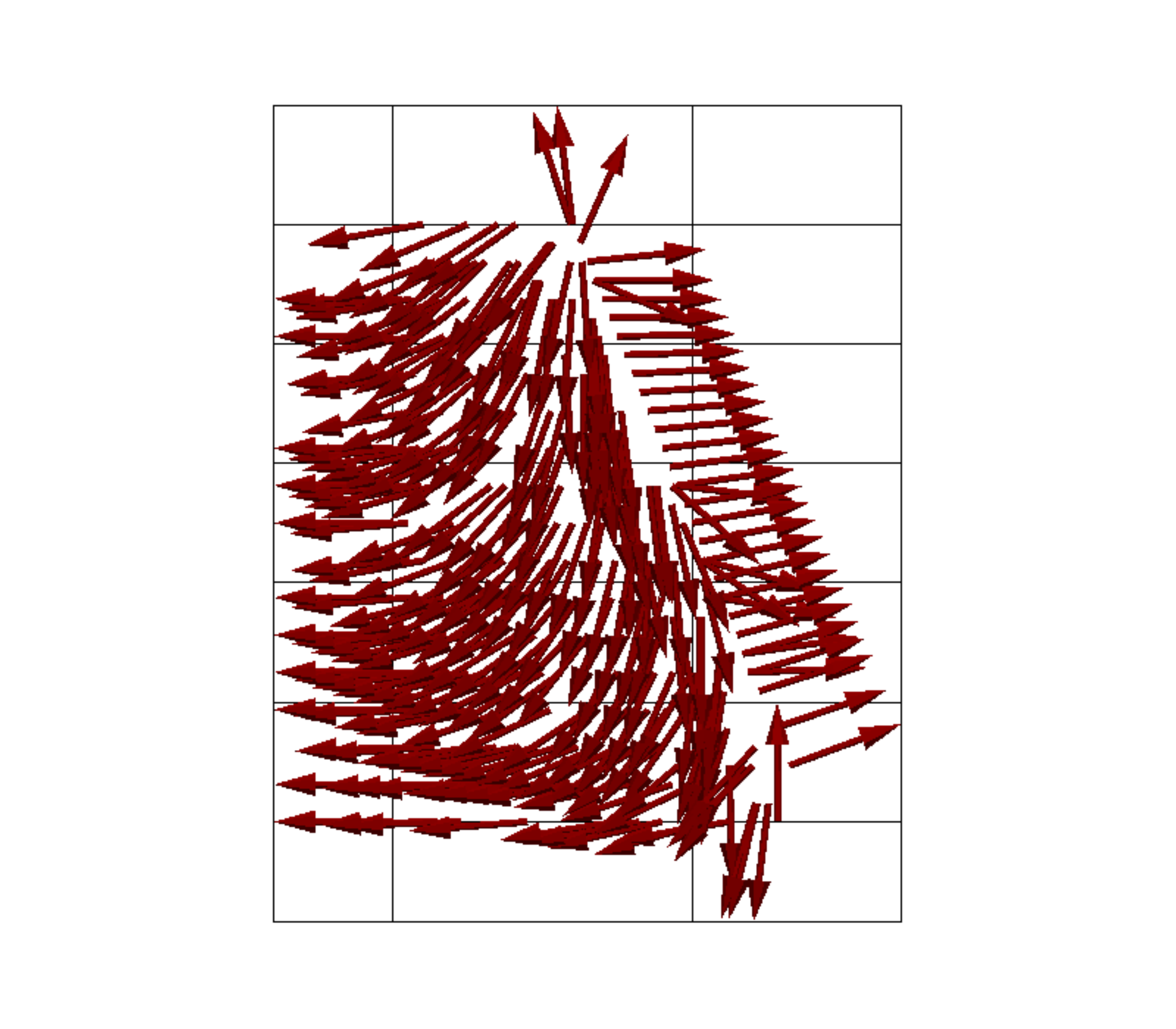} } 
\subfloat[$\textbf{J}_{o}$]{
\includegraphics[scale=0.2]{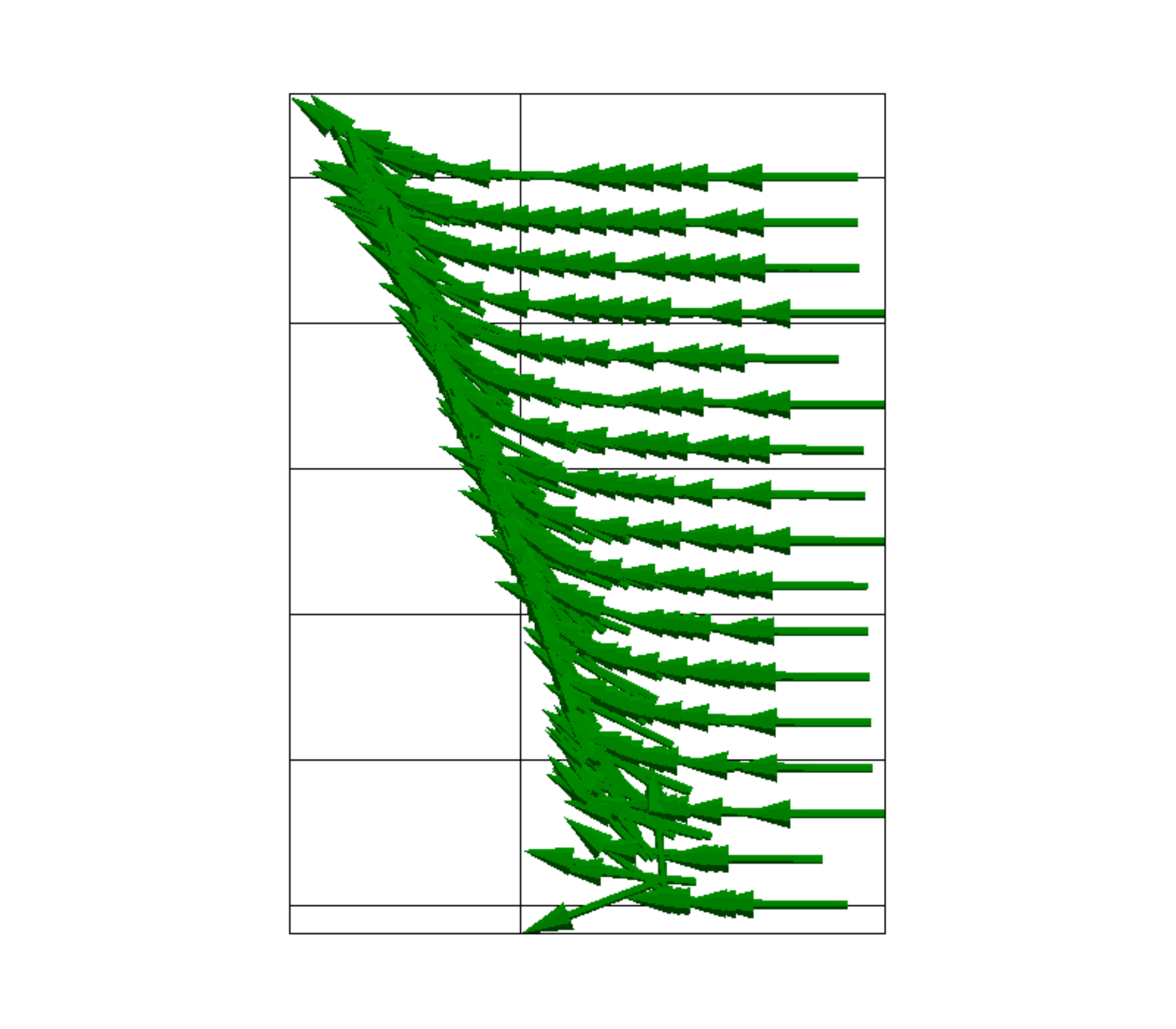} } 
\caption{Carrier currrents for device D-IX with $\Phi_{\text{app}} \, = \, 0.5 $.}
\label{fig:2D_slanted_current_50}
\end{figure}

\section{Conclusion and further remarks}
\label{SEC:Concl}

In this paper we proposed some numerical algorithms for solving systems of drift-diffusion-Poisson equations that arise from the modeling of charge transport in semiconductor-electrolyte based PEC solar cells. In our algorithm, we used a mixed finite element discretization for the Poisson equation and a local discontinuous Galerkin discretization for the corresponding drift-diffusion equations. For the temporal variable, we developed some implicit-explicit time stepping methods to alleviate the stability constraints imposed by the CFL conditions as well as to decouple the equations in the semiconductor domain from the equations in the electrolyte domain. We calibrated our implementation with manufactured solutions and presented numerical simulations to show the impact of various device parameters on the performance of the PEC cells.

The stability and convergence of some of the numerical algorithms we proposed can be proved in slightly simplified setting. Detailed stability and convergence analysis will appear in forthcoming work. On the application side, we are interested in the implementation of our algorithms on parallel processors for simulations in higher dimensional problems. Moreover, we are currently working on the usage of the model for optimal design of PEC cells that maximizing the performance of the cell over various achievable cell parameters.

Let us finish this paper by mentioning that even though our algorithms are designed for the simulation of PEC solar cells, they can be easily adapted to many other problems where charge transport through semiconductor-electrolyte interfaces have to be simulated. One such an application is in mathematical neuroscience where researchers are interested in studying performance of two-way communications between (sometimes organic) semiconductor chips and neural cells~\cite{AzGeBaDeCa-IEEE07,BaIvKuKu-Book09,PeFr-EPJE09,StMuFr-PRE97}. In this situation, a similar system of drift-diffusion-Poisson model can be derived for the charge transport process. Therefore, our numerical algorithms can be useful there.

\section*{Acknowledgments}

This work is partially supported by the National Science Foundation through grants DMS-1107465, DMS-1109625, and DMS-1321018. KR also acknowledges partial support from the Institute for Computational Engineering and Sciences (ICES) of UT Austin through a Moncrief Grand Challenge Faculty Award.


{\small

}


\end{document}